\documentclass[11pt]{amsart}

\usepackage[a4paper,margin=1in]{geometry}
\usepackage{amsmath,amssymb,amsfonts,amsthm,mathtools,mathrsfs,amscd}
\usepackage{tikz-cd}
\usepackage{cancel}
\usepackage{graphicx,multirow,booktabs,capt-of}
\usepackage{needspace}

\numberwithin{equation}{section}
\allowdisplaybreaks[1]
\newtheorem{theorem}{Theorem}[section]
\newtheorem{lemma}[theorem]{Lemma}
\newtheorem{proposition}[theorem]{Proposition}
\newtheorem{corollary}[theorem]{Corollary}
\theoremstyle{definition}
\newtheorem{definition}[theorem]{Definition}
\newtheorem{assumption}[theorem]{Assumption}

\theoremstyle{remark}
\newtheorem{remark}[theorem]{Remark}

\newcommand{\rd}{\mathrm d}
\newcommand{\epsi}{\epsilon}
\newcommand{\M}{\mathcal M}
\newcommand{\G}{\mathcal G}
\newcommand{\Ih}[1]{\mathcal{I}^{#1}_{r,h}}
\newcommand{\Th}{\mathcal T_h}
\newcommand{\Mh}{\mathcal M_h}
\newcommand{\Harm}{\mathfrak{H}}
\newcommand{\Pcal}{\mathcal P}
\newcommand{\Kcal}{\mathcal K}
\newcommand{\Vcal}{\mathcal V}
\newcommand{\Ncal}{\mathcal N}
\newcommand{\trc}{\operatorname{tr}}
\newcommand{\ntr}{\operatorname n}
\newcommand{\diam}{\operatorname{diam}}
\newcommand{\ip}[2]{\left\langle #1,#2\right\rangle}
\newcommand{\tracepair}[2]{%
  \mathopen{\langle\mkern-1.5mu\langle}#1,#2\mathclose{\rangle\mkern-1.5mu\rangle}}
\newcommand{\norm}[1]{\left\|#1\right\|}
\newcommand{\seminorm}[1]{\left|#1\right|}
\newcommand{\ringV}{\mathring{\Vcal}}
\newcommand{\ringN}{\mathring{\Ncal}}
\newcommand{\Proj}[3]{\Pi_{#1,#2}^{#3}}
\title{A Virtual Element Exterior Calculus framework}
\author[Xiaojing Dong, Yibing Han and Yunqing Huang]{Xiaojing Dong$^1$, Yibing Han$^{1\dagger}$ and Yunqing Huang$^1$}
\subjclass[2020]{Primary 65N30; Secondary 65N15}
\keywords{virtual element method; differential forms; de Rham complex;
polytopal meshes; norm equivalence}
\thanks{The research was supported by the National Natural Science Foundation of China
(Nos.~12671437, 12571393, 12071404, 12261131501), National Natural Science Foundation of China
Key Program (No.~12431014), Young Elite Scientist Sponsorship Program by CAST
(No.~2020QNRC001), the Science and Technology Innovation Program of Hunan Province
(No.~2024RC3158), Postgraduate Scientific Research Innovation Project of Hunan
Province (No.~LXBZZ2024112), the Project of Scientific Research Fund of the Hunan
Provincial Science and Technology Department (Nos.~2023GK2029, 2024ZL5017), and
Program for Science and Technology Innovative Research Team in Higher Educational
Institutions of Hunan Province of China.}
\thanks{$^\dagger$ Corresponding author.}
\address{$^1$Hunan Key Laboratory for Computation and Simulation in Science and Engineering,
Key Laboratory of Intelligent Computing \& Information Processing of Ministry of Education,
School of Mathematics and Computational Science, Xiangtan University,
Er huan Road, Xiangtan, 411105, Hunan, P.R. China}
\email[Xiaojing Dong]{dongxiaojing99@xtu.edu.cn}
\email[Yibing Han (corresponding author)]{202331510114@smail.xtu.edu.cn}
\email[Yunqing Huang]{huangyq@xtu.edu.cn}
\date{}

\begin{document}
\begin{abstract}
We develop a unified virtual element framework for compatible polytopal
discretizations of differential forms in arbitrary dimension and at arbitrary order.
On general polytopal meshes, a central challenge is to construct
de Rham subcomplexes and establish stability estimates with mesh-independent constants.
We construct a trimmed complex of conforming virtual differential forms recursively
from local Hodge--Laplace boundary value problems and obtain a complete subcomplex
by imposing additional moment constraints.
The local spaces of these two complexes contain the trimmed and
complete polynomial differential forms, respectively, and admit unisolvent
degrees of freedom from which polynomial projections can be computed.
The canonical interpolants commute with the exterior derivative, and each virtual
complex has cohomology isomorphic to that of the continuous de Rham complex.
For both virtual complexes and in every form degree, we establish negative-norm
inverse estimates and thereby prove a mesh-uniform equivalence between the $L^2$ norm
and the scaled degree-of-freedom norm.
Through the cochain isomorphisms induced by the degrees of freedom, these norm
equivalences yield uniformly bounded liftings of fully discrete data into the
corresponding continuous spaces, in particular for the Discrete de Rham method.
We also derive canonical interpolation error estimates and discrete Poincar\'e
inequalities. As an
application, we obtain an inf-sup stable mixed Hodge–Laplace method with harmonic variables and optimal-order convergence under sufficient regularity. 
\end{abstract}
\maketitle

\section{Introduction}
\label{sec:introduction}

In this work, we construct two families of arbitrary-order virtual de Rham complexes of
differential forms and establish mesh-uniform norm equivalences that bridge the
virtual-function and fully discrete approaches to the analysis of polytopal methods.
Let $\Omega\subset\mathbb R^n$ be a bounded Lipschitz polyhedral domain
with $n\ge2$.
The continuous de Rham complex on $\Omega$ is given by
\begin{equation}
\label{eq:domain-derham-complex}
  0
  \longrightarrow
  H\Lambda^0(\rd;\Omega)
  \xrightarrow{\rd^0}
  H\Lambda^1(\rd;\Omega)
  \xrightarrow{\rd^1}\cdots
  \xrightarrow{\rd^{n-1}}
  H\Lambda^n(\rd;\Omega)
  \longrightarrow 0.
\end{equation}
Here, $\rd^k$ denotes the exterior derivative on $k$-forms, and
$H\Lambda^k(\rd;\Omega)$ consists of $L^2$ differential $k$-forms with
$L^2$ exterior derivatives.

Differential complexes play a central role in the analysis of partial
differential equations (PDEs) and the construction of stable, structure-preserving
discretizations~\cite{ArnoldFalkWinther2010,ArnoldHu2021}.
The de Rham complex~\eqref{eq:domain-derham-complex} is a canonical
example: the analysis of well-posedness for the associated PDEs involves
its cohomology spaces $\ker\rd^k/\operatorname{im}\rd^{k-1}$, whose
classes have unique harmonic representatives.
This motivates \emph{compatible} discretizations that retain the complex
structure and reproduce its cohomology.
Finite element exterior calculus
(FEEC)~\cite{Arnold2006,ArnoldFalkWinther2010,Arnold2018} systematically
develops compatible discretizations through the construction and analysis
of conforming finite element subcomplexes of~\eqref{eq:domain-derham-complex}
on simplicial and cubical meshes.

In recent years, polytopal methods, including mimetic finite
differences~\cite{BrezziBuffaLipnikov2009,BeiraoLipnikovManzini2011},
finite volume methods~\cite{DroniouEymardGallouetHerbin2010},
virtual element methods (VEM)~\cite{BeiraoEtAl2013},
hybrid high-order methods~\cite{DiPietroDroniou2020HHO},
and discrete de Rham (DDR)
methods~\cite{DiPietroDroniou2023,DiPietroDroniouRapetti2020},
have attracted considerable attention owing to their geometric flexibility.
VEM is a representative approach to compatible discretization on such meshes.
It defines local spaces through boundary value problems and computes
suitable polynomial projections from the degrees of freedom (DOFs) without
explicit knowledge of the virtual basis
functions~\cite{BeiraoEtAl2014,BeiraoEtAl2016}.
Conforming VEM subcomplexes have been constructed for the
three-dimensional vector proxies of the de Rham
complex~\cite{BeiraoEtAl2016,BeiraoEtAl2018}.
Such complexes and their variants with enhanced regularity have
found applications in electromagnetism, incompressible flows,
magnetohydrodynamics, and fourth-order
problems~\cite{BeiraoEtAl2018,ChenGuoWei2025Maxwell,
BeiraoDassiVacca2020,MoraSilgado2025,
BeiraoDassiManziniMascotto2023MHD,DongHanHuang2026QuadDiv}.

Another route to compatible discretization is
DDR, which reconstructs
differential operators and polynomial potentials from polynomial data on
mesh entities. In three dimensions, VEM and DDR may be viewed as twin approaches:
higher-order potential reconstruction in DDR and VEM enhancement~\cite{AhmadEtAl2013}
play analogous roles in reducing the number of DOFs~\cite{DiPietroDroniou2023Serendipity},
and the DOF maps induce a cochain isomorphism between the virtual and fully
discrete complexes of each method~\cite{BeiraoDassiDiPietroDroniou2022}.
Although DDR operates in a fully discrete setting, conforming liftings
into continuous spaces have been used to establish adjoint-consistency
estimates needed in convergence
analysis~\cite{DiPietroDroniou2023,DiPietroDroniouPitassi2025}.
In conforming VEM, the virtual functions themselves provide continuous
representatives, while mesh-uniform norm equivalences play a central role in
stability analysis (see, e.g., the work on nodal, edge, and face elements
in~\cite{BeiraoMascotto2023,BeiraoMascottoMeng2022,ChenHuang2018}).

For differential forms in arbitrary dimension, two families of fully discrete
de Rham complexes inspired respectively by VEM and DDR have been
constructed~\cite{Bonaldi2025}.
The DDR-inspired complex also admits a conforming interpretation within
the finite element systems (FES) framework~\cite{ChristiansenRapetti2025}.
However, the uniformity of the norm-equivalence constants in this FES
interpretation remains conditional on a conjectural stability
hypothesis~\cite[Remark~2.1 and Section~5]{ChristiansenRapetti2025}.
For the local VEM spaces underlying our two conforming complexes, we prove that,
in every form degree, the $L^2$ norm of a virtual form is equivalent to the
scaled norm of its degrees of freedom, with constants controlled by explicit
mesh regularity parameters.
Together with the degree-of-freedom isomorphisms, these equivalences ensure that
stability and approximation estimates in the corresponding norms carry over
between the two formulations.
This link between the virtual-function and fully discrete approaches, highlighted
by Beir{\~a}o da Veiga et al.~\cite{BeiraoDiPietroDroniou2026}, provides a
common basis for the analysis of compatible polytopal methods.

Our construction starts with trimmed local spaces obtained by extending
the $\mathcal B^*$ problem in FEEC~\cite[Section~4.4.3]{Arnold2018}
to prescribed polynomial interior data and boundary traces assembled
from lower-dimen\-sional spaces.
For each polynomial degree $r\ge0$, these spaces contain the trimmed
polynomial differential forms of order $r+1$.
We extend the subspaces with vanishing tangential traces from mesh
subcells to obtain a geometric decomposition of the local spaces,
analogous to that in~\cite[Section~4.7]{Arnold2006}, and thereby
derive explicit dimension formulas.
Following VEM enhancement~\cite{AhmadEtAl2013}, we impose moment
constraints on forms and their exterior derivatives to define complete
subspaces containing all polynomial forms up to degree $r$.
The corresponding degree-$r$ $L^2$ projections remain computable from the
reduced DOFs.
Both families admit unisolvent moment DOFs and form conforming
complexes with commuting interpolants and cohomology isomorphic to
that of \eqref{eq:domain-derham-complex}.
In particular, at $r=0$ the complete complex coincides with the Whitney
complex on simplicial meshes~\cite{Arnold2006}.
The DOF maps of both families also induce cochain isomorphisms with the corresponding
fully discrete complexes in~\cite{Bonaldi2025}.

For the stability and approximation analysis, we assume that each
positive-dimensional mesh cell is uniformly star-shaped and has
diameter uniformly equivalent to that of any containing element.
We establish negative-norm inverse estimates by
induction using local Hodge decompositions, polynomial inverse and
Poincar\'e inequalities, and multiplicative Rellich-type trace estimates.
Recursive application of these inverse estimates and local upper bounds
yields the norm equivalences for the trimmed family, which
extend to the complete family through its defining moment constraints.
Canonical interpolation error bounds follow from scaled Sobolev trace
estimates and polynomial approximation.

To discretize the Hodge--Laplace problem, we equip the complete complex
with computable, polynomially consistent discrete $L^2$ inner products
whose stability follows from the norm equivalence.
We combine a cellular cochain estimate~\cite[Lemma~6]{pietro_uniform_2025}
with the local norm equivalence to prove the discrete Poincar\'e inequality
for the trimmed complex. A stable cochain retraction then transfers this
inequality to the complete complex.
The discrete Hodge decomposition and Poincar\'e inequalities yield a
uniformly inf--sup stable mixed method that retains harmonic variables
on domains with nontrivial cohomology.
Optimal-order error estimates follow under sufficient regularity of the
solution, source term, and harmonic forms.

Section~2 introduces the continuous and discrete settings.
Section~3 presents the construction, computability, commuting
interpolation, and cohomology of the two virtual complexes.
Section~4 establishes norm equivalences and
approximation estimates, and Section~5 analyzes the mixed Hodge--Laplace
method.
Appendix~A proves the scaled trace and Poincar\'e estimates stated in Section~4.
Appendix~B establishes the equivalence with fully discrete polytopal
methods and details the analytical transfer enabled by the norm equivalences.
Appendix~C proves the discrete Poincar\'e inequality for the trimmed
complex, used in the complete-family analysis of Section~5.

\setcounter{section}{1}

\section{Preliminaries}
\label{sec:preliminaries}
\subsection{Differential forms and Hilbert complexes}
\label{subsec:continuous-preliminaries}
We review smooth differential forms, the associated Hilbert spaces and
complexes, and their Hodge theory~\cite{Arnold2018}, together with weak
traces for boundary conditions~\cite{Schulz2022}.

\subsubsection{Smooth differential forms}
\label{subsec:exterior-calculus}
Let $\M\subset\overline\Omega$ be an oriented $m$-dimensional
Lipschitz domain relative to its affine hull, with $0\le m\le n$.
For a scalar function space $X(\M)$, such as $L^2(\M)$ or
$C_c^\infty(\M)$, let $X\Lambda^k(\M)$ denote the space of differential
$k$-forms whose coefficients in affine coordinates belong to $X(\M)$.
We write $\Lambda^k(\M)$ for smooth $k$-forms and understand all these
spaces to be zero for $k<0$ or $k>m$.
For $m=0$, point evaluation identifies $\Lambda^0(\M)$ with $\mathbb R$.

For $\omega\in\Lambda^k(\M)$ and $\zeta\in\Lambda^\ell(\M)$,
the wedge product satisfies
$\omega\wedge\zeta=(-1)^{k\ell}\zeta\wedge\omega$.
For a smooth vector field $Y$ on $\M$ and tangent vectors
$v_1,\ldots,v_{k-1}$ at $x\in\M$, contraction is defined by
\begin{equation}
\label{eq:form-contraction}
  (\omega\lrcorner Y)_x(v_1,\ldots,v_{k-1})
  :=\omega_x(Y(x),v_1,\ldots,v_{k-1}),
  \qquad k\ge1,
\end{equation}
with $\omega\lrcorner Y:=0$ for $k=0$, and satisfies
$(\omega\lrcorner Y)\lrcorner Y=0$.
The exterior derivative
$\rd^k:\Lambda^k(\M)\to\Lambda^{k+1}(\M)$
satisfies $\rd^{k+1}\rd^k=0$.
For $\omega\in\Lambda^k(\M)$ and
$\zeta\in\Lambda^{m-k-1}(\M)$ smooth up to the boundary,
Stokes' theorem gives the integration-by-parts formula
\begin{equation}
\label{eq:stokes-integration-by-parts}
  \int_{\M}\rd^k\omega\wedge\zeta
  +(-1)^k\int_{\M}\omega\wedge\rd^{m-k-1}\zeta
  =\int_{\partial\M}\trc_k\omega\wedge\trc_{m-k-1}\zeta,
\end{equation}
where $\trc_\ell:=\mathrm i^*$ denotes the tangential trace on $\ell$-forms,
$\mathrm i:\partial\M\hookrightarrow\overline\M$ is the inclusion,
and $\mathrm i^*$ denotes its pullback.
The boundary carries the induced orientation, and we set $\trc_m:=0$.

Let $\star_{\M}:\Lambda^k(\M)\to\Lambda^{m-k}(\M)$ denote the Hodge star
associated with the Euclidean metric and the chosen orientation.
It satisfies $\star_{\M}^2=(-1)^{k(m-k)}I$.
The formal codifferential
$\delta_k:\Lambda^k(\M)\to\Lambda^{k-1}(\M)$
is defined by
$\delta_k:=(-1)^k\star_{\M}^{-1}\rd^{m-k}\star_{\M}$
for $1\le k\le m$, with $\delta_0:=0$.
Using the boundary Hodge star $\star_{\partial\M}$ associated with
the induced metric and orientation, define the normal trace of
$\omega\in\Lambda^k(\M)$ smooth up to the boundary by
$\ntr_k\omega:=\star_{\partial\M}^{-1}\mathrm i^*(\star_{\M}\omega)$
for $1\le k\le m$, with $\ntr_0:=0$.
When unambiguous, we omit the domain subscript on $\star_{\M}$
and the degree indices on $\rd^k$, $\delta_k$, $\trc_k$, and $\ntr_k$.

\subsubsection{Hilbert spaces of differential forms}
\label{subsubsec:hilbert-form-spaces}
For $\omega,\zeta\in L^2\Lambda^k(\M)$, write
$\langle\omega,\zeta\rangle_{\M}
:=\int_{\M}\omega\wedge\star\zeta$
and $\norm{\omega}_{\M}^2:=\langle\omega,\omega\rangle_{\M}$
for the $L^2$ inner product and its induced norm.
We use analogous notation on $\partial\M$.
For a closed subspace $\mathcal X\subset L^2\Lambda^k(\M)$,
the $L^2$-orthogonal projection
$\Pi_{\mathcal X}:L^2\Lambda^k(\M)\to\mathcal X$
is characterized by
\begin{equation}
\label{eq:l2-orthogonal-projection}
  \langle\Pi_{\mathcal X}\omega,\zeta\rangle_{\M}
  =\langle\omega,\zeta\rangle_{\M}
  =\int_{\M}\omega\wedge\star\zeta,
  \qquad \zeta\in\mathcal X.
\end{equation}
Applying this identity to $\star\omega$, with
$\mathcal X\subset L^2\Lambda^{m-k}(\M)$, gives
\cite[Lemma~2.1]{Bonaldi2025}
\begin{equation}
\label{eq:hodge-projection-moments}
  \langle\Pi_{\mathcal X}(\star\omega),\zeta\rangle_{\M}
  =\int_{\M}\star^{-1}
       \Pi_{\mathcal X}(\star\omega)\wedge\zeta
  =\int_{\M}\omega\wedge\zeta,
\end{equation}
for all $\omega\in L^2\Lambda^k(\M)$ and $\zeta\in\mathcal X$.

For $s\ge0$, denote the coefficientwise $H^s$ norm by
$\norm{\cdot}_{s,\M}$ and, for $s>0$, the associated seminorm by
$\seminorm{\cdot}_{s,\M}$.
Let $H_0^s\Lambda^k(\M)$ be the closure of
$C_c^\infty\Lambda^k(\M)$ in $H^s\Lambda^k(\M)$.
For $s>0$, define
$H^{-s}\Lambda^k(\M):=(H_0^s\Lambda^k(\M))'$
with respect to the $L^2$ pivot, and denote its norm by
$\norm{\cdot}_{-s,\M}$.
Angle brackets also denote the corresponding duality pairings.
When $\M=\Omega$, we omit the domain subscript in inner products,
norms, and seminorms.

Interpreting $\rd$ and $\delta$ distributionally, define the graph spaces
\begin{equation*}
\begin{aligned}
  H\Lambda^k(\rd;\M)
  &:=\{\omega\in L^2\Lambda^k(\M):
                   \rd\omega\in L^2\Lambda^{k+1}(\M)\},\\
  H\Lambda^k(\delta;\M)
  &:=\{\omega\in L^2\Lambda^k(\M):
                   \delta\omega\in L^2\Lambda^{k-1}(\M)\}.
\end{aligned}
\end{equation*}
For $D=\rd$ or $\delta$, equip $H\Lambda^k(D;\M)$ with the graph norm
$\norm{\omega}_{H\Lambda^k(D;\M)}^2
:=\norm{\omega}_{\M}^2+\norm{D\omega}_{\M}^2$.
Define the homogeneous space as the closure in this graph norm:
\[
  \overset{\circ}{H}\Lambda^k(D;\M)
  :=\overline{C_c^\infty\Lambda^k(\M)}^{\,H\Lambda^k(D;\M)}.
\]
The spaces $H\Lambda^k(D;\partial\M)$ and their graph norms
are defined analogously using distributional derivatives on
$\partial\M$.
These spaces give the domain complexes
\[
  (H\Lambda^\bullet(D;\M),D),
  \qquad
  (H\Lambda^\bullet(D;\partial\M),D),
  \qquad
  (\overset{\circ}{H}\Lambda^\bullet(D;\M),D).
\]
These are Hilbert complexes: their differentials are closed and densely
defined, and their domains are Hilbert spaces in the graph norms.

\subsubsection{Weak traces}
\label{subsubsec:weak-traces}
\label{subsec:derham-trace}
The tangential and normal traces extend continuously from smooth forms
to $H^1$ forms.
Their ranges
\[
  H_{\parallel}^{1/2}\Lambda^k(\partial\M)
  :=\trc_kH^1\Lambda^k(\M),
  \qquad
  H_{\perp}^{1/2}\Lambda^k(\partial\M)
  :=\ntr_{k+1}H^1\Lambda^{k+1}(\M)
\]
are equipped with the quotient norms induced by $H^1$.
We denote their $L^2$-pivot duals by
$H_{\parallel}^{-1/2}\Lambda^k(\partial\M)$ and
$H_{\perp}^{-1/2}\Lambda^k(\partial\M)$, respectively.
The distributional boundary operators $\rd$ and $\delta$ are
closed, densely defined unbounded operators on
$H_{\perp}^{-1/2}\Lambda^\bullet(\partial\M)$ and
$H_{\parallel}^{-1/2}\Lambda^\bullet(\partial\M)$, respectively.
Their domains are the graph spaces
\begin{equation*}
\begin{aligned}
  H_{\perp}^{-1/2}\Lambda^k(\rd;\partial\M)
  &:=\{\omega\in H_{\perp}^{-1/2}\Lambda^k(\partial\M):
       \rd\omega\in H_{\perp}^{-1/2}\Lambda^{k+1}(\partial\M)\},\\
  H_{\parallel}^{-1/2}\Lambda^k(\delta;\partial\M)
  &:=\{\omega\in H_{\parallel}^{-1/2}\Lambda^k(\partial\M):
       \delta\omega\in H_{\parallel}^{-1/2}\Lambda^{k-1}(\partial\M)\}.
\end{aligned}
\end{equation*}
The associated complexes are
$(H_{\perp}^{-1/2}\Lambda^\bullet(\rd;\partial\M),\rd)$ and
$(H_{\parallel}^{-1/2}\Lambda^\bullet(\delta;\partial\M),\delta)$.
\par
The tangential and normal traces extend to bounded surjections
\cite[Eq.~(5.1.9)]{Schulz2022}
\begin{equation}
\label{eq:de-rham-trace}
  \trc_k:H\Lambda^k(\rd;\M)
  \twoheadrightarrow H_{\perp}^{-1/2}\Lambda^k(\rd;\partial\M),
  \qquad
  \ntr_k:H\Lambda^k(\delta;\M)
  \twoheadrightarrow H_{\parallel}^{-1/2}\Lambda^{k-1}(\delta;\partial\M).
\end{equation}
They satisfy the commutation and anticommutation relations
\begin{equation}
\label{eq:trace-commutation}
  \trc_{k+1}\rd^k=\rd^k\trc_k,
  \qquad
  \ntr_{k-1}\delta_k=-\delta_{k-1}\ntr_k.
\end{equation}
The trace-kernel characterizations follow from
\cite[Theorem~6.4]{Arnold2018} and Hodge duality:
\begin{equation}
\label{eq:trace-kernels}
  \overset{\circ}{H}\Lambda^k(\rd;\M)=\ker\trc_k,
  \qquad
  \overset{\circ}{H}\Lambda^k(\delta;\M)
  =\star_{\M}^{-1}\overset{\circ}{H}\Lambda^{m-k}(\rd;\M)
  =\ker\ntr_k.
\end{equation}

The integration-by-parts formula \eqref{eq:stokes-integration-by-parts},
expressed using the Hodge star, extends to
$\omega\in H\Lambda^k(\rd;\M)$ and
$\zeta\in H\Lambda^{k+1}(\delta;\M)$ in the form
\cite[Eq.~(5.1.10)]{Schulz2022}
\begin{equation}
\label{eq:green-formula}
  \langle\rd^k\omega,\zeta\rangle_{\M}
  =\langle\omega,\delta_{k+1}\zeta\rangle_{\M}
   +\tracepair{\trc_k\omega}{\ntr_{k+1}\zeta}_{\partial\M}.
\end{equation}
Here $\tracepair{\cdot}{\cdot}_{\partial\M}$
denotes the continuous duality pairing between the trace spaces
$H_{\perp}^{-1/2}\Lambda^k(\rd;\partial\M)$ and
$H_{\parallel}^{-1/2}\Lambda^k(\delta;\partial\M)$,
recovering the usual $H^{-1/2}$--$H^{1/2}$ duality and the
$L^2$ inner product as regularity permits.
Together with \eqref{eq:trace-kernels}, this formula yields the Hilbert
adjoint identities
\cite[Theorem~6.5]{Arnold2018}
\begin{equation}
\label{eq:adjoint-realizations}
  (\rd,H\Lambda^k(\rd;\M))^*
  =(\delta,\overset{\circ}{H}\Lambda^{k+1}(\delta;\M)),
  \qquad
  (\rd,\overset{\circ}{H}\Lambda^k(\rd;\M))^*
  =(\delta,H\Lambda^{k+1}(\delta;\M)).
\end{equation}
Taking adjoints gives the reverse identities.

\subsubsection{Hodge decompositions, Poincar\'e inequality, and Hodge Laplacian}
\label{subsubsec:hodge-theory}
We define the spaces of exact, coexact, and closed forms by
\begin{equation}
\label{eq:exact-coexact-spaces}
\begin{gathered}
  \mathcal B^k(\M):=\rd H\Lambda^{k-1}(\rd;\M),
  \qquad
  \mathcal B_k^*(\M):=\delta H\Lambda^{k+1}(\delta;\M),\\
  \mathcal Z^k(\M):=\ker(\rd|_{H\Lambda^k(\rd;\M)}).
\end{gathered}
\end{equation}
Define the harmonic space by
$\Harm^k(\M):=\mathcal Z^k(\M)\cap\mathcal B^k(\M)^\perp$,
where orthogonality is taken in $L^2$.
A ring denotes the analogous spaces for the homogeneous complex.
These spaces are closed in $L^2\Lambda^k(\M)$, and the adjoint
identities \eqref{eq:adjoint-realizations} yield
\begin{align}
\label{eq:hodge-decomposition}
  L^2\Lambda^k(\M)
  &=\mathcal B^k(\M)\oplus\Harm^k(\M)
    \oplus\mathring{\mathcal B}_k^*(\M),\\
\label{eq:hodge-decomposition-homogeneous}
  L^2\Lambda^k(\M)
  &=\mathring{\mathcal B}^k(\M)
    \oplus\mathring{\Harm}^k(\M)
    \oplus\mathcal B_k^*(\M).
\end{align}
Define the $L^2$-orthogonal complements of the closed forms by
\[
  \mathcal Z^{k\perp}(\M)
  :=H\Lambda^k(\rd;\M)\cap\mathring{\mathcal B}_k^*(\M),
  \qquad
  \mathring{\mathcal Z}^{k\perp}(\M)
  :=\overset{\circ}{H}\Lambda^k(\rd;\M)
       \cap\mathcal B_k^*(\M).
\]
Intersecting the $L^2$ decompositions with the respective graph spaces yields
\begin{align}
\label{eq:domain-hodge-decomposition}
  H\Lambda^k(\rd;\M)
  &=\mathcal B^k(\M)\oplus\Harm^k(\M)
    \oplus\mathcal Z^{k\perp}(\M),\\
\label{eq:domain-hodge-decomposition-homogeneous}
  \overset{\circ}{H}\Lambda^k(\rd;\M)
  &=\mathring{\mathcal B}^k(\M)
    \oplus\mathring{\Harm}^k(\M)
    \oplus\mathring{\mathcal Z}^{k\perp}(\M).
\end{align}

The Poincar\'e inequality \cite[Theorem~4.6]{Arnold2018}
holds with a constant $C_{\mathrm P}\ge1$ depending on $\M$
and the form degree:
\begin{equation}
\label{eq:gaffney}
  \norm{\omega}_{H\Lambda^k(\rd;\M)}
  \le C_{\mathrm P}\norm{\rd\omega}_{\M},
  \qquad
  \omega\in\mathcal Z^{k\perp}(\M)
       \cup\mathring{\mathcal Z}^{k\perp}(\M).
\end{equation}

Let $(V^\bullet(\M),\rd)$ denote either
$(H\Lambda^\bullet(\rd;\M),\rd)$ or
$(\overset{\circ}{H}\Lambda^\bullet(\rd;\M),\rd)$,
with Hilbert adjoint complex $(V_\bullet^*(\M),\delta)$ given by
\eqref{eq:adjoint-realizations}.
The corresponding Hodge Laplacian $L^k\omega:=\rd\delta\omega+\delta\rd\omega$
has domain
\begin{equation}
\label{eq:hodge-laplacian-domain}
  \mathcal D(L^k):=\bigl\{\omega\in V^k(\M)\cap V_k^*(\M):
    \delta\omega\in V^{k-1}(\M),\quad
    \rd\omega\in V_{k+1}^*(\M)\bigr\}.
\end{equation}
Writing $\Harm_V^k(\M)$ for the corresponding harmonic space,
the strong problem for $j\in L^2\Lambda^k(\M)$ is
\begin{equation}
\label{eq:strong-hodge-laplace}
  L^k u=j-\Pi_{\Harm_V^k(\M)}j,
  \qquad u\in\mathcal D(L^k),
  \qquad u\perp\Harm_V^k(\M).
\end{equation}
This problem is well posed \cite[Theorem~4.8]{Arnold2018},
and its equivalent mixed formulation is presented and discretized
in Section~\ref{sec:hodge-laplace-vem}.
The terms $\rd\delta u$ and $\delta\rd u$ are the exact and
coexact components of $j$ in the corresponding Hodge decomposition.

\subsection{Discrete setting}
\label{subsec:discrete-setting}
We introduce polytopal meshes and trace compatibility, recall local
polynomial differential forms~\cite{Arnold2018,Arnold2006},
and formulate the local boundary value problems used in the recursive
construction of Section~\ref{sec:veec}.

\subsubsection{Polytopal meshes and trace compatibility}
\label{subsec:polytopal-mesh-trace-compatibility}
\label{subsubsec:trace-compatibility}
For an integer $d$ with $1\le d\le n$, a \emph{flat $d$-cell}
is a bounded polytopal domain in a $d$-dimensional affine space,
with Lipschitz boundary and closure homeomorphic to the closed
$d$-ball. A \emph{flat $0$-cell} is a point.
A \emph{polytopal mesh} $\Mh$ of $\Omega$ is a finite partition of
$\overline\Omega$ into such cells. Cell boundaries are unions of
lower-dimensional mesh cells, and intersections of cell closures
are unions of mesh cells. Every $d$-cell with $d<n$ is a
\emph{facet} of some $(d+1)$-cell.
In the terminology of algebraic topology, these conditions make
$\overline\Omega$ a finite regular cellular complex whose open
cells are precisely the members of $\Mh$; see \cite{Hatcher2002}.
\par
For $0\le d\le n$, let $\Delta_d(\Mh):=\{f\in\Mh:\dim f=d\}$
denote the set of $d$-cells, and set $\Th:=\Delta_n(\Mh)$.
For $f\in\Delta_d(\Mh)$ and $0\le\ell\le d$, define
$\Delta_\ell(f):=\{e\in\Delta_\ell(\Mh):e\subseteq\overline f\}$.
\par
Fix an orientation on each mesh cell. For $f\in\Delta_d(\Mh)$
with $d\ge1$ and $e\in\Delta_{d-1}(f)$, let $\varepsilon_{fe}=1$
if the fixed orientation of $e$ agrees with the boundary orientation
induced by $f$, and $-1$ otherwise.
For $\zeta\in L^2\Lambda^{d-1}(\partial f)$, boundary integration
is given by
\begin{equation}
\label{eq:integral.relative.orientation}
  \int_{\partial f}\zeta
  =
  \sum_{e\in\Delta_{d-1}(f)}
  \varepsilon_{fe}\int_e\zeta|_e.
\end{equation}
\par

Let $\omega$ be a $k$-form defined piecewise on two $d$-cells
$f_1,f_2$ sharing a $(d-1)$-cell $e$, with
$\omega|_{f_i}\in H\Lambda^k(\rd;f_i)$ for $i=1,2$.
The weak Stokes formula shows that
$\omega\in H\Lambda^k(\rd;f_1\cup f_2\cup e)$ if and only if
its tangential trace is \emph{single-valued} on $e$:
\begin{equation}
\label{eq:single-valued-tangential-trace}
  \trc_e^{f_1}\omega
  =
  \trc_e^{f_2}\omega,
\end{equation}
Here $\trc_e^f\omega$ denotes the weak tangential trace of
$\omega|_f$ on $e$, expressed in the fixed orientation of $e$.
We set $\trc_f^f:=\operatorname{id}$ on forms defined on $f$
and omit the superscript when unambiguous.
For $e'\in\Delta_\ell(f)$ with $\ell\le d-2$, consider a chain
$f_\ell=e',\ldots,f_d=f$ with $f_j\in\Delta_j(f_{j+1})$
for $\ell\le j<d$.
If all successive facet traces are defined and their composition
is independent of the chosen chain, we set
\[
  \trc_{e'}^f
  :=\trc_{f_\ell}^{f_{\ell+1}}\circ\cdots\circ
    \trc_{f_{d-1}}^{f_d}.
\]
Under these assumptions on both $f_1$ and $f_2$,
applying $\trc_{e'}^e$ to \eqref{eq:single-valued-tangential-trace} gives
$\trc_{e'}^{f_1}\omega
 =\trc_{e'}^{f_2}\omega$
for $e'\in\Delta_\ell(e)$ with $k\le\ell<d-1$.
This leads to the following characterization; see
\cite[Lemma~5.1]{Arnold2006} for the piecewise smooth case.
\par
\begin{lemma}[Trace compatibility]
\label{lem:trace-gluing}
Let $\G=\Omega$ or $\G=\partial f$ for a cell
$f\in\Delta_d(\Mh)$ with $d\ge2$.
Set $m:=\dim\G$ and
$\Mh|_{\G}:=\{e\in\Mh:e\subseteq\overline\G\}$.
Let $\omega\in L^2\Lambda^k(\G)$.
Assume that, for each $f'\in\Delta_m(\Mh|_{\G})$,
$\omega|_{f'}$ admits chain-independent iterated tangential traces
satisfying
\[
  \trc_e^{f'}\omega
  \in H\Lambda^k(\rd;e),
  \qquad e\in\Delta_\ell(f'),\quad k\le \ell\le m.
\]
Then the following statements are equivalent:
\begin{enumerate}
\item[(i)] $\omega\in H\Lambda^k(\rd;\G)$.
\item[(ii)] $\omega$ has single-valued tangential traces
on every $(m-1)$-cell of $\Mh|_{\G}$.
\item[(iii)] $\omega$ has single-valued tangential traces
on every $\ell$-cell of $\Mh|_{\G}$ with $k\le \ell<m$.
\end{enumerate}
\end{lemma}

\subsubsection{Polynomial differential forms}
\label{subsec:polynomial-forms}
For integers $r\ge0$ and $0\le k\le d$ and a cell
$f\in\Delta_d(\Mh)$, let $\mathcal H_r\Lambda^k(f)$ denote the space
of $k$-forms with homogeneous polynomial coefficients of degree $r$,
with $\mathcal H_j\Lambda^k(f):=\{0\}$ for $j<0$.
Let $x_f$ denote the centroid of $f$, given by
$x_f=|f|^{-1}\int_f x\,dx$ for $d\ge1$, and set $x_f=f$ for a vertex.
Here $|f|$ denotes the $d$-dimensional measure of $f$.
With $x_f$ as the origin, we define the Koszul operator $\kappa$
by contraction \eqref{eq:form-contraction} with the radial vector field
$x\mapsto x-x_f$.
For $r+k>0$, the homogeneous spaces admit the decomposition
\cite[(3.10)]{Arnold2006}
\begin{equation}
\label{eq:homogeneous-koszul-decomposition}
\mathcal H_r\Lambda^k(f)
=
\kappa\mathcal H_{r-1}\Lambda^{k+1}(f)
\oplus
\rd\mathcal H_{r+1}\Lambda^{k-1}(f).
\end{equation}
\par
Let $\Pcal_r\Lambda^k(f)$ denote the space of $k$-forms with
polynomial coefficients of degree at most $r$, with
$\Pcal_j\Lambda^k(f):=\{0\}$ for $j<0$.
Then $\Pcal_r\Lambda^k(f)=\bigoplus_{j=0}^{r}\mathcal H_j\Lambda^k(f)$,
and summing \eqref{eq:homogeneous-koszul-decomposition} over
$j=0,\ldots,r$ gives, for $1\le k\le d$,
\begin{equation}
\label{eq:complete-polynomial-decomposition}
\Pcal_r\Lambda^k(f)
=
\kappa\Pcal_{r-1}\Lambda^{k+1}(f)
\oplus
\rd\Pcal_{r+1}\Lambda^{k-1}(f).
\end{equation}
For $k=0$, the complete space decomposes as
\begin{equation}
\label{eq:polynomial-koszul-decomposition-k0}
\Pcal_r\Lambda^0(f)=\Pcal_0\Lambda^0(f)
\oplus\kappa\Pcal_{r-1}\Lambda^1(f).
\end{equation}
Define the Koszul image spaces by
$\Kcal_r^k(f):=\kappa\Pcal_{r-1}\Lambda^{k+1}(f)$,
with $\Kcal_0^k(f)=\{0\}$ and $\Kcal_r^d(f)=\{0\}$.
Splitting by polynomial degree gives, for $r\ge1$,
\begin{equation}
\label{eq:koszul-degree-splitting}
  \Kcal_r^k(f)=\Kcal_{r-1}^k(f)
  \oplus\kappa\mathcal H_{r-1}\Lambda^{k+1}(f).
\end{equation}
For $r\ge0$ and $1\le k\le d$, decomposing
$\Pcal_{r+1}\Lambda^{k-1}(f)$ by
\eqref{eq:complete-polynomial-decomposition}--\eqref{eq:polynomial-koszul-decomposition-k0}
and applying $\rd$ gives
$\rd\Pcal_{r+1}\Lambda^{k-1}(f)=\rd\Kcal_{r+1}^{k-1}(f)$.
Substituting this identity into
\eqref{eq:complete-polynomial-decomposition} yields
\[
  \Pcal_r\Lambda^k(f)
  =\rd\Kcal_{r+1}^{k-1}(f)\oplus\Kcal_r^k(f).
\]
Together with the injectivity of $\rd$ on $\Kcal_{r+1}^{k-1}(f)$
\cite[Theorem~3.2]{Arnold2006}, this gives the isomorphism
\begin{equation}
\label{eq:polynomial-koszul-isomorphism}
\Kcal_{r+1}^{k-1}(f)\times\Kcal_r^k(f)
\xrightarrow{\cong}\Pcal_r\Lambda^k(f),
\qquad (\mu,\nu)\mapsto\rd\mu+\nu.
\end{equation}
\par
For $r\ge1$, the trimmed polynomial space is
\cite[Section~7.3]{Arnold2018}
\begin{equation}
\label{eq:trimmed-polynomial-space}
\Pcal_r^-\Lambda^k(f)
:=\Pcal_{r-1}\Lambda^k(f)
\oplus\kappa\mathcal H_{r-1}\Lambda^{k+1}(f).
\end{equation}
At $r=0$, we set $\Pcal_0^-\Lambda^0(f):=\Pcal_0\Lambda^0(f)$ and
$\Pcal_0^-\Lambda^k(f):=\{0\}$ for $k\ge1$.
For $r\ge1$ and $1\le k\le d$, decomposing the first summand in
\eqref{eq:trimmed-polynomial-space} by
\eqref{eq:polynomial-koszul-isomorphism} and combining the Koszul
components using \eqref{eq:koszul-degree-splitting} gives the isomorphism
\begin{equation}
\label{eq:trimmed-koszul-isomorphism}
\Kcal_r^{k-1}(f)\times\Kcal_r^k(f)
\xrightarrow{\cong}\Pcal_r^-\Lambda^k(f),
\qquad (\mu,\nu)\mapsto\rd\mu+\nu,
\end{equation}
which extends to $r=0$ by convention.
In particular, for $r\ge0$,
\begin{equation}
\label{eq:trimmed-polynomial-top-degree}
\Pcal_r^-\Lambda^0(f)=\Pcal_r\Lambda^0(f),
\qquad
\Pcal_r^-\Lambda^d(f)=\Pcal_{r-1}\Lambda^d(f)\quad(d\ge1).
\end{equation}
For $r\ge0$ and $0\le k<d$, applying $\rd$ to the
isomorphisms in \eqref{eq:polynomial-koszul-isomorphism} and
\eqref{eq:trimmed-koszul-isomorphism}, and using $\rd^2=0$, gives
\begin{equation}
\label{eq:polynomial-differential-images}
\rd\Pcal_r^-\Lambda^k(f)=\rd\Pcal_r\Lambda^k(f)=\rd\Kcal_r^k(f)
\subseteq\Pcal_{r-1}\Lambda^{k+1}(f)
\subseteq\Pcal_r^-\Lambda^{k+1}(f).
\end{equation}
For $1\le k\le d$, the defining polynomial inclusions,
\eqref{eq:polynomial-differential-images}, and Hodge duality give
\begin{equation}
\label{eq:polynomial-codifferential-image}
\delta\Pcal_r^-\Lambda^k(f)
\subseteq\delta\Pcal_r\Lambda^k(f)
=\star^{-1}\rd\Kcal_r^{d-k}(f)
\subseteq\Pcal_{r-1}\Lambda^{k-1}(f)
\subseteq\Pcal_r^-\Lambda^{k-1}(f).
\end{equation}
\par
We abbreviate the $L^2$-orthogonal projections defined by
\eqref{eq:l2-orthogonal-projection} as
\[
  \Proj{r}{f}{k}
  :=
  \Pi_{\Pcal_r\Lambda^k(f)},
  \qquad
  \Pi_{r,f}^{\kappa,k}
  :=
  \Pi_{\Kcal_r^k(f)},
  \qquad
  \Pi_{r,f}^{-,k}
  :=
  \Pi_{\Pcal_r^-\Lambda^k(f)}.
\]
\par

\begin{lemma}[Trace preservation {\cite[Lemma~2.4]{Bonaldi2025}}]
\label{lem:traces-trimmed}
Let $f\in\Delta_d(\Mh)$ and
$e\in\Delta_{\ell}(f)$, with $0\le k\le \ell\le d$. Then, for every
$r\ge0$,
\[
  \trc_e\Pcal_r^-\Lambda^k(f)
  \subset
  \Pcal_r^-\Lambda^k(e).
\]
\end{lemma}

\subsubsection{Local boundary value problems}
\label{subsubsec:local-boundary-problems}
Let $f\in\Delta_d(\Mh)$ with $1\le d\le n$.
We take $V^\bullet(f)=\overset{\circ}{H}\Lambda^\bullet(\rd;f)$,
whose coexact spaces are $\mathcal B_k^*(f)$ from
\eqref{eq:exact-coexact-spaces}.
For a source $\mathfrak{f}^k\in\mathcal B_k^*(f)$, the solution $u^k$
of \eqref{eq:strong-hodge-laplace} with $\M=f$
also belongs to $\mathcal B_k^*(f)$ and satisfies
\[
  \delta\rd u^k=\mathfrak{f}^k,\qquad
  \delta u^k=0,\qquad \trc u^k=0.
\]
This is the $\mathcal B^*$ problem
\cite[Section~4.4.3]{Arnold2018}.
The local problem below generalizes it by prescribing both the
codifferential and the tangential trace.

\begin{definition}[$\mathcal{NB}_k^*$ problem]
\label{def:NBkstar}
For $0\le k\le d-1$ and data
\[
  (\mathfrak{f}^k,\mathfrak{f}^{k-1},\mathfrak{b}^k)
  \in
  \mathcal B_k^*(f)
  \times
  \mathcal B_{k-1}^*(f)
  \times
  H_{\perp}^{-1/2}
  \Lambda^k(\rd;\partial f),
\]
seek
$u^k\in H\Lambda^k(\rd;f)
\cap H\Lambda^k(\delta;f)$ with
$\rd u^k\in
H\Lambda^{k+1}(\delta;f)$ satisfying
\begin{equation}
\label{eq:NBkstar}
  \delta\rd u^k=\mathfrak{f}^k
  \quad\text{in }f,
  \qquad
  \delta u^k=\mathfrak{f}^{k-1}
  \quad\text{in }f,
  \qquad
  \trc u^k=\mathfrak{b}^k
  \quad\text{on }\partial f.
\end{equation}
\par
For $k=d$ and
$(\mathfrak{f}^{d-1},\mathfrak{b}^d)
\in\mathcal B_{d-1}^*(f)\times\mathbb R$,
seek $u^d\in H\Lambda^d(\delta;f)$
satisfying
\[
  \delta u^d=\mathfrak{f}^{d-1},
  \qquad
  \int_{f}u^d=\mathfrak{b}^d.
\]
\end{definition}
\begin{theorem}[Well-posedness]
\label{thm:NBkstar-wellposed}
The $\mathcal{NB}_k^*$ problem has a unique solution for every
admissible data set in Definition~\ref{def:NBkstar}. For
$0\le k\le d-1$, there is a constant $C_{\mathrm{NB}}>0$
depending only on $f$ such that
\begin{equation}
\label{eq:NBkstar-stability}
  \|u^k\|_{H\Lambda^k(\rd;f)}
  \le
  C_{\mathrm{NB}}
  \left(
    \|\mathfrak{f}^k\|_{f}
    +\|\mathfrak{f}^{k-1}\|_{f}
    +\|\mathfrak{b}^k\|_{
      H_{\perp}^{-1/2}\Lambda^k(\rd;\partial f)}
  \right).
\end{equation}
\end{theorem}
\begin{proof}
Let $0\le k\le d-1$. The bounded surjection $\trc_k$ in
\eqref{eq:de-rham-trace} admits a bounded linear right inverse
$\mathcal L_k$. Set
\begin{equation}
\label{eq:NBkstar-trace-lifting}
  \psi^k:=\mathcal L_k \mathfrak{b}^k,
  \qquad
  \trc \psi^k=\mathfrak{b}^k,
  \qquad
  \|\psi^k\|_{H\Lambda^k(\rd;f)}
  \le
  C_{\mathrm{tr}}
  \|\mathfrak{b}^k\|_{
    H_{\perp}^{-1/2}\Lambda^k(\rd;\partial f)}.
\end{equation}
Since $f$ is contractible, Hodge duality gives
$\mathring{\Harm}^\ell(f)=\{0\}$ for $0\le\ell<d$.
Thus, \eqref{eq:domain-hodge-decomposition-homogeneous} reduces to
\begin{equation}
\label{eq:NBkstar-hodge-decomposition}
  \overset{\circ}{H}\Lambda^\ell(\rd;f)
  =\mathring{\mathcal B}^\ell(f)
    \mathbin{\overset{\perp}{\oplus}}\mathring{\mathcal Z}^{\ell\perp}(f),
  \qquad 0\le\ell<d.
\end{equation}
For $\ell=k,k-1$, the space $\mathring{\mathcal Z}^{\ell\perp}(f)$ is closed in the norm
$\|\cdot\|_{H\Lambda^\ell(\rd;f)}$,
and the Poincar\'e inequality \eqref{eq:gaffney} makes
$\langle\rd\cdot,\rd\cdot\rangle_{f}$ coercive on it.
Hence the Lax--Milgram theorem yields unique forms
$\phi^\ell\in \mathring{\mathcal Z}^{\ell\perp}(f)$ satisfying
\begin{align}
\label{eq:NBkstar-correction-km1}
  \langle\rd \phi^{k-1},\rd \zeta\rangle_{f}
  &=
  \langle \mathfrak{f}^{k-1},\zeta\rangle_{f}
  -
  \langle \psi^k,\rd \zeta\rangle_{f}
  &&\forall \zeta\in \mathring{\mathcal Z}^{k-1\perp}(f),
  \\
\label{eq:NBkstar-correction-k}
  \langle\rd \phi^k,\rd \zeta\rangle_{f}
  &=
  \langle \mathfrak{f}^k,\zeta\rangle_{f}
  -
  \langle\rd \psi^k,\rd \zeta\rangle_{f}
  &&\forall \zeta\in \mathring{\mathcal Z}^{k\perp}(f).
\end{align}
Testing \eqref{eq:NBkstar-correction-km1}--\eqref{eq:NBkstar-correction-k}
with $\zeta=\phi^{k-1}$ and $\zeta=\phi^k$, respectively, and using
\eqref{eq:gaffney} gives
\begin{equation}
\label{eq:NBkstar-correction-estimate}
  \|\phi^k\|_{H\Lambda^k(\rd;f)}
  +
  \|\phi^{k-1}\|_{H\Lambda^{k-1}(\rd;f)}
  \le
  C_{\mathrm B}
  \left(
    \|\mathfrak{f}^k\|_{f}
    +
    \|\mathfrak{f}^{k-1}\|_{f}
    +
    \|\psi^k\|_{H\Lambda^k(\rd;f)}
  \right),
\end{equation}
where $C_{\mathrm B}>0$ depends only on
$f$.
Since $\mathfrak{f}^\ell\in\mathcal B_\ell^*(f)$ and $\rd^2=0$, the orthogonal
decomposition \eqref{eq:NBkstar-hodge-decomposition} extends the test spaces in
\eqref{eq:NBkstar-correction-km1} and
\eqref{eq:NBkstar-correction-k} from $\mathring{\mathcal Z}^{\ell\perp}(f)$ to
$\overset{\circ}{H}\Lambda^\ell(\rd;f)$, for $\ell=k-1,k$, respectively.
\par
Set $u^k:=\psi^k+\phi^k+\rd \phi^{k-1}$.
Since $\phi^k\in\mathcal B_k^*(f)$, we have $\delta \phi^k=0$.
The adjoint identities \eqref{eq:adjoint-realizations}, together with
$\rd^2=0$, then give
\[
  \delta u^k=\delta(\psi^k+\rd \phi^{k-1})=\mathfrak{f}^{k-1},
  \qquad
  \delta\rd u^k=\delta\rd(\psi^k+\phi^k)=\mathfrak{f}^k.
\]
Since $\trc \phi^k=0$ and
$\trc(\rd \phi^{k-1})=\rd\trc \phi^{k-1}=0$,
we have $\trc u^k=\mathfrak{b}^k$.
The definition of $u^k$ and the estimates
\eqref{eq:NBkstar-trace-lifting} and
\eqref{eq:NBkstar-correction-estimate} yield
\eqref{eq:NBkstar-stability}.
For uniqueness, let $u^k$ solve the homogeneous problem.
Since $\trc u^k=0$ and $\delta u^k=0$, the decomposition
\eqref{eq:NBkstar-hodge-decomposition} gives
$u^k\in \mathring{\mathcal Z}^{k\perp}(f)$. Testing $\delta\rd u^k=0$ with $u^k$ and applying
\eqref{eq:gaffney} yield $u^k=0$.
\par

For $k=d$, the equation $\delta u^d=\mathfrak{f}^{d-1}$ is solvable by the
definition of $\mathcal B_{d-1}^*(f)$. Since $f$ is connected,
its solutions differ by multiples of $\star1$, and the prescribed
integral $\int_{f}u^d=\mathfrak{b}^d$ selects a unique solution.
\end{proof}
\begin{corollary}[Compatible potentials]
\label{cor:NBkstar-differential-compatibility}
Let $0\le k<d$ and let $u^k$ solve the $\mathcal{NB}_k^*$ problem
with data $(\mathfrak{f}^k,\mathfrak{f}^{k-1},\mathfrak{b}^k)$. Let $u^{k+1}$ solve the
$\mathcal{NB}_{k+1}^*$ problem with data
$(0,\mathfrak{f}^k,\rd\mathfrak{b}^k)$ for $k+1<d$ and
$\bigl(\mathfrak{f}^k,\langle\mathfrak{b}^k,
\star_{\partial f}^{-1}1\rangle_{\partial f}\bigr)$ for $k+1=d$.
Then $u^{k+1}=\rd u^k$.
\end{corollary}
\begin{proof}
By \eqref{eq:NBkstar} and $\rd^2=0$, $\rd u^k$ satisfies the interior
equations. Trace commutation gives
$\trc(\rd u^k)=\rd\mathfrak{b}^k$ for $k+1<d$,
while the weak Stokes formula gives
$\int_f\rd u^k
=\langle\mathfrak{b}^k,\star_{\partial f}^{-1}1\rangle_{\partial f}$
for $k+1=d$. Uniqueness then proves the claim.
\end{proof}

\section{Virtual element de Rham complexes}
\label{sec:veec}
Fix a polynomial degree $r\ge0$. We
construct two conforming virtual
element subcomplexes of the de Rham complex
\eqref{eq:domain-derham-complex}: a
trimmed family $\Vcal_{r+1}^{-,\bullet}$ containing
$\Pcal_{r+1}^{-}\Lambda^\bullet$ and a
complete family $\Vcal_{r}^{\bullet}$ containing
$\Pcal_r\Lambda^\bullet$.
\par

\subsection{Trimmed virtual element spaces of differential forms}
\label{subsec:trimmed-spaces}
We first construct the local spaces for the trimmed family, derive their
geometric decomposition and dimension formulas, and establish unisolvence
of the DOFs and computability of polynomial moments.
\subsubsection{Recursive local construction}
\label{subsubsec:trimmed-definition}
For all integers $0\le k\le d\le n$ and
$f\in\Delta_d(\Mh)$, we construct
$\Vcal_{r+1}^{-,k}(f)$ using the
$\mathcal{NB}_k^*$ problem of Definition~\ref{def:NBkstar}
with polynomial interior data and recursively defined boundary data,
so that:
\begin{itemize}
\item[(T1)] Conformity and polynomial inclusion:
$\Pcal_{r+1}^{-}\Lambda^k(f)
\subset\Vcal_{r+1}^{-,k}(f)
\subset H\Lambda^k(\rd;f)$.

\item[(T2)] Trace closure: if
$e\in\Delta_\ell(f)$ with $k\le\ell<d$, then
$\trc_e\Vcal_{r+1}^{-,k}(f)
\subset\Vcal_{r+1}^{-,k}(e)$.

\item[(T3)] Complex property:
$\rd\Vcal_{r+1}^{-,k}(f)
\subset
\Vcal_{r+1}^{-,k+1}(f)$,
$0\le k<d$.
\end{itemize}
The diagram below organizes the induction
by the diagonal index $\lambda:=d-k$
along a flag satisfying
$f_d\in\Delta_d(f_{d+1})$ for $0\le d<n$ and
$f_n\in\Delta_n(\Mh)$.
Its rows are the local complexes of \textup{(T3)},
and its vertical arrows encode the trace closure \textup{(T2)}.
Each diagonal collects the spaces with the same value of $\lambda$.
\begin{equation*}
\begin{tikzcd}[column sep=1.9em,row sep=1.8em]
\Vcal_{r+1}^{-,0}(f_n)
  \arrow[r,"\rd"] \arrow[d,"\trc_{f_{n-1}}"'] &
\Vcal_{r+1}^{-,1}(f_n)
  \arrow[r,"\rd"] \arrow[d,"\trc_{f_{n-1}}"'] &
\cdots
  \arrow[r,"\rd"] &
\Vcal_{r+1}^{-,n-1}(f_n)
  \arrow[r,"\rd"] \arrow[d,"\trc_{f_{n-1}}"'] &
\Vcal_{r+1}^{-,n}(f_n)
\\
\Vcal_{r+1}^{-,0}(f_{n-1})
  \arrow[r,"\rd"] \arrow[d,"\trc_{f_{n-2}}"'] &
\Vcal_{r+1}^{-,1}(f_{n-1})
  \arrow[r,"\rd"] \arrow[d,"\trc_{f_{n-2}}"'] &
\cdots
  \arrow[r,"\rd"] &
\Vcal_{r+1}^{-,n-1}(f_{n-1})
\\
\vdots
  \arrow[d,"\trc_{f_2}"'] &
\vdots
  \arrow[d,"\trc_{f_2}"'] &
\ddots &
\\
\Vcal_{r+1}^{-,0}(f_2)
  \arrow[r,"\rd"] \arrow[d,"\trc_{f_1}"'] &
\Vcal_{r+1}^{-,1}(f_2)
  \arrow[r,"\rd"] \arrow[d,"\trc_{f_1}"'] &
\Vcal_{r+1}^{-,2}(f_2)
\\
\Vcal_{r+1}^{-,0}(f_1)
  \arrow[r,"\rd"] \arrow[d,"\trc_{f_0}"'] &
\Vcal_{r+1}^{-,1}(f_1)
\\
\Vcal_{r+1}^{-,0}(f_0).
\end{tikzcd}
\end{equation*}
Accordingly, the diagonal $\lambda=0$ initializes the
construction, $\lambda=1$ is the first
nontrivial level, and the remaining
diagonals are obtained recursively.
\par
For $\lambda=0$, the form degree is $k=d$.
Using the top-degree case of \eqref{eq:trimmed-polynomial-top-degree}
with polynomial index $r+1$, we set
\begin{equation}
\label{eq:trimmed-top}
  \Vcal_{r+1}^{-,d}(f)
  :=
  \Pcal_{r+1}^{-}\Lambda^d(f)
  =
  \Pcal_r\Lambda^d(f).
\end{equation}
In particular, for $d=0$,
$\Vcal_{r+1}^{-,0}(f)=\mathbb R$ by convention.
\par
\begin{remark}[Endpoint formulation]
\label{rem:trimmed-endpoint-nb}
The space \eqref{eq:trimmed-top} is also
the solution space of the endpoint
$\mathcal{NB}_d^*$ problem
\[
  \delta\omega=\mathfrak{f}^{d-1},
  \qquad
  \int_f\omega=c,
  \qquad
  \mathfrak{f}^{d-1}\in\star^{-1}\rd\Kcal_r^0(f),
  \quad c\in\mathbb R .
\]
Indeed,
\eqref{eq:polynomial-codifferential-image} gives
$\delta\Pcal_r\Lambda^d(f)=\star^{-1}\rd\Kcal_r^0(f)$, so a
polynomial solution exists after fixing
the mean. Uniqueness in
Theorem~\ref{thm:NBkstar-wellposed}
identifies the endpoint solution space
with $\Pcal_r\Lambda^d(f)$.
\end{remark}
To choose data spaces for constructing the remaining virtual spaces
through the $\mathcal{NB}_k^*$ problem while ensuring the polynomial
inclusion in \textup{(T1)}, we make the following observation.
For $0\le k<d$, let $\eta\in\Pcal_{r+1}^{-}\Lambda^k(f)$.
By \eqref{eq:polynomial-differential-images},
$\rd\eta\in\Pcal_r\Lambda^{k+1}(f)$.
Applying \eqref{eq:polynomial-codifferential-image} and
Lemma~\ref{lem:traces-trimmed}, we obtain, for every $e\in\Delta_{d-1}(f)$,
\begin{equation}
\label{eq:trimmed-polynomial-data}
  \bigl(\delta\rd\eta,\delta\eta,\trc_e\eta\bigr)
  \in
  \star^{-1}\rd\Kcal_r^{d-k-1}(f)
  \times
  \star^{-1}\rd\Kcal_{r+1}^{d-k}(f)
  \times
  \Pcal_{r+1}^{-}\Lambda^k(e).
\end{equation}
On the diagonal $\lambda=1$, we specialize \eqref{eq:trimmed-polynomial-data}
to $k=d-1$ and use the resulting data spaces to define the local solution
space of the $\mathcal{NB}_{d-1}^{*}$ problem by
\begin{equation}
\label{eq:trimmed-codim-one}
\begin{aligned}
  \Vcal_{r+1}^{-,d-1}(f)
  :=
  \bigl\{\omega\in H\Lambda^{d-1}(\rd;f):\;&
  \delta\rd\omega\in\star^{-1}\rd\Kcal_r^0(f),\quad
  \delta\omega\in\star^{-1}\rd\Kcal_{r+1}^{1}(f),\\
  &\trc\omega
  \in\Vcal_{r+1}^{-,d-1}(\partial f)\bigr\},
\end{aligned}
\end{equation}
where the top-degree facet spaces in \eqref{eq:trimmed-top} give the boundary space
\begin{equation}
\label{eq:trimmed-codim-one-boundary}
  \Vcal_{r+1}^{-,d-1}(\partial f)
  :=
  \left\{
  \mathfrak{b}^{d-1}\in L^2\Lambda^{d-1}(\partial f):
  \mathfrak{b}^{d-1}|_e\in\Pcal_r\Lambda^{d-1}(e)
  \quad
  \forall e\in\Delta_{d-1}(f)
  \right\}.
\end{equation}
The interior data are coexact by Hodge duality, and the boundary data in
\eqref{eq:trimmed-codim-one-boundary} are $L^2$ top-degree forms.
Theorem~\ref{thm:NBkstar-wellposed} therefore gives well-posedness.
The two inclusions in \textup{(T1)} follow from \eqref{eq:trimmed-polynomial-data}
and \eqref{eq:trimmed-codim-one}, respectively. The facetwise condition in
\eqref{eq:trimmed-codim-one-boundary} gives trace closure \textup{(T2)}.
For each $\omega\in\Vcal_{r+1}^{-,d-1}(f)$,
Corollary~\ref{cor:NBkstar-differential-compatibility} shows that
$\rd\omega$ solves the corresponding endpoint $\mathcal{NB}_{d}^{*}$ problem.
Since $\delta\rd\omega\in\star^{-1}\rd\Kcal_r^0(f)$ by
\eqref{eq:trimmed-codim-one}, Remark~\ref{rem:trimmed-endpoint-nb} gives
$\rd\omega\in\Pcal_r\Lambda^d(f)=\Vcal_{r+1}^{-,d}(f)$, proving \textup{(T3)}.
\par
Assume that all spaces through induction level $\lambda\ge1$
have been constructed and satisfy \textup{(T1)}--\textup{(T3)}.
For $f\in\Delta_d(\Mh)$ with $d-k=\lambda+1$, define
the compatible boundary space by
\begin{equation}
\label{eq:trimmed-boundary}
\begin{aligned}
  \Vcal_{r+1}^{-,k}(\partial f)
  :=
  \bigl\{
  \mathfrak{b}\in L^2\Lambda^k(\partial f):\;&
  \mathfrak{b}|_{e_1}
  \in\Vcal_{r+1}^{-,k}(e_1),\quad
  \trc_{e'}^{e_1}\mathfrak{b}
  =
  \trc_{e'}^{e_2}\mathfrak{b},\\
  &\forall e_1,e_2\in\Delta_{d-1}(f),\quad
  e'\in\Delta_{d-2}(e_1)\cap\Delta_{d-2}(e_2)
  \bigr\}.
\end{aligned}
\end{equation}
By the induction hypotheses \textup{(T1)} and \textup{(T2)},
Lemma~\ref{lem:trace-gluing} gives conformity and single-valued traces:
\begin{equation}
\label{eq:trimmed-boundary-traces}
  \Vcal_{r+1}^{-,k}(\partial f)\subset H\Lambda^k(\rd;\partial f),
  \qquad
  \trc_e\Vcal_{r+1}^{-,k}(\partial f)\subset\Vcal_{r+1}^{-,k}(e),
  \quad e\in\Delta_\ell(f),\ k\le\ell<d.
\end{equation}
Using the boundary space above and coexact polynomial interior data,
we define the local trimmed space through the
$\mathcal{NB}_k^*$ problem:
\begin{equation}
\label{eq:trimmed-general}
\begin{aligned}
  \Vcal_{r+1}^{-,k}(f)
  :=
  \bigl\{\omega\in H\Lambda^k(\rd;f):\;&
  \delta\rd\omega
  \in\star^{-1}\rd\Kcal_{r+1}^{d-k-1}(f),
  \quad
  \delta\omega
  \in\star^{-1}\rd\Kcal_{r+1}^{d-k}(f),
  \\
  &
  \trc\omega
  \in\Vcal_{r+1}^{-,k}(\partial f)
  \bigr\}.
\end{aligned}
\end{equation}
Theorem~\ref{thm:NBkstar-wellposed} therefore gives well-posedness and conformity.
Property \textup{(T2)} follows directly from
\eqref{eq:trimmed-boundary-traces}.
Moreover, on every facet $e\in\Delta_{d-1}(f)$,
the induction hypotheses \textup{(T1)} and \textup{(T3)} give
\[
  \Pcal_{r+1}^{-}\Lambda^k(e)
  \subset\Vcal_{r+1}^{-,k}(e),
  \qquad
  \rd\Vcal_{r+1}^{-,k}(e)
  \subset\Vcal_{r+1}^{-,k+1}(e).
\]
Applying $\Kcal_r^{d-k-1}(f)\subseteq\Kcal_{r+1}^{d-k-1}(f)$
and the first inclusion above to \eqref{eq:trimmed-polynomial-data}
establishes the polynomial inclusion in \textup{(T1)}.
The second inclusion and trace commutation,
$\trc_e\rd=\rd\trc_e$, give
$\rd\Vcal_{r+1}^{-,k}(\partial f)
\subset\Vcal_{r+1}^{-,k+1}(\partial f)$.
The choice of $r+1$ in \eqref{eq:trimmed-general} makes the
$\delta\rd$ data space coincide with the codifferential data space
at degree $k+1$.
Corollary~\ref{cor:NBkstar-differential-compatibility} therefore gives
$\rd\Vcal_{r+1}^{-,k}(f)\subset\Vcal_{r+1}^{-,k+1}(f)$,
completing the recursive construction.
\par
\begin{remark}[Relation to existing virtual element constructions]
\label{rem:trimmed-vector-proxies}
Table~\ref{tab:low-dimensional-vector-proxies} gives the vector proxies for
the differential and trace operators, and \eqref{eq:polynomial-codifferential-image}
identifies the interior data spaces for $\delta\rd\omega$ and $\delta\omega$
as images of polynomial differential form spaces under $\delta$.
Under these identifications, the trimmed virtual spaces recover
the differential and trace structure of the nodal, edge, and face
virtual elements in \cite{BeiraoEtAl2018},
with slightly different polynomial degrees and without
serendipity constraints.
For $k=0$, the construction also follows the recursive
boundary-value approach used in the $m=1$ case of the
arbitrary-dimensional scalar $H^m$-conforming virtual
elements in \cite{ChenHuangWei2022}.
In particular, on a one-dimensional cell $f$,
\begin{equation}
\label{eq:trimmed-one-dimensional-scalar}
  \Vcal_{r+1}^{-,0}(f)=\Pcal_{r+1}\Lambda^0(f).
\end{equation}
Indeed, for $\omega\in\Vcal_{r+1}^{-,0}(f)$,
\eqref{eq:trimmed-codim-one} gives
$\omega''\in\Pcal_{r-1}(f)$, where primes denote ordinary
derivatives along $f$.
Thus $\omega\in\Pcal_{r+1}(f)$, and the reverse inclusion
follows from \textup{(T1)}.
\end{remark}
\begin{table}[htbp]
\centering
\setlength{\belowcaptionskip}{5pt}
\caption{Vector proxies for $(\rd^k,\delta_k,\trc_k)$
in two and three dimensions.}
\label{tab:low-dimensional-vector-proxies}
\small
\renewcommand{\arraystretch}{1.12}
\setlength{\tabcolsep}{4pt}
\begin{tabular*}{\linewidth}{@{\extracolsep{\fill}}ccccc@{}}
\toprule
$d$ & $k=0$ & $k=1$ & $k=2$ & $k=3$ \\
\midrule
$2$ & $(\operatorname{grad},0,\gamma_0)$
    & $(\operatorname{rot},-\operatorname{div},\boldsymbol t\cdot)$
    & $(0,\operatorname{curl},0)$ & --- \\
$3$ & $(\operatorname{grad},0,\gamma_0)$
    & $(\operatorname{curl},-\operatorname{div},\Pi_\tau)$
    & $(\operatorname{div},\operatorname{curl},\boldsymbol n\cdot)$
    & $(0,-\operatorname{grad},0)$ \\
\bottomrule
\end{tabular*}
\par\smallskip
\begin{minipage}{\linewidth}
\footnotesize
Here $\boldsymbol n$ is the outward unit normal,
$\gamma_0$ denotes the scalar boundary trace, and
$\Pi_\tau\boldsymbol v=\boldsymbol v-(\boldsymbol v\cdot\boldsymbol n)\boldsymbol n$ is the tangential projection.
For $d=2$, we use $\boldsymbol t=(-n_2,n_1)$,
$\operatorname{rot}\boldsymbol v=\partial_1v_2-\partial_2v_1$, and
$\operatorname{curl}v=(\partial_2v,-\partial_1v)$.
\end{minipage}
\end{table}

\subsubsection{Geometric decomposition and dimension formulas}
\label{subsubsec:trimmed-geometric-decomposition}
Fix a form degree $k$ with $0\le k\le n$.
For each $f\in\Delta_d(\Mh)$ with $k\le d\le n$, we construct
a geometric decomposition of $\Vcal_{r+1}^{-,k}(f)$ by extending
forms with vanishing tangential trace from its subcells.
To this end, for $e\in\Delta_\ell(\Mh)$ with $k\le\ell\le n$, define
\[
  \ringV_{r+1}^{-,k}(e)
  :=\{\omega\in\Vcal_{r+1}^{-,k}(e):\trc\omega=0\}.
\]
The extensions are constructed through the local $\mathcal{NB}_k^*$ problems.
For barycentric extensions of simplicial polynomial forms, see
\cite[Section~4.7]{Arnold2006}.

For each $\ell$-cell $e$ and each $d$-cell $f$, we seek a linear lifting
\[
  \mathcal E_e^f:
  \ringV_{r+1}^{-,k}(e)
  \longrightarrow
  \Vcal_{r+1}^{-,k}(f).
\]
We first set $\mathcal E_e^e=\operatorname{id}$ and
$\mathcal E_e^f=0$ whenever $e\not\subseteq\overline f$.
These initial liftings satisfy, for every
$\omega_e\in\ringV_{r+1}^{-,k}(e)$, the trace compatibility identity
\begin{equation}
\label{eq:trimmed-extension-single-valued-traces}
  \trc_{e'}(\mathcal E_e^f\omega_e)
  =\mathcal E_e^{e'}\omega_e,
  \qquad
  e'\in\Delta_i(f),\quad k\le i\le d.
\end{equation}
Indeed, both sides vanish for the zero liftings, while the case $f=e$
follows from $\trc\omega_e=0$.

We construct the remaining liftings by induction on $d$, preserving
\eqref{eq:trimmed-extension-single-valued-traces}.
For the inductive step, let $d>\ell$ and $e\in\Delta_\ell(f)$.
Assume that the liftings to cells of dimension less than $d$ have been
constructed and satisfy \eqref{eq:trimmed-extension-single-valued-traces}.
Define the boundary data by
\[
  \mathfrak b|_{e_1}:=\mathcal E_e^{e_1}\omega_e,
  \qquad e_1\in\Delta_{d-1}(f).
\]
By the induction hypothesis \eqref{eq:trimmed-extension-single-valued-traces},
the traces match for any two facets
$e_1,e_2$ on every common subcell $e'$ of dimension at least $k$:
\begin{equation}
\label{eq:trimmed-extension-boundary-compatibility}
  \trc_{e'}^{e_1}\mathfrak b
  =\mathcal E_e^{e'}\omega_e
  =\trc_{e'}^{e_2}\mathfrak b
  \quad\Longrightarrow\quad
  \mathfrak b\in\Vcal_{r+1}^{-,k}(\partial f).
\end{equation}
With the compatible boundary data $\mathfrak b$ established in
\eqref{eq:trimmed-extension-boundary-compatibility} and zero interior data,
Theorem~\ref{thm:NBkstar-wellposed} defines a unique lifting
$\mathcal E_e^f\omega_e\in\Vcal_{r+1}^{-,k}(f)$ satisfying
\begin{equation}
\label{eq:trimmed-extension-local-problem}
  \delta\rd(\mathcal E_e^f\omega_e)=0,\quad
  \delta(\mathcal E_e^f\omega_e)=0\quad\text{in }f,\qquad
  \trc(\mathcal E_e^f\omega_e)=\mathfrak b.
\end{equation}
The equation involving $\delta$ alone is omitted for $k=0$.
For $e'\in\Delta_i(f)$ with $k\le i<d$, the boundary condition in
\eqref{eq:trimmed-extension-local-problem} and the compatibility established in
\eqref{eq:trimmed-extension-boundary-compatibility} give
\[
  \trc_{e'}(\mathcal E_e^f\omega_e)
  =\trc_{e'}\mathfrak b
  =\mathcal E_e^{e'}\omega_e.
\]
The case $e'=f$ is immediate, so
\eqref{eq:trimmed-extension-single-valued-traces} holds on $f$,
completing the induction.

For $e\in\Delta_\ell(f)$, identity
\eqref{eq:trimmed-extension-single-valued-traces}
and the zero-lifting convention imply
\begin{equation}
\label{eq:trimmed-entity-extension-traces}
  \trc_{e'}\mathcal E_e^f\omega_e
  =
  \begin{cases}
    \omega_e, & e'=e,\\
    0, & e'\ne e,
  \end{cases}
  \qquad
  e'\in\Delta_i(f),\quad k\le i\le\ell.
\end{equation}
In particular, $\trc_e\mathcal E_e^f=\operatorname{id}$,
so $\mathcal E_e^f$ is injective.

\begin{theorem}[Geometric decomposition]
\label{prop:trimmed-geometric-decomposition}
For every $f\in\Delta_d(\Mh)$ and $0\le k\le d$,
\begin{equation}
\label{eq:trimmed-local-geometric-decomposition}
  \Vcal_{r+1}^{-,k}(f)
  =\bigoplus_{\ell=k}^{d}\bigoplus_{e\in\Delta_\ell(f)}
    \mathcal E_e^f\ringV_{r+1}^{-,k}(e).
\end{equation}
Consequently, the following dimension formula holds:
\begin{equation}
\label{eq:trimmed-local-dimension}
  \dim\Vcal_{r+1}^{-,k}(f)
  =\sum_{\ell=k}^{d}\sum_{e\in\Delta_\ell(f)}
    \dim\ringV_{r+1}^{-,k}(e).
\end{equation}
\end{theorem}
\begin{proof}
Given $\omega\in\Vcal_{r+1}^{-,k}(f)$, we construct its geometric
decomposition by successively removing tangential traces in increasing
cell dimension.
Set $\zeta_{k-1}:=\omega$ and, for $\ell=k,\ldots,d$, define
\[
  \omega_e:=\trc_e\zeta_{\ell-1}
  \quad(e\in\Delta_\ell(f)),
  \qquad
  \zeta_\ell:=\zeta_{\ell-1}
  -\sum_{e\in\Delta_\ell(f)}\mathcal E_e^f\omega_e.
\]
We start with $\ell=k$.
For each $e\in\Delta_k(f)$, the trace $\omega_e=\trc_e\omega$
belongs to $\ringV_{r+1}^{-,k}(e)$, since its boundary trace vanishes by degree.
Thus $\zeta_k$ is well defined, and
\eqref{eq:trimmed-entity-extension-traces} gives
\[
  \trc_e\zeta_k
  =\trc_e\omega-\omega_e=0,
  \qquad e\in\Delta_k(f).
\]

Now let $k<\ell\le d$ and assume that
$\zeta_{\ell-1}\in\Vcal_{r+1}^{-,k}(f)$ has vanishing traces on subcells
of dimensions $k,\ldots,\ell-1$.
For each $e\in\Delta_\ell(f)$, property \textup{(T2)}, with
$\trc_f=\operatorname{id}$ when $e=f$, gives
$\omega_e\in\Vcal_{r+1}^{-,k}(e)$. Moreover,
\[
  \trc_{e'}\omega_e
  =\trc_{e'}\zeta_{\ell-1}
  =0,
  \qquad e'\in\Delta_{\ell-1}(e).
\]
Hence $\omega_e\in\ringV_{r+1}^{-,k}(e)$, so $\zeta_\ell$ is well defined.
For $e'\in\Delta_i(f)$ with $k\le i\le\ell$,
\eqref{eq:trimmed-entity-extension-traces} yields
\[
  \trc_{e'}\zeta_\ell
  =\trc_{e'}\zeta_{\ell-1}
   -\sum_{e\in\Delta_\ell(f)}
     \trc_{e'}(\mathcal E_e^f\omega_e)
  =
  \begin{cases}
    \trc_{e'}\zeta_{\ell-1}=0, & k\le i<\ell,\\
    \trc_{e'}\zeta_{\ell-1}-\omega_{e'}=0, & i=\ell.
  \end{cases}
\]
At $\ell=d$, we have $\omega_f=\zeta_{d-1}$ and
$\mathcal E_f^f=\operatorname{id}$, so $\zeta_d=0$.
Summing the residual identities gives
\begin{equation}
\label{eq:trimmed-geometric-expansion}
  \omega
  =\zeta_{k-1}-\zeta_d
  =\sum_{\ell=k}^{d}(\zeta_{\ell-1}-\zeta_\ell)
  =\sum_{\ell=k}^{d}\sum_{e\in\Delta_\ell(f)}
     \mathcal E_e^f\omega_e.
\end{equation}
Hence the spaces on the right-hand side of
\eqref{eq:trimmed-local-geometric-decomposition}
span $\Vcal_{r+1}^{-,k}(f)$.

For directness, suppose the final sum in
\eqref{eq:trimmed-geometric-expansion} vanishes for an arbitrary
family $\omega_e\in\ringV_{r+1}^{-,k}(e)$.
By \eqref{eq:trimmed-entity-extension-traces}, taking the trace on
each $e\in\Delta_k(f)$ gives
\[
  0
  =\trc_e\left(
    \sum_{\ell=k}^{d}\sum_{e'\in\Delta_\ell(f)}
    \mathcal E_{e'}^f\omega_{e'}
  \right)
  =\sum_{e'\in\Delta_k(f)}
    \trc_e(\mathcal E_{e'}^f\omega_{e'})
  =\omega_e.
\]
Taking traces successively on cells of dimensions $k+1,\ldots,d$,
each time using that the coefficients on lower-dimensional cells vanish,
gives $\omega_e=0$ for every $e$.
Thus the sum is direct.
The dimension formula \eqref{eq:trimmed-local-dimension}
follows from the injectivity of the liftings.
\end{proof}
\subsubsection{Unisolvent and computable DOFs}
\label{subsubsec:trimmed-dofs}
\label{subsubsec:trimmed-computable-moments}

We first give unisolvent DOFs for the subspaces with
vanishing tangential trace.
\begin{lemma}[Interior unisolvence]
\label{lem:trace-zero-unisolvence}
For $e\in\Delta_\ell(\Mh)$ and $0\le k\le\ell$, the
following moments form a unisolvent set
$\ringN_{r+1}^{-,k}(e)$ for $\ringV_{r+1}^{-,k}(e)$:
\begin{equation}
\label{eq:trace-zero-moments}
\begin{aligned}
&\int_e\omega\wedge\mu,
\quad
\mu\in\Pcal_r\Lambda^0(e),
&&
k=\ell;
\\
&\int_e\rd\omega\wedge\mu,
\quad
\mu\in\Kcal_r^0(e),
\qquad
\int_e\omega\wedge\nu,
\quad
\nu\in\Kcal_{r+1}^1(e),
&&
k=\ell-1;
\\
&\int_e\rd\omega\wedge\mu,
\quad
\mu\in\Kcal_{r+1}^{\ell-k-1}(e),
\qquad
\int_e\omega\wedge\nu,
\quad
\nu\in\Kcal_{r+1}^{\ell-k}(e),
&&
0\le k\le\ell-2.
\end{aligned}
\end{equation}
The moments involving $\nu$ are omitted when $k=0$.
\end{lemma}
\par
\begin{proof}
For $k=\ell$, the spaces
$\ringV_{r+1}^{-,\ell}(e)=\Pcal_r\Lambda^\ell(e)$
and $\Pcal_r\Lambda^0(e)$ have the same dimension.
If the moments vanish, choosing $\mu=\star\omega$ gives
$\|\omega\|_e^2=0$, proving unisolvence.
\par

For $0\le k<\ell$, Theorem~\ref{thm:NBkstar-wellposed}
with zero boundary trace, together with the injectivity
of $\rd$ on the Koszul spaces recalled in
Section~\ref{subsec:polynomial-forms}, gives
\begin{equation}
\label{eq:trimmed-trace-zero-dimension}
  \dim\ringV_{r+1}^{-,k}(e)
  =
  \begin{cases}
    \dim\Kcal_r^0(e)+\dim\Kcal_{r+1}^1(e),
    & k=\ell-1,\\[1mm]
    \dim\Kcal_{r+1}^{\ell-k-1}(e)
    +\dim\Kcal_{r+1}^{\ell-k}(e),
    & 0\le k\le\ell-2.
  \end{cases}
\end{equation}
By \eqref{eq:trace-zero-moments}, the right-hand side
equals $\#\ringN_{r+1}^{-,k}(e)$.
It therefore remains to show that vanishing moments
imply $\omega=0$.
Let $\omega\in\ringV_{r+1}^{-,k}(e)$ have vanishing moments.
By \eqref{eq:trimmed-codim-one} and \eqref{eq:trimmed-general},
choose $\mu_0$ in the derivative test space in
\eqref{eq:trace-zero-moments} such that
$\delta\rd\omega=\star^{-1}\rd\mu_0$.
Since $\trc\omega=0$, the integration-by-parts formula
\eqref{eq:green-formula} and the vanishing derivative
moments yield
\begin{equation*}
  \|\rd\omega\|_e^2
  =
  \langle\omega,\delta\rd\omega\rangle_e
  =\int_e\omega\wedge\rd\mu_0
  =
  (-1)^{k+1}\int_e\rd\omega\wedge\mu_0
  =
  0.
\end{equation*}
Hence $\rd\omega=0$.
If $k=0$, the zero trace and $\rd\omega=0$ imply $\omega=0$.
Since $e$ is contractible, the homogeneous Hodge decomposition
\eqref{eq:NBkstar-hodge-decomposition} yields $\omega=\rd\phi$ for some
$\phi\in\overset{\circ}{H}\Lambda^{k-1}(\rd;e)$ whenever $1\le k<\ell$.
The second interior condition gives
$\delta\omega=\star^{-1}\rd\nu_0$ for some
$\nu_0\in\Kcal_{r+1}^{\ell-k}(e)$.
Using $\trc\phi=0$, \eqref{eq:green-formula},
and the vanishing value moments, we obtain
\[
  \|\omega\|_e^2
  =
  \langle\phi,\delta\omega\rangle_e
  =
  (-1)^k\int_e\omega\wedge\nu_0
  =
  0.
\]
This proves unisolvence.
\end{proof}
For $f\in\Delta_d(\Mh)$ and $0\le k\le d$, we define the local
DOFs for $\Vcal_{r+1}^{-,k}(f)$ by applying the moments
\eqref{eq:trace-zero-moments} to the tangential traces on its subcells.
Splitting the $k$-cell moments by
\eqref{eq:polynomial-koszul-decomposition-k0} gives the following
three classes, whose union is denoted by $\Ncal_{r+1}^{-,k}(f)$:
\begin{align}
\mathsf C^{k,-}(f):\quad&
  \int_e\trc_e\omega,\qquad \ell=k,
\label{eq:trimmed-dof-constant}\\[1mm]
\mathsf K_\ell^{k,-}(f):\quad&
  \int_e\trc_e\omega\wedge\nu,\qquad
  \nu\in
  \begin{cases}
    \Kcal_r^0(e), & \ell=k,\\
    \Kcal_{r+1}^{\ell-k}(e), & k+1\le\ell\le d,
  \end{cases}
\label{eq:trimmed-dof-koszul}\\[1mm]
\mathsf D_\ell^{k,-}(f):\quad&
  \int_e\trc_e\rd\omega\wedge\mu,\qquad
  \mu\in
  \begin{cases}
    \Kcal_r^0(e), & \ell=k+1,\\
    \Kcal_{r+1}^{\ell-k-1}(e), & k+2\le\ell\le d.
  \end{cases}
\label{eq:trimmed-dof-derivative}
\end{align}
Here $e\in\Delta_\ell(f)$, and the moments are well defined
by \textup{(T1)}--\textup{(T2)}.

\begin{theorem}[Trimmed unisolvence]
\label{thm:trimmed-local-unisolvence}
The family $\Ncal_{r+1}^{-,k}(f)$ is unisolvent for
$\Vcal_{r+1}^{-,k}(f)$.
\end{theorem}
\begin{proof}
From
\eqref{eq:trimmed-local-dimension},
Lemma~\ref{lem:trace-zero-unisolvence},
and the definition of
$\Ncal_{r+1}^{-,k}(f)$, we have
\[
\begin{aligned}
  \dim\Vcal_{r+1}^{-,k}(f)
  &=
  \sum_{\ell=k}^{d}
  \sum_{e\in\Delta_\ell(f)}
  \dim\ringV_{r+1}^{-,k}(e)
  =
  \sum_{\ell=k}^{d}
  \sum_{e\in\Delta_\ell(f)}
  \#\ringN_{r+1}^{-,k}(e)
  =
  \#\Ncal_{r+1}^{-,k}(f).
\end{aligned}
\]
Let $\omega\in\Vcal_{r+1}^{-,k}(f)$ have vanishing moments.
For $e\in\Delta_k(f)$,
$\trc_e\omega\in\Vcal_{r+1}^{-,k}(e)=\ringV_{r+1}^{-,k}(e)$,
and all its local moments vanish, so
Lemma~\ref{lem:trace-zero-unisolvence} gives $\trc_e\omega=0$.
Using \textup{(T2)}, apply the same lemma successively on cells
of dimensions $k+1,\ldots,d$, whose facet traces have already vanished.
Thus all cell traces are zero, and taking $e=f$ gives $\omega=0$.
Together with the dimension identity, this proves unisolvence.
\end{proof}
\begin{remark}[Equivalent forms of DOFs]
\label{rem:equivalent-dofs}
\label{rem:trimmed-moment-shift}
For $\omega\in\ringV_{r+1}^{-,k}(e)$,
the integration-by-parts formula \eqref{eq:stokes-integration-by-parts} gives
\begin{equation}
\label{eq:trace-zero-moment-alternative}
  \int_e\rd\omega\wedge\mu
  =(-1)^{k+1}\int_e\omega\wedge\rd\mu.
\end{equation}
An alternative unisolvent set on $\Vcal_{r+1}^{-,k}(f)$ is
obtained by replacing $\mathsf D_\ell^{k,-}(f)$ with the value moments
\begin{equation}
\label{eq:trimmed-dof-derivative-equivalent}
  \int_e\trc_e\omega\wedge\rd\mu,
  \qquad e\in\Delta_\ell(f),\quad k+1\le\ell\le d,
\end{equation}
where $\mu$ ranges over the same test spaces as in
\eqref{eq:trimmed-dof-derivative}.
Unisolvence follows by the same induction as in
Theorem~\ref{thm:trimmed-local-unisolvence}, since
\eqref{eq:trace-zero-moment-alternative} recovers the original derivative
moments on each subcell once its boundary trace vanishes.
These alternative moments avoid the additional regularity needed
to take traces of $\rd\omega$ in canonical interpolation.
We retain the derivative moments because the
$\mathsf D_\ell^{k,-}(f)$ moments of $\omega$ are precisely the
$\mathsf K_\ell^{k+1,-}(f)$ moments of $\rd\omega$,
making the commuting property transparent.
\end{remark}

\begin{proposition}[Recursive computability of moments]
\label{prop:trimmed-computable-moments}
Let $f\in\Delta_d(\Mh)$, $0\le k\le d$, and
$\omega\in\Vcal_{r+1}^{-,k}(f)$.
The moments
\begin{equation}
\label{eq:trimmed-computable-moments}
  \int_f\omega\wedge\eta,
  \qquad \eta\in\Pcal_{r+1}^{-}\Lambda^{d-k}(f),
\end{equation}
are computable recursively in $d$ from $\Ncal_{r+1}^{-,k}(f)$.
Consequently, the complete degree-$r$ projection
$\Proj{r}{f}{k}\omega$ is also computable.
\end{proposition}
\begin{proof}
We argue by induction on $d$.
For $d=k$, \eqref{eq:trimmed-top} gives
$\omega\in\Pcal_r\Lambda^d(f)$, which is determined by the DOFs,
so all these moments are computable.
For $d>k$, assume that the required moments are computable on
lower-dimensional cells. The isomorphism
\eqref{eq:trimmed-koszul-isomorphism} gives $\eta=\rd\mu+\nu$ with
$\mu\in\Kcal_{r+1}^{d-k-1}(f)$ and
$\nu\in\Kcal_{r+1}^{d-k}(f)$.
The integration-by-parts formula
\eqref{eq:stokes-integration-by-parts} yields
\[
  (-1)^{k+1}\int_f\omega\wedge\eta
  =\int_f\rd\omega\wedge\mu
    +(-1)^{k+1}\int_f\omega\wedge\nu
    -\int_{\partial f}\trc\omega\wedge\trc\mu.
\]
The first volume term is a
$\mathsf D_d^{k,-}(f)$ moment for $d\ge k+2$.
For $d=k+1$, \textup{(T3)} and \eqref{eq:trimmed-top}
give $\rd\omega\in\Pcal_r\Lambda^d(f)$.
By \eqref{eq:polynomial-koszul-decomposition-k0},
this polynomial is determined by the
$\mathsf D_d^{k,-}(f)$ moments and the constant moment
$\int_f\rd\omega=\int_{\partial f}\trc\omega$
supplied by the $\mathsf C^{k,-}(f)$ data.
Thus $\int_f\rd\omega\wedge\mu$ is computable
also for $\mu\in\Kcal_{r+1}^0(f)$.
The second volume term is a $\mathsf K_d^{k,-}(f)$ moment.
For the boundary term, Lemma~\ref{lem:traces-trimmed} gives
$\trc_e\mu\in\Pcal_{r+1}^-\Lambda^{d-k-1}(e)$
on each facet $e$ of $f$, and \eqref{eq:integral.relative.orientation} gives
\[
  \int_{\partial f}\trc\omega\wedge\trc\mu
  =\sum_{e\in\Delta_{d-1}(f)}\varepsilon_{fe}
    \int_e\trc_e\omega\wedge\trc_e\mu.
\]
By trace closure \textup{(T2)},
$\trc_e\omega\in\Vcal_{r+1}^{-,k}(e)$.
Each facet moment is therefore computable by the induction hypothesis,
since the facet DOFs are included in $\Ncal_{r+1}^{-,k}(f)$.
\par
Note that, by \eqref{eq:hodge-projection-moments}, the computed moments
in \eqref{eq:trimmed-computable-moments}
determine the trimmed projection $\Pi_{r+1,f}^{-,d-k}(\star\omega)$.
Since $\star\Pcal_r\Lambda^k(f)=\Pcal_r\Lambda^{d-k}(f)
\subset\Pcal_{r+1}^-\Lambda^{d-k}(f)$,
these moments also determine $\Proj{r}{f}{k}\omega$
by \eqref{eq:l2-orthogonal-projection}.
\end{proof}
\subsection{Complete virtual element spaces of differential forms}
\label{subsec:complete-spaces}
The trimmed spaces contain higher-order polynomial forms and provide,
through their DOFs, the moment information needed to compute
$\Proj{r}{f}{k}\omega$.
In this section, we follow the VEM enhancement strategy
\cite{AhmadEtAl2013} to reduce these DOFs while retaining this
computability.
\subsubsection{Construction by moment constraints}
\label{subsubsec:complete-construction}

For $f\in\Delta_d(\Mh)$ and
$0\le k\le d$, define
\begin{equation}
\label{eq:complete-space}
\begin{aligned}
  \Vcal_r^k(f):=
  \biggl\{
  \omega\in\Vcal_{r+1}^{-,k}(f):
  {}\;&
  \int_e
  \bigl(
  \trc_e\omega
  -
  \Proj{r}{e}{k}\trc_e\omega
  \bigr)
  \wedge\nu
  =0,
  \\
  &\quad
  \nu\in
  \kappa\mathcal H_r\Lambda^{\ell-k+1}(e),
  \quad
  e\in\Delta_\ell(f),
  \quad
  k+1\le\ell\le d;
  \\
  &\int_e
  \bigl(
  \trc_e\rd\omega
  -
  \Proj{r}{e}{k+1}
  \trc_e\rd\omega
  \bigr)
  \wedge\mu
  =0,
  \\
  &\quad
  \mu\in
  \kappa\mathcal H_r\Lambda^{\ell-k}(e),
  \quad
  e\in\Delta_\ell(f),
  \quad
  k+2\le\ell\le d
  \biggr\},
\end{aligned}
\end{equation}
Empty index ranges impose no conditions.
In particular, $\Vcal_r^d(f)=\Pcal_r\Lambda^d(f)$;
for $k=d-1$, only the first moment constraint in
\eqref{eq:complete-space}, with $e=f$, remains.
\par

\begin{proposition}[Local complete-space
properties]
\label{prop:complete-local-properties}
\leavevmode\par
For each $f\in\Delta_d(\Mh)$, the family
$\{\Vcal_r^k(f)\}_{k=0}^d$ satisfies the
following properties.
\begin{itemize}
\item[{\rm (C1)}]
Conformity and polynomial inclusion:
$\Vcal_r^k(f)\subset H\Lambda^k(\rd;f)$,
$\Pcal_r\Lambda^k(f)\subset\Vcal_r^k(f)$
for $1\le k\le d$, and
$\Pcal_{r+1}\Lambda^0(f)\subset
\Vcal_r^0(f)$.
\item[{\rm (C2)}]
Trace closure: if
$e\in\Delta_\ell(f)$ and
$0\le k\le\ell<d$, then  
$
\trc_e\Vcal_r^k(f)
\subset\Vcal_r^k(e).
$
\item[{\rm (C3)}]
Complex property:
$
  \rd\Vcal_r^k(f)
  \subset
  \Vcal_r^{k+1}(f),
  \qquad
  0\le k<d.
$
\end{itemize}
\end{proposition}
\par
\begin{proof}
The inclusion $\Vcal_r^k(f)\subset\Vcal_{r+1}^{-,k}(f)$
and \textup{(T1)} give conformity.
For $1\le k\le d$, \textup{(T1)} ensures
$\Pcal_r\Lambda^k(f)\subset\Vcal_{r+1}^{-,k}(f)$,
and both residuals in \eqref{eq:complete-space}
vanish on $\Pcal_r\Lambda^k(f)$.
For $k=0$, \textup{(T1)} gives
$\Pcal_{r+1}\Lambda^0(f)\subset\Vcal_{r+1}^{-,0}(f)$.
The first moment constraints are void because
$\Lambda^{\ell+1}(e)=\{0\}$ for $e\in\Delta_\ell(f)$,
and the second residual vanishes since
$\rd\Pcal_{r+1}\Lambda^0(f)\subset\Pcal_r\Lambda^1(f)$.
Thus \textup{(C1)} holds.
\par

Let $0\le k<d$ and $\omega\in\Vcal_r^k(f)$.
For $e\in\Delta_\ell(f)$ with $k\le\ell<d$,
\textup{(T2)} gives
$\trc_e\omega\in\Vcal_{r+1}^{-,k}(e)$.
Since traces commute with $\rd$, the moment constraints
for $\trc_e\omega$ on the subcells of $e$ follow from
\eqref{eq:complete-space}. Hence
$\trc_e\omega\in\Vcal_r^k(e)$, proving \textup{(C2)}.
\par
By \textup{(T3)}, $\rd\omega\in\Vcal_{r+1}^{-,k+1}(f)$.
Its first moment constraints are exactly the second
moment constraints for $\omega$ in \eqref{eq:complete-space},
and its second moment constraints hold because
$\rd^2\omega=0$.
Thus $\rd\omega\in\Vcal_r^{k+1}(f)$, proving \textup{(C3)}.
\end{proof}
From the Koszul splitting \eqref{eq:koszul-degree-splitting}, the constraints \eqref{eq:complete-space} reduce the trimmed moment classes
\eqref{eq:trimmed-dof-constant}--\eqref{eq:trimmed-dof-derivative}
to the following retained data, for $e\in\Delta_\ell(f)$:
\begin{equation}
\label{eq:complete-retained-moments}
\begin{aligned}
 &\int_e\trc_e\omega\wedge(\nu +\mu),
 &&(\nu,\mu)\in \Kcal_r^0(e)\times \Pcal_0\Lambda^0(e),
 &&\ell=k,\\
 &\int_e\trc_e\omega\wedge\nu,\quad
   \int_e\trc_e\rd\omega\wedge\mu,
 &&(\nu,\mu)\in\Kcal_r^{\ell-k}(e)\times\Kcal_r^{\ell-k-1}(e),
 &&k+1\le\ell\le d.
\end{aligned}
\end{equation}
For $\ell=k$, the scalar decomposition
\eqref{eq:polynomial-koszul-decomposition-k0} and
\eqref{eq:trimmed-polynomial-top-degree} combine the constant
and Koszul test functions into $\Pcal_r^-\Lambda^0(e)$.
For $\ell>k$, passing to the equivalent value moments
\eqref{eq:trimmed-dof-derivative-equivalent} gives the combined
test forms $\rd\mu+\nu$, which range over
$\Pcal_r^-\Lambda^{\ell-k}(e)$ by the isomorphism
\eqref{eq:trimmed-koszul-isomorphism}.
These observations motivate the following unified DOF set
$\Ncal_r^k(f)$:
\begin{equation}
\label{eq:complete-dofs}
  \int_e\trc_e\omega\wedge\mu,
  \qquad \mu\in\Pcal_r^-\Lambda^{\ell-k}(e),
  \quad e\in\Delta_\ell(f),\quad k\le\ell\le d.
\end{equation}
We prove below that these moments suffice to compute
the $L^2$ projections and form a unisolvent set.

\subsubsection{Recursive computation of \(L^2\) projections}
\label{subsubsec:complete-computable-projections}
\begin{proposition}[Computability of $L^2$ projections]
\label{thm:projection-computable}
Let $f\in\Delta_d(\Mh)$, $0\le k\le d$, and
$\omega\in\Vcal_r^k(f)$.
The $L^2$ projections $\Proj{r}{f}{k}\omega$ and,
for $k<d$, $\Proj{r}{f}{k+1}\rd\omega$ are computable
from the DOFs $\Ncal_r^k(f)$ by recursion on the cell dimension $d$.
\end{proposition}
\par
\begin{proof}
By \eqref{eq:l2-orthogonal-projection} and
$\star\Pcal_r\Lambda^\ell(f)=\Pcal_r\Lambda^{d-\ell}(f)$,
it suffices to compute recursively in $d$ the following wedge moments:
\[
  \int_f\rd\omega\wedge\eta,\quad
  \eta\in\Pcal_r\Lambda^{d-k-1}(f)\ (k<d),
  \qquad
  \int_f\omega\wedge\eta,\quad
  \eta\in\Pcal_r\Lambda^{d-k}(f).
\]
For $d=k$, \eqref{eq:complete-space} and
\eqref{eq:trimmed-top} give
$\Vcal_r^k(f)=\Pcal_r\Lambda^k(f)$.
Since $\Pcal_r^-\Lambda^0(f)=\Pcal_r\Lambda^0(f)$,
the moments \eqref{eq:complete-dofs} determine $\omega$,
so the required moments are computable.
\par
For $d>k$, assume that the required moments are computable on
lower-dimensional cells.
Write $\Pi_{r,\partial f}^k$ for the facetwise $L^2$ projector,
defined by
\[
  \left(\Pi_{r,\partial f}^k\trc\omega\right)|_e
  :=\Proj{r}{e}{k}\trc_e\omega,
  \qquad e\in\Delta_{d-1}(f).
\]
Boundary integrals are understood as oriented facet sums
according to \eqref{eq:integral.relative.orientation}.
We compute the two families of moments in the following order.
\begin{itemize}
\item \emph{Moments of $\rd\omega$.}
For $\eta\in\Pcal_r\Lambda^{d-k-1}(f)$,
the integration-by-parts formula \eqref{eq:stokes-integration-by-parts} gives
\begin{equation}
\label{eq:projection-computability-derivative}
  \int_f\rd\omega\wedge\eta
  =(-1)^{k+1}\int_f\omega\wedge\rd\eta
    +\int_{\partial f}\trc\omega\wedge\trc\eta.
\end{equation}
The volume term is determined by \eqref{eq:complete-dofs}, since
$\rd\eta\in\Pcal_{r-1}\Lambda^{d-k}(f)
\subset\Pcal_r^-\Lambda^{d-k}(f)$.
On each facet $e$ of $f$, trace closure \textup{(C2)} gives
$\trc_e\omega\in\Vcal_r^k(e)$.
Since $\trc_e\eta\in\Pcal_r\Lambda^{d-k-1}(e)$,
the boundary term depends only on
$\Pi_{r,\partial f}^k\trc\omega$, which is computable
by the induction hypothesis.
\item \emph{Moments of $\omega$.}
For $\eta\in\Pcal_r\Lambda^{d-k}(f)$, the isomorphism
\eqref{eq:polynomial-koszul-isomorphism} gives
$\eta=\rd\mu+\nu$, with
$\mu\in\Kcal_{r+1}^{d-k-1}(f)$ and
$\nu\in\Kcal_r^{d-k}(f)$.
Integration by parts gives
\begin{equation}
\label{eq:projection-computability-stokes}
\begin{aligned}
  (-1)^{k+1}\int_f\omega\wedge\eta
  &=\int_f\rd\omega\wedge\mu
    +(-1)^{k+1}\int_f\omega\wedge\nu
    -\int_{\partial f}\trc\omega\wedge\trc\mu\\
  &=\int_f\Proj{r}{f}{k+1}\rd\omega\wedge\mu
    +(-1)^{k+1}\int_f\omega\wedge\nu
    -\int_{\partial f}
       \Pi_{r,\partial f}^k\trc\omega\wedge\trc\mu.
\end{aligned}
\end{equation}
For the second equality, the projections
$\Proj{r}{f}{k+1}\rd\omega$ and
$\Pi_{r,\partial f}^k\trc\omega$
have already been computed using
\eqref{eq:projection-computability-derivative}
and the induction hypothesis, respectively.
For $d\ge k+2$, the isomorphism
\eqref{eq:trimmed-koszul-isomorphism} and the decomposition
\eqref{eq:trimmed-polynomial-space} give
\begin{equation*}
  \Kcal_{r+1}^{d-k-1}(f)
  \subset\Pcal_{r+1}^-\Lambda^{d-k-1}(f)
  =\Pcal_r\Lambda^{d-k-1}(f)
   \oplus\kappa\mathcal H_r\Lambda^{d-k}(f).
\end{equation*}
By Lemma~\ref{lem:traces-trimmed},
$\trc_e\mu\in\Pcal_{r+1}^-\Lambda^{d-k-1}(e)$
on each facet $e$.
Using this splitting on $f$ and its facets, the definitions
of the $L^2$ projections and the two constraints in
\eqref{eq:complete-space} give the replacements in both
the derivative volume term and the boundary term.
When $d=k+1$, $\rd\omega$ and the facet traces of $\omega$
already belong to the corresponding degree-$r$ polynomial
spaces by \textup{(C2)}--\textup{(C3)} and
\eqref{eq:trimmed-top}, so both replacements hold directly.
The remaining volume moment is one of the prescribed DOFs, since
$\nu\in\Kcal_r^{d-k}(f)\subset\Pcal_r^-\Lambda^{d-k}(f)$.
\end{itemize}
Thus both families of moments are computable, completing the induction.
\end{proof}
Under the three-dimensional vector proxies in Table~\ref{tab:low-dimensional-vector-proxies},
our complete spaces $\Vcal_r^k(f)$, $k=0,1,2$, and the recursive
$L^2$-projection computations in Proposition~\ref{thm:projection-computable}
have nodal, edge, and face counterparts
in~\cite[Sections~6.2.1--6.2.2]{BeiraoDassiDiPietroDroniou2022}.
In particular, the degree-$(r+1)$ nodal $L^2$ projection computed there
is recovered by the scalar modification in
Remark~\ref{rem:degree-r-plus-one-scalar-potential} of
Appendix~\ref{subsec:appendix-fully-discrete-realizations}, using the same moment data.
\subsubsection{Unisolvent DOFs}
\label{subsubsec:complete-dofs}
\begin{theorem}[Complete unisolvence]
\label{thm:complete-local-unisolvence}
For every $f\in\Delta_d(\Mh)$ and
$0\le k\le d$, the family $\Ncal_r^k(f)$ is unisolvent for
$\Vcal_r^k(f)$.
\end{theorem}
\par
\begin{proof}
For $k=d$, no additional restrictions are imposed, so
$\Vcal_r^d(f)$ coincides with $\Vcal_{r+1}^{-,d}(f)$
and has the same DOFs.
Unisolvence therefore follows from
Theorem~\ref{thm:trimmed-local-unisolvence}. Hence assume
$0\le k\le d-1$.
\par

Define the product of the test spaces for the higher-order
moment constraints in \eqref{eq:complete-space} by
\begin{equation}
\label{eq:complete-constraint-test-space}
  \mathcal R_r^k(f)
  :=
  \bigtimes_{e\in\Delta_{k+1}(f)}
  \kappa\mathcal H_r\Lambda^2(e)
  \times
  \bigtimes_{\ell=k+2}^{d}
  \bigtimes_{e\in\Delta_\ell(f)}
  \left(
  \kappa\mathcal H_r\Lambda^{\ell-k}(e)
  \times
  \kappa\mathcal H_r\Lambda^{\ell-k+1}(e)
  \right)
\end{equation}
The constraints \eqref{eq:complete-space} impose at most
$\dim\mathcal R_r^k(f)$ independent linear conditions on
$\Vcal_{r+1}^{-,k}(f)$.
Hence,
\[
\begin{aligned}
  \dim\Vcal_r^k(f)
  &\ge
  \dim\Vcal_{r+1}^{-,k}(f)
  -\dim\mathcal R_r^k(f)\\
  &=\sum_{e\in\Delta_k(f)}\dim\Pcal_r\Lambda^0(e)
    +\sum_{\ell=k+1}^{d}\sum_{e\in\Delta_\ell(f)}
      \bigl(\dim\Kcal_r^{\ell-k-1}(e)
      +\dim\Kcal_r^{\ell-k}(e)\bigr)\\
  &=
  \sum_{\ell=k}^d
  \sum_{e\in\Delta_\ell(f)}
  \dim\Pcal_r^-
  \Lambda^{\ell-k}(e).
\end{aligned}
\]
The first equality follows from
\eqref{eq:trimmed-local-dimension},
\eqref{eq:trimmed-trace-zero-dimension}, and the Koszul splitting
\eqref{eq:koszul-degree-splitting}.
The second equality follows from
\eqref{eq:trimmed-polynomial-top-degree} for $\ell=k$ and
the isomorphism \eqref{eq:trimmed-koszul-isomorphism} for $\ell>k$.
Thus, $\dim\Vcal_r^k(f)$ is at least the number of prescribed DOFs.
\par
In view of the dimension bound, it remains to prove that
vanishing DOFs imply $\omega=0$.
Let $\omega\in\Vcal_r^k(f)$ have vanishing DOFs in
$\Ncal_r^k(f)$.
Since $\Vcal_r^k(f)\subset\Vcal_{r+1}^{-,k}(f)$,
it suffices by Theorem~\ref{thm:trimmed-local-unisolvence}
to show that every DOF in $\Ncal_{r+1}^{-,k}(f)$
vanishes on $\omega$.
The Koszul splitting \eqref{eq:koszul-degree-splitting} decomposes their test spaces
into the retained spaces in
\eqref{eq:complete-retained-moments} and the higher-order
spaces in \eqref{eq:complete-constraint-test-space}.
The retained moments vanish by trace closure \textup{(C2)}
and Proposition~\ref{thm:projection-computable},
while the higher-order moments vanish by the constraints
\eqref{eq:complete-space}.
Thus these DOFs vanish on $\omega$, giving $\omega=0$.
\end{proof}

\begin{remark}[Lowest-order spaces and Whitney forms]
\label{rem:lowest-order-whitney}
At $r=0$, the complete spaces satisfy
\begin{equation}
\label{eq:lowest-order-whitney-inclusion}
  \Pcal_1^-\Lambda^k(f)\subset\Vcal_0^k(f),
  \qquad f\in\Delta_d(\Mh),\quad 0\le k\le d.
\end{equation}
Indeed, for $\eta\in\Pcal_1^-\Lambda^k(f)$,
Lemma~\ref{lem:traces-trimmed} gives
$\trc_e \eta=\mu_0+\kappa \nu_0$ on each subcell $e\in\Delta_\ell(f)$
with $\ell\ge k$, where $\mu_0$ and $\nu_0$ are constant forms
of degrees $k$ and $k+1$.
Since $\nu_0$ is constant, each coefficient of $\kappa \nu_0$ is a
linear combination of the components of $x-x_e$.
These have zero mean on $e$ because $x_e$ is its centroid.
Thus $\Proj{0}{e}{k}\kappa \nu_0=0$ and
$\Proj{0}{e}{k}\trc_e \eta=\mu_0$.
For $k+1\le\ell\le d$ and a test form
$\nu\in\kappa\mathcal H_0\Lambda^{\ell-k+1}(e)$ as in
\eqref{eq:complete-space}, we have $\kappa\nu=0$ by $\kappa^2=0$
and $\nu_0\wedge\nu=0$ because its form degree is $\ell+1$.
The graded Leibniz rule \eqref{eq:contraction-wedge-rule}
therefore gives the first constraints directly:
\[
  \bigl(\trc_e \eta-\Proj{0}{e}{k}\trc_e \eta\bigr)\wedge\nu
  =(\kappa \nu_0)\wedge\nu
  =\kappa(\nu_0\wedge\nu)-(-1)^{k+1}\nu_0\wedge\kappa\nu=0.
\]
The second constraints hold because $\rd\trc_e \eta$ is constant.
On a simplex, the DOFs \eqref{eq:complete-dofs} are precisely the Whitney DOFs
\cite[Section~4.3]{Arnold2006}. Thus
Theorem~\ref{thm:complete-local-unisolvence} and
\eqref{eq:lowest-order-whitney-inclusion} give
\[
  \Vcal_0^k(f)=\Pcal_1^-\Lambda^k(f)
  \qquad\text{if $f$ is a simplex}.
\]
\end{remark}

\subsection{Global complexes, commuting interpolation, and de Rham cohomology}
\label{subsec:fully-discrete-cochain-isomorphisms}

\subsubsection{Global assembly}
\label{subsubsec:global-conforming-assembly}

For \(0\le k\le n\), define the global trimmed space
\begin{equation}
\label{eq:global-trimmed-space}
  \Vcal_{r+1,h}^{-,k}
  :=
  \left\{
    \omega_h\in H\Lambda^k(\rd;\Omega):
    \omega_h|_T\in\Vcal_{r+1}^{-,k}(T)
    \quad\forall T\in\Th
  \right\}.
\end{equation}
Its global DOF maps
$\mathsf C_h^{k,-}$, $\mathsf K_{\ell,h}^{k,-}$, and
$\mathsf D_{\ell,h}^{k,-}$ assign to each
$\omega_h\in\Vcal_{r+1,h}^{-,k}$ the moments
\eqref{eq:trimmed-dof-constant}--\eqref{eq:trimmed-dof-derivative},
with $e\in\Delta_\ell(\Mh)$ and $d$ replaced by $n$.
Theorem~\ref{thm:trimmed-local-unisolvence} determines a unique
$\omega_h|_T\in\Vcal_{r+1}^{-,k}(T)$ on each $T\in\Th$
from the prescribed moments.
For $k<n$ and any two mesh elements $T_1,T_2\in\Th$ sharing an interior facet
$f\in\Delta_{n-1}(T_1)\cap\Delta_{n-1}(T_2)$, both traces
$\trc_f^{T_i}\omega_h$, $i=1,2$, belong to
$\Vcal_{r+1}^{-,k}(f)$ by \textup{(T2)}.
They have the same local moments and hence coincide by local unisolvence.
Lemma~\ref{lem:trace-gluing} then gives
$\omega_h\in H\Lambda^k(\rd;\Omega)$.
For $k=n$, the space consists of piecewise polynomial $n$-forms
with no interelement continuity conditions.
Thus local unisolvence gives global unisolvence in every degree,
and \textup{(T3)} yields the conforming trimmed complex
\begin{equation}
\label{eq:trimmed-virtual-derham-complex}
  0
  \longrightarrow
  \Vcal_{r+1,h}^{-,0}
  \xrightarrow{\rd}
  \Vcal_{r+1,h}^{-,1}
  \xrightarrow{\rd}
  \cdots
  \xrightarrow{\rd}
  \Vcal_{r+1,h}^{-,n}
  \longrightarrow0.
\end{equation}
Similarly, define the global complete space
\begin{equation}
\label{eq:global-complete-space}
  \Vcal_{r,h}^{k}
  :=
  \left\{
    \omega_h\in H\Lambda^k(\rd;\Omega):
    \omega_h|_T\in\Vcal_r^k(T)
    \quad\forall T\in\Th
  \right\},
\end{equation}
and equip it with the moment DOFs
\begin{equation}
\label{eq:global-complete-dofs}
  \int_f\trc_f\omega_h\wedge\mu,
  \qquad
  \mu\in\Pcal_r^-\Lambda^{d-k}(f),
  \quad
  f\in\Delta_d(\Mh),
  \quad
  k\le d\le n.
\end{equation}
The same assembly argument, using
Theorem~\ref{thm:complete-local-unisolvence} and properties
\textup{(C1)}--\textup{(C2)}, proves global unisolvence of
\eqref{eq:global-complete-dofs}.
Property \textup{(C3)} yields the conforming complete complex
\begin{equation}
\label{eq:complete-virtual-derham-complex}
  0
  \longrightarrow
  \Vcal_{r,h}^{0}
  \xrightarrow{\rd}
  \Vcal_{r,h}^{1}
  \xrightarrow{\rd}
  \cdots
  \xrightarrow{\rd}
  \Vcal_{r,h}^{n}
  \longrightarrow0.
\end{equation}
\begin{remark}[Fully discrete realizations]
\label{rem:fully-discrete-identification}
For $\omega_h^-\in\Vcal_{r+1,h}^{-,k}$, define the trimmed
DOF map by its components on each
$f\in\Delta_d(\Mh)$, $d\ge k$:
\begin{equation}
\label{eq:global-fully-discrete-isomorphisms}
  \bigl(\mathcal J_{r+1,h}^{-,k}\omega_h^-\bigr)_f
  :=\begin{cases}
    \Pi_{r,f}^{0}(\star_f\trc_f\omega_h^-),&d=k,\\[2pt]
    \bigl(\Pi_{r+1,f}^{\kappa,1}(\star_f\trc_f\omega_h^-),
          \Pi_{r,f}^{\kappa,0}(\star_f\trc_f\rd\omega_h^-)\bigr),
       &d=k+1,\\[2pt]
    \bigl(\Pi_{r+1,f}^{\kappa,d-k}(\star_f\trc_f\omega_h^-),
          \Pi_{r+1,f}^{\kappa,d-k-1}(\star_f\trc_f\rd\omega_h^-)\bigr),
       &d\ge k+2.
  \end{cases}
\end{equation}
For $\omega_h\in\Vcal_{r,h}^{k}$, define the complete
DOF map by
\begin{equation}
\label{eq:global-complete-isomorphism}
  \bigl(\mathcal J_{r,h}^{k}\omega_h\bigr)_f
  :=\Pi_{r,f}^{-,d-k}(\star_f\trc_f\omega_h),
  \qquad f\in\Delta_d(\Mh),\quad d\ge k.
\end{equation}
Then define the fully discrete spaces by
\begin{align}
\label{eq:global-trimmed-fully-discrete-space}
  \underline{\Vcal}_{r+1,h}^{-,k}
  &:=\mathcal J_{r+1,h}^{-,k}\bigl(\Vcal_{r+1,h}^{-,k}\bigr),\\
\label{eq:global-complete-fully-discrete-space}
  \underline{\Vcal}_{r,h}^{k}
  &:=\mathcal J_{r,h}^{k}\bigl(\Vcal_{r,h}^{k}\bigr).
\end{align}
Each projected polynomial represents the corresponding DOF moments
through the $L^2$ inner product.
Local unisolvence (Theorems~\ref{thm:trimmed-local-unisolvence}
and~\ref{thm:complete-local-unisolvence}), together with global
assembly, then makes these maps isomorphisms.
In our degree convention, these data spaces and moment maps are
precisely those of the VEM-inspired and DDR-inspired constructions
in \cite[Sections~4.1 and~3.1]{Bonaldi2025}, respectively.
Theorem~\ref{thm:fully-discrete-cochain-isomorphisms}
in Appendix~\ref{subsec:appendix-fully-discrete-realizations}
establishes that these DOF maps are cochain isomorphisms.
Equivalently, for $0\le k<n$, the following diagrams commute,
with isomorphisms as vertical maps:
\begin{equation}
\label{eq:fully-discrete-commuting-diagrams}
\begin{tikzcd}[column sep=4em,row sep=2.6em]
\Vcal_{r+1,h}^{-,k}
  \arrow[r,"\rd"]
  \arrow[d,"\mathcal J_{r+1,h}^{-,k}"']
& \Vcal_{r+1,h}^{-,k+1}
  \arrow[d,"\mathcal J_{r+1,h}^{-,k+1}"] \\
\underline{\Vcal}_{r+1,h}^{-,k}
  \arrow[r,"\underline{\rd}_{r+1,h}^{-,k}"']
& \underline{\Vcal}_{r+1,h}^{-,k+1}
\end{tikzcd}
\qquad
\begin{tikzcd}[column sep=4em,row sep=2.6em]
\Vcal_{r,h}^{k}
  \arrow[r,"\rd"]
  \arrow[d,"\mathcal J_{r,h}^{k}"']
& \Vcal_{r,h}^{k+1}
  \arrow[d,"\mathcal J_{r,h}^{k+1}"] \\
\underline{\Vcal}_{r,h}^{k}
  \arrow[r,"\underline{\rd}_{r,h}^{k}"']
& \underline{\Vcal}_{r,h}^{k+1}
\end{tikzcd}
\end{equation}
where $\underline{\rd}_{r+1,h}^{-,k}$ and
$\underline{\rd}_{r,h}^{k}$ are the discrete exterior derivatives
defined in \cite{Bonaldi2025} for the VEM-inspired and
DDR-inspired complexes, respectively.
\end{remark}

\begin{remark}[Serendipity reduction of DOFs]
\label{rem:comparison-serendipity}
Unlike higher-order reconstruction, which is inherent to DDR
and analogous to enhancement in VEM, serendipity reduces DOFs
according to the geometry of mesh cells
\cite{BeiraoBrezziMariniRusso2016,
BeiraoBrezziMariniRusso2017,BeiraoEtAl2018}.
Recently, these two strategies were combined in a framework for
polytopal complexes \cite{DiPietroDroniou2023Serendipity}.
Its application to the three-dimensional DDR complex yields the
serendipity DDR (SDDR) complex with more competitive DOF counts;
see \cite[Table~3]{DiPietroDroniou2023Serendipity}.
Extending this construction to the present complete complex
of differential forms is beyond the scope of this work.
\end{remark}

\subsubsection{Commuting canonical interpolators}
\label{subsubsec:commuting-canonical-interpolators}

Fix \(\epsi>0\) and define
\begin{equation}
\label{eq:commuting-regularity-domain}
  \mathscr S_\epsi^k
  :=
  H^{\frac{n-k}{2}+\epsi}\Lambda^k(\Omega),
  \qquad
  \mathscr D_\epsi^k
  :=
  \left\{
    \omega\in\mathscr S_\epsi^k:
    \rd\omega\in\mathscr S_\epsi^{k+1}
  \right\},
  \quad 0\le k\le n-1,
\end{equation}
with \(\mathscr D_\epsi^n:=\mathscr S_\epsi^n\).
Under these regularity assumptions, matching the global versions of
\eqref{eq:trimmed-dof-constant}--\eqref{eq:trimmed-dof-derivative}
and the moments in \eqref{eq:global-complete-dofs}, respectively,
defines the canonical interpolators
\[
  \mathcal I_{r+1,h}^{-,k}:
  \mathscr D_\epsi^k\longrightarrow\Vcal_{r+1,h}^{-,k},
  \qquad
  \mathcal I_{r,h}^{k}:
  \mathscr S_\epsi^k\longrightarrow\Vcal_{r,h}^{k}.
\]

\begin{theorem}[Commuting canonical interpolators]
\label{thm:commutativity}
For \(0\le k\le n-1\) and \(\omega\in\mathscr D_\epsi^k\),
\begin{equation}
\label{eq:virtual-interpolator-commutativity}
  \rd\mathcal I_{r+1,h}^{-,k}\omega
  =
  \mathcal I_{r+1,h}^{-,k+1}\rd\omega,
  \qquad
  \rd\mathcal I_{r,h}^{k}\omega
  =
  \mathcal I_{r,h}^{k+1}\rd\omega.
\end{equation}
\end{theorem}

\begin{proof}
By \eqref{eq:trimmed-virtual-derham-complex} and the definition of the trimmed interpolator,
both sides of the first identity in
\eqref{eq:virtual-interpolator-commutativity} belong to
\(\Vcal_{r+1,h}^{-,k+1}\).
By global unisolvence, it suffices to compare the three classes of moments
\eqref{eq:trimmed-dof-constant}--\eqref{eq:trimmed-dof-derivative},
with \(k\) replaced by \(k+1\), on every
\(f\in\Delta_\ell(\Mh)\), \(k+1\le\ell\le n\).
The derivative moments vanish by the interpolation conditions
and \(\rd^2=0\):
\[
  \mathsf D_{\ell,h}^{k+1,-}
  \bigl(\rd\mathcal I_{r+1,h}^{-,k}\omega\bigr)
  =0=
  \mathsf D_{\ell,h}^{k+1,-}
  \bigl(\mathcal I_{r+1,h}^{-,k+1}\rd\omega\bigr),
  \qquad k+2\le \ell\le n.
\]
As noted in Remark~\ref{rem:equivalent-dofs}, the
\(\mathsf K_{\ell,h}^{k+1,-}\) moments of \(\rd\omega\) are precisely
the \(\mathsf D_{\ell,h}^{k,-}\) moments of \(\omega\).
Thus, for each test form \(\nu\) in \eqref{eq:trimmed-dof-koszul},
with \(k\) replaced by \(k+1\), the interpolation conditions give
\[
\begin{aligned}
  \int_f\trc_f\bigl(\rd\mathcal I_{r+1,h}^{-,k}\omega\bigr)\wedge\nu
  =\int_f\trc_f\rd\omega\wedge\nu=\int_f\trc_f\bigl(\mathcal I_{r+1,h}^{-,k+1}\rd\omega\bigr)\wedge\nu.
\end{aligned}
\]
Finally, for \(\ell=k+1\), Stokes' theorem
and the interpolation conditions for the constant moments
\eqref{eq:trimmed-dof-constant} give
\[
  \int_f\trc_f\bigl(\rd\mathcal I_{r+1,h}^{-,k}\omega
    -\mathcal I_{r+1,h}^{-,k+1}\rd\omega\bigr)
  =\int_{\partial f}\trc_{\partial f}
    \bigl(\mathcal I_{r+1,h}^{-,k}\omega-\omega\bigr)=0.
\]
Thus all three classes of moments agree, proving the first identity.

For the complete family, let
\(\mu\in\Pcal_r^-\Lambda^{\ell-k-1}(f)\).
By \eqref{eq:polynomial-differential-images} and
Lemma~\ref{lem:traces-trimmed}, \(\rd\mu\) and the facet traces
of \(\mu\) belong to the corresponding test spaces in
\eqref{eq:global-complete-dofs}.
The integration-by-parts formula \eqref{eq:stokes-integration-by-parts} and the interpolation
conditions therefore give
\[
\begin{aligned}
  \int_f\trc_f\bigl(\rd\mathcal I_{r,h}^{k}\omega
    -\mathcal I_{r,h}^{k+1}\rd\omega\bigr)\wedge\mu
  &=\int_{\partial f}\trc_{\partial f}
    \bigl(\mathcal I_{r,h}^{k}\omega-\omega\bigr)
    \wedge\trc_{\partial f}\mu\\
  &\quad+(-1)^{k+1}\int_f\trc_f
    \bigl(\mathcal I_{r,h}^{k}\omega-\omega\bigr)\wedge\rd\mu=0.
\end{aligned}
\]
Property \textup{(C3)} and global unisolvence prove the second identity.
\end{proof}

The same moment-matching argument defines cochain reductions directly
on the virtual spaces.
Write $\mathcal X_{r,h}^{\bullet}$ for either
$\Vcal_{r+1,h}^{-,\bullet}$ or $\Vcal_{r,h}^{\bullet}$.

\begin{corollary}[Cochain reductions induced by interpolation]
\label{cor:cochain-reductions}
For either choice of $\mathcal X_{r,h}^{\bullet}$ and $0\le k\le n$,
define $Q_{r,h}^k:\mathcal X_{r,h}^k\to\Vcal_{0,h}^k$ by
\begin{equation*}
  \int_f\trc_f Q_{r,h}^k\omega_h
  =\int_f\trc_f\omega_h
  \qquad
  \forall f\in\Delta_k(\Mh),
\end{equation*}
and define
$P_{r,h}^k:\Vcal_{r+1,h}^{-,k}\to\Vcal_{r,h}^k$ by
\begin{equation}
\label{eq:complete-moment-reduction}
  \int_e\trc_e(P_{r,h}^k \omega_h)\wedge\mu
  =\int_e\trc_e \omega_h\wedge\mu,
  \quad
  \mu\in\Pcal_r^-\Lambda^{d-k}(e),
  \quad
  e\in\Delta_d(\Mh),
  \quad k\le d\le n.
\end{equation}
Then $Q_{r,h}^{\bullet}$ is a surjective cochain map and
$P_{r,h}^{\bullet}$ is a cochain retraction:
\begin{align}
\label{eq:Q-cochain-reduction}
  Q_{r,h}^{k+1}\rd
  &=\rd Q_{r,h}^k,
  &
  Q_{r,h}^k(\mathcal X_{r,h}^k)
  &=\Vcal_{0,h}^k,
  \\
\label{eq:P-cochain-reduction}
  P_{r,h}^{k+1}\rd
  &=\rd P_{r,h}^k,
  &
  P_{r,h}^k\big|_{\Vcal_{r,h}^k}
  &=\operatorname{id}.
\end{align}
The cochain identities hold for $0\le k<n$.
\end{corollary}

\begin{proof}
The cochain identities follow from Stokes' theorem and global
unisolvence, as in the proof of Theorem~\ref{thm:commutativity}.
The $k$-cell integrals (vertex values when $k=0$) can be prescribed arbitrarily in either
family, so global unisolvence makes $Q_{r,h}^k$ surjective.
Finally, for $\omega_h\in\Vcal_{r,h}^k\subset\Vcal_{r+1,h}^{-,k}$,
$P_{r,h}^k\omega_h$ and $\omega_h$ have the same complete DOFs, so
unisolvence gives $P_{r,h}^k\omega_h=\omega_h$.
\end{proof}

\subsubsection{De Rham cohomology}
\label{subsubsec:de-rham-cohomology}
In \cite[Sections~3.6 and~4.6]{Bonaldi2025}, cohomology isomorphisms
with the continuous de Rham complex are established separately for
the two fully discrete families by reduction to lowest order and
exactness of their zero-average subcomplexes.
These results can be transferred directly to our virtual complexes
through the cochain isomorphisms of
Theorem~\ref{thm:fully-discrete-cochain-isomorphisms}.
Here we give a unified proof within the conforming framework using
the cochain reductions of Corollary~\ref{cor:cochain-reductions},
with $P_{r,h}^{\bullet}$ connecting the trimmed and complete complexes.

For a cochain complex $(A^{\bullet},\mathrm d_A)$, write
\begin{equation*}
  H^k(A^{\bullet})
  :=\frac{\ker(\mathrm d_A:A^k\to A^{k+1})}
  {\operatorname{im}(\mathrm d_A:A^{k-1}\to A^k)}.
\end{equation*}
All complexes are extended by zero outside degrees $0,\ldots,n$.
A cochain map $F^{\bullet}:A^{\bullet}\to B^{\bullet}$ induces
$H^k(F^{\bullet}):H^k(A^{\bullet})\to H^k(B^{\bullet})$ by
$[a]\mapsto[F^ka]$.

Following the argument in~\cite[Sections~2.2 and~5.5]{Arnold2006},
the de Rham map, defined by integration over oriented cells,
is a cochain isomorphism from $\Vcal_{0,h}^{\bullet}$ onto the
real cellular cochain complex of $\Mh$: it is bijective by
unisolvence and commutes with the differentials by Stokes' theorem.
The cellular cohomology spaces are, in turn, isomorphic to those
of~\eqref{eq:domain-derham-complex}, yielding
\begin{equation}
\label{eq:lowest-order-cohomology}
  H^k(\Vcal_{0,h}^{\bullet})
  \cong H^k\bigl(H\Lambda^{\bullet}(\rd;\Omega)\bigr),
  \qquad 0\le k\le n.
\end{equation}

For either virtual family, define the spaces with vanishing averages
on $k$-cells by
\begin{equation}
\label{eq:zero-moment-space}
  \mathcal X_{r,h,\#}^k
  :=\ker Q_{r,h}^k
  =\left\{
    \omega_h\in\mathcal X_{r,h}^k:
    \int_f\trc_f\omega_h=0
    \quad\forall f\in\Delta_k(\Mh)
  \right\}.
\end{equation}

\begin{lemma}[Exactness of the zero-average subcomplex]
\label{lem:zero-moment-exactness}
For either choice of $\mathcal X_{r,h}^{\bullet}$, the zero-average
spaces form the complex
\begin{equation}
\label{eq:zero-moment-exactness}
  0\longrightarrow\mathcal X_{r,h,\#}^{0}
  \xrightarrow{\rd}\mathcal X_{r,h,\#}^{1}
  \xrightarrow{\rd}\cdots
  \xrightarrow{\rd}\mathcal X_{r,h,\#}^{n}
  \longrightarrow0.
\end{equation}
with $H^k(\mathcal X_{r,h,\#}^{\bullet})=0$ for $0\le k\le n$.
\end{lemma}

\begin{proof}
By \eqref{eq:Q-cochain-reduction}, the zero-average spaces form a
subcomplex. In degree zero, $\rd\omega_h=0$ implies that $\omega_h$
is constant on each connected component. Its vanishing vertex
values therefore give $\omega_h=0$.

For the trimmed family, let $1\le k\le n$ and let
$\omega_h\in\Vcal_{r+1,h,\#}^{-,k}$ be closed.
Using the moment correspondence in Remark~\ref{rem:trimmed-moment-shift},
prescribe $\phi_h\in\Vcal_{r+1,h}^{-,k-1}$ by
\[
\begin{aligned}
  &\mathsf C_h^{k-1,-}(\phi_h)=0,\qquad
  \mathsf K_{\ell,h}^{k-1,-}(\phi_h)=0,
  &&k-1\le\ell\le n,\\
  &\mathsf D_{\ell,h}^{k-1,-}(\phi_h)
  =\mathsf K_{\ell,h}^{k,-}(\omega_h),
  &&k\le\ell\le n.
\end{aligned}
\]
Global unisolvence determines a unique $\phi_h$, and the first
condition gives $\phi_h\in\Vcal_{r+1,h,\#}^{-,k-1}$.

We compare $\rd\phi_h$ with $\omega_h$. Their derivative classes
vanish because $\rd^2\phi_h=0$ and $\rd\omega_h=0$.
The prescribed derivative moments give
\[
  \mathsf K_{\ell,h}^{k,-}(\rd\phi_h)
  =\mathsf D_{\ell,h}^{k-1,-}(\phi_h)
  =\mathsf K_{\ell,h}^{k,-}(\omega_h),
  \qquad k\le\ell\le n.
\]
For $f\in\Delta_k(\Mh)$, Stokes' theorem
and the vanishing $\mathsf C_h^{k-1,-}$ moments give
\[
  \int_f\trc_f\rd\phi_h
  =\int_{\partial f}\trc_{\partial f}\phi_h
  =0
  =\int_f\trc_f\omega_h.
\]
By \textup{(T3)}, $\rd\phi_h\in\Vcal_{r+1,h}^{-,k}$.
Since all three classes of moments agree, global unisolvence
gives $\rd\phi_h=\omega_h$.

For the complete family, let $\omega_h\in\Vcal_{r,h,\#}^{k}$ be
closed, with $1\le k\le n$. By \eqref{eq:complete-space},
$\omega_h\in\Vcal_{r+1,h,\#}^{-,k}$, so trimmed exactness gives
$\phi_h^-\in\Vcal_{r+1,h,\#}^{-,k-1}$ with
$\rd\phi_h^-=\omega_h$. Set
$\phi_h:=P_{r,h}^{k-1}\phi_h^-$. By
\eqref{eq:complete-moment-reduction} and \eqref{eq:P-cochain-reduction},
\[
  \phi_h\in\Vcal_{r,h,\#}^{k-1},
  \qquad
  \rd\phi_h
  =P_{r,h}^{k}\rd\phi_h^-
  =P_{r,h}^{k}\omega_h
  =\omega_h.
\]
\end{proof}

\begin{theorem}[De Rham cohomology of the virtual complexes]
\label{thm:veec-cohomology}
For every $r\ge0$, the trimmed and complete virtual de Rham
complexes \eqref{eq:trimmed-virtual-derham-complex} and
\eqref{eq:complete-virtual-derham-complex} and the $L^2$ de Rham
complex \eqref{eq:domain-derham-complex} satisfy, for $0\le k\le n$,
\[
\begin{aligned}
  H^k(\Vcal_{r+1,h}^{-,\bullet})
  \cong H^k(\Vcal_{r,h}^{\bullet})
  \cong H^k\bigl(H\Lambda^{\bullet}(\rd;\Omega)\bigr).
\end{aligned}
\]
\end{theorem}

\begin{proof}
For either virtual family, \eqref{eq:Q-cochain-reduction} and
\eqref{eq:zero-moment-space} give the short exact sequence of
cochain complexes
\begin{equation}
\label{eq:cohomology-short-exact-sequence}
  0\longrightarrow
  \mathcal X_{r,h,\#}^{\bullet}
  \longrightarrow
  \mathcal X_{r,h}^{\bullet}
  \xrightarrow{Q_{r,h}^{\bullet}}
  \Vcal_{0,h}^{\bullet}
  \longrightarrow0.
\end{equation}
By Lemma~\ref{lem:zero-moment-exactness}, the kernel complex
\eqref{eq:zero-moment-exactness} satisfies
$H^j(\mathcal X_{r,h,\#}^{\bullet})=0$ for every $j$.
The long exact cohomology sequence associated with
\eqref{eq:cohomology-short-exact-sequence}, given by
\cite[Theorem~1.3.1]{Weibel1994}, contains the exact segment
\begin{equation*}
  H^k(\mathcal X_{r,h,\#}^{\bullet})
  \longrightarrow H^k(\mathcal X_{r,h}^{\bullet})
  \xrightarrow{H^k(Q_{r,h}^{\bullet})}
  H^k(\Vcal_{0,h}^{\bullet})
  \xrightarrow{\partial^k}
  H^{k+1}(\mathcal X_{r,h,\#}^{\bullet}).
\end{equation*}
Here $\partial^k$ is the connecting homomorphism, and the two outer
groups vanish.
Exactness therefore makes $H^k(Q_{r,h}^{\bullet})$ both injective
and surjective, hence an isomorphism for every $0\le k\le n$.
Equation~\eqref{eq:lowest-order-cohomology} now proves the assertion
for both virtual families.
\end{proof}

\par
\section{Mesh-uniform estimates for virtual differential forms}
\label{sec:local-estimates}
We establish inverse estimates and norm equivalences, which are central
to the stability and convergence analysis of virtual element methods,
for both virtual families in arbitrary dimension.
Using the upper bounds in these equivalences, we further derive
canonical interpolation estimates.
For a family of polytopal meshes $\{\Mh\}_{h>0}$, let
$f\in\Delta_d(\Mh)$, $1\le d\le n$, denote a generic
polytopal $d$-cell throughout this section.
After normalizing the length scale so that $\diam(\Omega)\le1$,
we set $h_f:=\diam(f)$ and $h:=\max_{T\in\Th}h_T$.
We impose the following mesh regularity assumptions
\cite{BeiraoEtAl2013,BeiraoEtAl2018,BeiraoMascottoMeng2022}.
\begin{assumption}[Mesh regularity]
\label{ass:mesh-regularity}
\leavevmode\par
\begin{itemize}
\item[(M1)]
Every $f$ is star-shaped with respect to a ball
$B_f\subset f$ in $\operatorname{aff}(f)$ of radius at least
$\rho h_f$, where $\rho>0$ is independent of the mesh.
\item[(M2)]
There exists a mesh-independent constant $\gamma>0$ such that
$h_f\ge\gamma h_T$ for every $T\in\Th$ and
$f\in\Delta_d(T)$, $1\le d\le n-1$.
\end{itemize}
\end{assumption}
Condition {\rm (M1)} prevents cells from becoming arbitrarily thin
relative to their diameters, while {\rm (M2)} excludes
proper subcells of positive dimension that are arbitrarily small relative to
their containing elements. Together, these conditions imply a
uniform bound on the number of subcells per cell.

For nonnegative quantities $a$ and $b$, we write $a\lesssim b$,
or equivalently $b\gtrsim a$, if $a\le Cb$, and $a\eqsim b$
if $a\lesssim b$ and $b\lesssim a$.
Unless otherwise stated, hidden constants are independent of $h$
but may depend on $n$, $r$, $\rho$, $\gamma$, and the Sobolev orders.
Constants in global estimates may also depend on $\Omega$.

\subsection{Scaled estimates on star-shaped domains}
\label{subsec:mesh-auxiliary-estimates}
The estimates below require only the star-shapedness condition
{\rm (M1)}, with constants depending only on $\rho$ for fixed
dimension and Sobolev orders, and additionally on $r$ for polynomial
estimates and approximation bounds.
Proofs of the trace and Poincar\'e estimates are given in
Appendix~\ref{sec:appendix-local}.

\begin{lemma}[Scaled trace estimates]
\label{lem:scaled-traces}
Let $0\le k\le d$.
The following estimates hold.
\begin{enumerate}
\renewcommand{\theenumi}{\roman{enumi}}
\renewcommand{\labelenumi}{\textup{(\theenumi)}}
\item\label{item:trace-rellich}
\textup{(Rellich-type trace estimates)} Let $\omega\in H\Lambda^k(\rd;f)\cap H\Lambda^k(\delta;f)$.
If $\trc\omega$ is square-integrable on $\partial f$,
then so is $\ntr\omega$, and the following multiplicative estimate holds:
\begin{equation}
\label{eq:cell-div-curl-trace}
\begin{aligned}
  \norm{\ntr\omega}_{\partial f}^2
  \lesssim
  h_f^{-1}\norm{\omega}_f^2
  +\norm{\omega}_f
  \bigl(\norm{\rd\omega}_f
  +\norm{\delta\omega}_f\bigr)
  +\norm{\trc\omega}_{\partial f}^2.
\end{aligned}
\end{equation}
Consequently, the following additive estimate holds:
\begin{equation}
\label{eq:cell-div-curl-trace-additive}
\begin{aligned}
  \norm{\ntr\omega}_{\partial f}
  \lesssim
  h_f^{-1/2}\norm{\omega}_f
  +h_f^{1/2}\bigl(\norm{\rd\omega}_f
  +\norm{\delta\omega}_f\bigr)
  +\norm{\trc\omega}_{\partial f}.
\end{aligned}
\end{equation}
The same estimates hold with $\trc$ and $\ntr$ interchanged.

\item\label{item:trace-weak}
\textup{(Negative-norm trace estimates)} Let $d\ge2$ and $k<d$.
For every $\omega\in H\Lambda^k(\rd;f)$ and every facet
$e\in\Delta_{d-1}(f)$,
\begin{equation}
\label{eq:weak-trace}
  h_f^{-1/2}\norm{\trc_e\omega}_{-1,e}
  \lesssim\norm{\trc\omega}_{-1/2,\partial f}
  \lesssim\norm{\omega}_f+h_f\norm{\rd\omega}_f.
\end{equation}
Here $\norm{\cdot}_{-1/2,\partial f}$ denotes the scaled dual
trace norm defined in Appendix~\ref{app:facet-negative-trace}.
\item\label{item:trace-facets}
\textup{(Sobolev trace estimates on facets)} Let $k<d$,
$s>1/2$, and $e\in\Delta_{d-1}(f)$.
For every $\omega\in H^s\Lambda^k(f)$,
\begin{equation}
\label{eq:facet-sobolev-trace}
\begin{aligned}
  \norm{\trc_e\omega}_e
  +h_f^{s-1/2}\seminorm{\trc_e\omega}_{s-1/2,e}
  \lesssim
  h_f^{-1/2}\norm{\omega}_f
  +h_f^{s-1/2}\seminorm{\omega}_{s,f}.
\end{aligned}
\end{equation}
For $d=1$, the left-hand side denotes the magnitude of the
vertex value, $|\omega(e)|$.
\end{enumerate}
\end{lemma}
\begin{lemma}[Scaled Poincar\'e inequality]
\label{lem:cell-poincare}
Let $0\le k\le d$.
Every $\omega$ in
$\overset{\circ}{H}\Lambda^k(\rd;f)
\cap H\Lambda^k(\delta;f)$ or
$H\Lambda^k(\rd;f)
\cap\overset{\circ}{H}\Lambda^k(\delta;f)$
satisfies
\begin{equation}
\label{eq:cell-poincare}
  \norm{\omega}_f
  \lesssim
  h_f\bigl(\norm{\rd\omega}_f
  +\norm{\delta\omega}_f\bigr),
\end{equation}
provided that $\int_f\omega=0$ in the first space when $k=d$,
and $\ip{\omega}{1}_f=0$ in the second space when $k=0$.
\end{lemma}
\par

We also collect the polynomial estimates used below.

\begin{lemma}[Polynomial estimates]
\label{lem:polynomial-estimates}
Let $r\ge0$ and $0\le k\le d$.
Every
$\mu\in\Pcal_r\Lambda^k(f)$ satisfies
\begin{equation}
\label{eq:polynomial-estimates}
  \norm{\trc\mu}_{\partial f}
  +\norm{\ntr\mu}_{\partial f}
  \lesssim h_f^{-1/2}\norm{\mu}_f,
  \qquad
  \seminorm{\mu}_{1,f}
  \lesssim h_f^{-1}\norm{\mu}_f,
\end{equation}
and
\begin{equation}
\label{eq:polynomial-negative-inverse}
  \norm{\mu}_f\lesssim h_f^{-j}\norm{\mu}_{-j,f},
  \qquad
  j=1,2.
\end{equation}
\end{lemma}
\begin{proof}
The star-shapedness condition {\rm (M1)} yields the tangential and
normal trace bounds and the $H^1$ inverse estimate as in
\cite[Lemma~2.3]{BrennerSung2018}, and the negative-norm bounds as in
\cite[Lemma~10]{Huang2020}, applied componentwise.
\end{proof}
\par

For the interpolation analysis, we also require a polynomial
approximation that commutes with the exterior derivative.

\begin{lemma}[Commuting polynomial approximation]
\label{lem:polynomial-approximation}
Let $r\ge0$ and $0\le k\le d$.
There exists a polynomial projector
$Q_{r,f}^k:L^2\Lambda^k(f)\to\Pcal_r\Lambda^k(f)$.
For real $s,\ell$ with $0<\ell\le s\le r+1$, the projection
$\omega_r^\pi:=Q_{r,f}^k\omega$ of every
$\omega\in H^s\Lambda^k(f)$ satisfies
\begin{equation}
\label{eq:section4-poly-approx}
  \norm{\omega-\omega_r^\pi}_f
  +
  h_f^\ell
  \seminorm{\omega-\omega_r^\pi}_{\ell,f}
  \lesssim
  h_f^s\seminorm{\omega}_{s,f}.
\end{equation}
For $r\ge1$ and $k<d$, it satisfies
$\rd Q_{r,f}^k\omega=Q_{r-1,f}^{k+1}\rd\omega$
for every $\omega\in H\Lambda^k(\rd;f)$.
Consequently, for real $t,\ell$ with $0<\ell\le t\le r$ and
$\rd\omega\in H^t\Lambda^{k+1}(f)$, the same polynomial
$\omega_r^\pi$ satisfies
\begin{equation}
\label{eq:section4-poly-approx-derivative}
  \norm{
    \rd\omega-\rd\omega_r^\pi
  }_f
  +
  h_f^\ell
  \seminorm{
    \rd\omega-\rd\omega_r^\pi
  }_{\ell,f}
  \lesssim
  h_f^t
  \seminorm{\rd\omega}_{t,f}.
\end{equation}
\end{lemma}
\begin{proof}
Take $Q_{r,f}^k$ to be the componentwise averaged Taylor projector
on $L^2$ of \cite[Section~4.1]{BrennerScott2008}, using the same
smooth averaging function scaled to $B_f$ for all polynomial and form degrees.
Under the star-shapedness assumption {\rm (M1)}, the
Bramble--Hilbert estimate
\cite[Lemma~4.3.8]{BrennerScott2008} and Sobolev interpolation
\cite[Chapter~7]{AdamsFournier2003}, applied modulo polynomials,
yield \eqref{eq:section4-poly-approx}.
The scalar differentiation identity
\cite[Proposition~4.1.17]{BrennerScott2008}
extends to distributions via the smooth-kernel representation
\cite[Proposition~4.1.12]{BrennerScott2008}, yielding
$\rd Q_{r,f}^k=Q_{r-1,f}^{k+1}\rd$ componentwise.
Applying \eqref{eq:section4-poly-approx} to $\rd\omega$ with
degree $r-1$ gives \eqref{eq:section4-poly-approx-derivative}.
\end{proof}
\par

\subsection{Inverse estimates}
\label{subsec:appendix-virtual-gaffney}
Using the recursive structure of the local spaces, we establish
negative-norm inverse estimates for all form degrees.
In particular, for $k=d-1$, the corresponding estimate for
three-dimensional face virtual elements is proved in
\cite[Lemma~5.1, Eq.~(5.7)]{BeiraoMascottoMeng2022}.
We begin with an orthogonal Hodge decomposition into exact and coexact
components.
\par

\begin{lemma}[Orthogonal Hodge decomposition]
\label{lem:relative-hodge-decomposition}
Let $0\le k<d$ and
$\omega\in H\Lambda^k(\rd;f)\cap H\Lambda^k(\delta;f)$.
Then
\begin{equation}
\label{eq:relative-hodge-decomposition}
  \omega=\rd\phi+\delta\xi,
  \qquad \ip{\rd\phi}{\delta\xi}_f=0,
\end{equation}
with the exact term omitted when $k=0$.
For $k\ge1$, $\phi$ is the unique element of
$\overset{\circ}{H}\Lambda^{k-1}(\rd;f)
\cap H\Lambda^{k-1}(\delta;f)$ satisfying
\begin{equation}
\label{eq:relative-exact-potential}
  \delta\rd\phi=\delta\omega,
  \qquad \delta\phi=0.
\end{equation}
The potential $\xi$ is the unique element of
$\overset{\circ}{H}\Lambda^{k+1}(\rd;f)
\cap H\Lambda^{k+1}(\delta;f)$ satisfying
\begin{equation}
\label{eq:relative-coexact-potential}
  \rd\delta\xi=\rd\omega,
  \qquad \rd\xi=0,
  \qquad \trc\delta\xi=\trc\omega,
\end{equation}
and $\int_f\xi=0$ when $k=d-1$.
If $\delta\omega=0$, then $\rd\phi=0$ and $\omega=\delta\xi$.
\end{lemma}

\begin{proof}
Since $f$ is contractible and $k<d$, $\mathring{\Harm}^k(f)=\{0\}$.
Thus the Hodge--Laplace problem \eqref{eq:strong-hodge-laplace}
for the relative complex
$V^\bullet(f)=\overset{\circ}{H}\Lambda^\bullet(\rd;f)$
with right-hand side $\omega$ has a unique solution $u$.
Set
\[
  \phi:=\delta u,
  \qquad
  \xi:=\rd u,
\]
with the exact potential omitted when $k=0$.
The homogeneous boundary conditions encoded in
\eqref{eq:hodge-laplacian-domain}, together with the complex identities,
give $\trc\phi=\trc\xi=0$, $\delta\phi=0$, and $\rd\xi=0$.
The strong equation \eqref{eq:strong-hodge-laplace} and the orthogonal
Hodge decomposition \eqref{eq:hodge-decomposition-homogeneous} yield
\eqref{eq:relative-hodge-decomposition}.

To verify the equations for the potentials, apply $\delta$ and $\rd$
to \eqref{eq:relative-hodge-decomposition} in the distributional sense.
Since $\delta^2=\rd^2=0$ and the assumed regularity gives
$\delta\omega,\rd\omega\in L^2$, we obtain
\[
  \delta\rd\phi=\delta\omega,
  \qquad
  \rd\delta\xi=\rd\omega.
\]
Consequently, $\rd\phi\in H\Lambda^k(\delta;f)$ and
$\delta\xi\in H\Lambda^k(\rd;f)$.
The homogeneous condition $\trc\phi=0$ and trace commutativity give
$\trc\rd\phi=\rd\trc\phi=0$.
Taking the tangential trace of \eqref{eq:relative-hodge-decomposition}
therefore yields $\trc\delta\xi=\trc\omega$.
When $k=d-1$, Stokes' theorem gives
$\int_f\xi=\int_f\rd u=0$.
\par
For uniqueness, let $\zeta$ be the difference of two solutions for either
potential. The integration-by-parts formula
\eqref{eq:green-formula} gives $\rd\zeta=\delta\zeta=0$;
\eqref{eq:cell-poincare} then implies $\zeta=0$.
Finally, if $\delta\omega=0$ and $k\ge1$, integration by parts gives
\[
  \norm{\rd\phi}_f^2
  =\ip{\delta\rd\phi}{\phi}_f
  =\ip{\delta\omega}{\phi}_f
  =0.
\]
Thus $\omega=\delta\xi$.
For $k=0$, this conclusion follows directly from the decomposition.
\end{proof}
\par

We also use the standard negative-norm bound
\begin{equation}
\label{eq:negative-differentiation}
  \norm{\rd\omega}_{-j,f}+\norm{\delta\omega}_{-j,f}
  \lesssim\norm{\omega}_{1-j,f},
  \qquad j=1,2,
\end{equation}
for $\omega\in H^{1-j}\Lambda^k(f)$, $0\le k\le d$.
In fact, for $\zeta\in H_0^j\Lambda^{k+1}(f)$, we have
$\delta\zeta\in H_0^{j-1}\Lambda^k(f)$, where $H_0^0=L^2$, and
\[
\begin{aligned}
  |\ip{\rd\omega}{\zeta}_f|
  &=|\ip{\omega}{\delta\zeta}_f|
  \le\norm{\omega}_{1-j,f}\norm{\delta\zeta}_{j-1,f}\\
  &\lesssim\norm{\omega}_{1-j,f}\norm{\zeta}_{j,f}.
\end{aligned}
\]
Taking the supremum over $\zeta$ gives the $\rd$-bound.
The $\delta$-bound follows analogously by testing against
$H_0^j\Lambda^{k-1}(f)$ and interchanging $\rd$ and $\delta$.

\begin{lemma}[Negative-norm inverse estimate]
\label{lem:negative-inverse}
Under conditions {\rm (M1)} and {\rm (M2)}, let $0\le k\le d$.
Every
$\omega_h\in\Vcal_{r+1}^{-,k}(f)$ satisfies
\begin{equation}
\label{eq:negative-inverse}
  \norm{\omega_h}_f\lesssim h_f^{-1}\norm{\omega_h}_{-1,f}.
\end{equation}
\end{lemma}

\begin{proof}
\emph{Step 1: induction and decomposition.}
We argue by induction on $\lambda:=d-k$.
For $\lambda=0$, \eqref{eq:trimmed-top} gives
$\omega_h\in\Pcal_r\Lambda^d(f)$, so
\eqref{eq:polynomial-negative-inverse} with $j=1$ proves
\eqref{eq:negative-inverse}.

Assume that the result holds whenever $d-k\le\lambda$,
and consider $d-k=\lambda+1$.
For $(d,k)=(1,0)$, \eqref{eq:trimmed-one-dimensional-scalar}
and the polynomial negative-norm inverse estimate
\eqref{eq:polynomial-negative-inverse}, applied with
polynomial degree $r+1$ and $j=1$, prove the claim.
We may therefore assume $d\ge2$.
Applying Lemma~\ref{lem:relative-hodge-decomposition}, with the
exact term omitted when $k=0$, we write
\begin{equation}
\label{eq:negative-inverse-decomposition}
  \omega_h=\rd\phi+\delta\xi.
\end{equation}
We estimate the two components separately as follows.

\emph{Step 2: the exact component.}
For $k\ge1$, the defining conditions
\eqref{eq:trimmed-codim-one} and \eqref{eq:trimmed-general}
give
\[
  \mu_h:=\delta\omega_h=\delta\rd\phi
  \in\star^{-1}\rd\Kcal_{r+1}^{d-k}(f)
  \subset\Pcal_r\Lambda^{k-1}(f).
\]
Since $\trc\phi=0$ and $\delta\phi=0$,
\eqref{eq:cell-poincare} gives
$\norm{\phi}_f\lesssim h_f\norm{\rd\phi}_f$.
The integration-by-parts formula \eqref{eq:green-formula}
therefore yields
\begin{equation}
\label{eq:exact-by-polynomial-data}
  \norm{\rd\phi}_f^2
  =\ip{\mu_h}{\phi}_f
  \lesssim h_f\norm{\mu_h}_f\norm{\rd\phi}_f.
\end{equation}
Since $\mu_h=\delta\omega_h$, combining
\eqref{eq:exact-by-polynomial-data} with  \eqref{eq:polynomial-negative-inverse} and the
$\delta$-estimate in \eqref{eq:negative-differentiation}
, both with $j=2$,
gives
\begin{equation}
\label{eq:exact-negative-bound}
  \norm{\rd\phi}_f
  \lesssim h_f\norm{\mu_h}_f
  \lesssim h_f^{-1}\norm{\mu_h}_{-2,f}
  \lesssim h_f^{-1}\norm{\omega_h}_{-1,f}.
\end{equation}
For $k=0$, this bound is understood with $\rd\phi=0$.

\emph{Step 3: the coexact component.}
The identities in \eqref{eq:relative-coexact-potential} give
\[
  \rd\delta\xi=\rd\omega_h,
  \qquad
  \trc\delta\xi=\trc\omega_h.
\]
Properties {\rm (T3)} and {\rm (T2)} therefore place
$\rd\delta\xi$ in $\Vcal_{r+1}^{-,k+1}(f)$ and
$\trc_e\delta\xi$ in $\Vcal_{r+1}^{-,k}(e)$ for each
$e\in\Delta_{d-1}(f)$.
Both have induction parameter $\lambda$.
The induction hypothesis and the $\rd$-estimate in
\eqref{eq:negative-differentiation} with $j=1$ yield
\begin{equation}
\label{eq:induction-delta-xi}
  \norm{\rd\delta\xi}_f
  \lesssim h_f^{-1}\norm{\rd\delta\xi}_{-1,f}
  \lesssim h_f^{-1}\norm{\delta\xi}_f.
\end{equation}
On each facet $e$, the induction hypothesis and the left-hand
inequality of \eqref{eq:weak-trace}, applied to $\delta\xi$, give
\[
  \norm{\trc_e\delta\xi}_e
  \lesssim h_e^{-1}\norm{\trc_e\delta\xi}_{-1,e}
  \lesssim h_e^{-1}h_f^{1/2}
    \norm{\trc\delta\xi}_{-1/2,\partial f}.
\]
Conditions {\rm (M1)} and {\rm (M2)} give $h_e\eqsim h_f$
and uniformly bound the number of facets.
Summing over $e$ therefore yields
\[
  \norm{\trc\delta\xi}_{\partial f}^2
  =\sum_{e\in\Delta_{d-1}(f)}\norm{\trc_e\delta\xi}_e^2
  \lesssim h_f^{-1}\norm{\trc\delta\xi}_{-1/2,\partial f}^2.
\]
The right-hand inequality of \eqref{eq:weak-trace} and
\eqref{eq:induction-delta-xi} therefore imply
\begin{equation}
\label{eq:delta-xi-trace}
  \norm{\trc\delta\xi}_{\partial f}
  \lesssim h_f^{-1/2}
    \bigl(\norm{\delta\xi}_f+h_f\norm{\rd\delta\xi}_f\bigr)
  \lesssim h_f^{-1/2}\norm{\delta\xi}_f.
\end{equation}
Since $\rd\xi=0$ and $\trc\xi=0$, the multiplicative trace
estimate \eqref{eq:cell-div-curl-trace} gives
\begin{equation}
\label{eq:negative-xi-normal-trace}
  \norm{\ntr\xi}_{\partial f}
  \lesssim h_f^{-1/2}\norm{\xi}_f
  +\norm{\xi}_f^{1/2}\norm{\delta\xi}_f^{1/2}.
\end{equation}
Integration by parts and
\eqref{eq:induction-delta-xi}--\eqref{eq:negative-xi-normal-trace} yield
\[
\begin{aligned}
  \norm{\delta\xi}_f^2
  &=\ip{\rd\delta\xi}{\xi}_f
    -\ip{\trc\delta\xi}{\ntr\xi}_{\partial f}\\
  &\lesssim h_f^{-1}\norm{\delta\xi}_f\norm{\xi}_f
    +h_f^{-1/2}\norm{\delta\xi}_f^{3/2}\norm{\xi}_f^{1/2}.
\end{aligned}
\]
Applying Young's inequality to the two terms gives
\[
\begin{aligned}
  \norm{\delta\xi}_f^2
  &\le\tfrac14\norm{\delta\xi}_f^2
    +C\bigl(h_f^{-1}\norm{\xi}_f\bigr)^2
    +\tfrac14\bigl(\norm{\delta\xi}_f^{3/2}\bigr)^{4/3}
    +C\bigl(h_f^{-1/2}\norm{\xi}_f^{1/2}\bigr)^4\\
  &\le\tfrac12\norm{\delta\xi}_f^2+C h_f^{-2}\norm{\xi}_f^2.
\end{aligned}
\]
Absorbing the first term on the right proves
\begin{equation}
\label{eq:delta-xi-by-xi}
  \norm{\delta\xi}_f\lesssim h_f^{-1}\norm{\xi}_f.
\end{equation}

\emph{Step 4: testing with a regular potential.}
Since $\rd\xi=0$ and $\trc\xi=0$, with $\int_f\xi=0$
when $k=d-1$, the regular-potential result
\cite[Proposition~4.1(ii)]{CostabelMcIntosh2010}, together with
the uniform operator bound \eqref{eq:uniform-bogovskii-bound}, gives
$\zeta:=\mathcal T_f^{k+1}\xi\in H_0^1\Lambda^k(f)$ satisfying
\[
  \rd\zeta=\xi,
  \qquad
  h_f^{-1}\norm{\zeta}_f+\seminorm{\zeta}_{1,f}
  \lesssim\norm{\xi}_f.
\]
Using $\trc\zeta=0$, the integration-by-parts formula
\eqref{eq:green-formula}, the Hodge decomposition
\eqref{eq:negative-inverse-decomposition}, duality, and
\eqref{eq:exact-negative-bound}, we obtain
\[
\begin{aligned}
  \norm{\xi}_f^2
  &=\ip{\rd\zeta}{\xi}_f
   =\ip{\zeta}{\delta\xi}_f
   =\ip{\zeta}{\omega_h}_f-\ip{\zeta}{\rd\phi}_f\\
  &\le\norm{\omega_h}_{-1,f}\norm{\zeta}_{1,f}
    +\norm{\rd\phi}_f\norm{\zeta}_f\\
  &\lesssim
    \bigl(\norm{\omega_h}_{-1,f}+h_f\norm{\rd\phi}_f\bigr)
    \norm{\xi}_f
   \lesssim\norm{\omega_h}_{-1,f}\norm{\xi}_f.
\end{aligned}
\]
Thus $\norm{\xi}_f\lesssim\norm{\omega_h}_{-1,f}$.
Together with \eqref{eq:delta-xi-by-xi}, this gives
\begin{equation}
\label{eq:coexact-negative-inverse}
  \norm{\delta\xi}_f\lesssim h_f^{-1}\norm{\omega_h}_{-1,f}.
\end{equation}
Finally, the triangle inequality, \eqref{eq:exact-negative-bound},
and \eqref{eq:coexact-negative-inverse} give
\[
  \norm{\omega_h}_f
  \le\norm{\rd\phi}_f+\norm{\delta\xi}_f
  \lesssim h_f^{-1}\norm{\omega_h}_{-1,f},
\]
which completes the induction and proves \eqref{eq:negative-inverse}.
\end{proof}

\begin{corollary}
[Differential and trace inverse estimates]
\label{cor:local-differential-trace-inverse}
Under conditions {\rm (M1)} and {\rm (M2)}, let $0\le k\le d-1$.
For every
$\omega_h\in
\Vcal_{r+1}^{-,k}(f)$,
\begin{equation}
\label{eq:local-differential-trace-inverse}
  h_f\norm{\rd\omega_h}_f
  +
  h_f^{1/2}
  \norm{\trc\omega_h}_{\partial f}
  \lesssim
  \norm{\omega_h}_f .
\end{equation}
\end{corollary}
\begin{proof}
For $d\ge2$, apply \eqref{eq:negative-inverse} to $\rd\omega_h$
and to each facet trace, and repeat the arguments for
\eqref{eq:induction-delta-xi}--\eqref{eq:delta-xi-trace}.
For $d=1$, $k=0$ and $\omega_h\in\Pcal_{r+1}\Lambda^0(f)$,
so the polynomial inverse and trace estimates
\eqref{eq:polynomial-estimates} apply.
\end{proof}

\subsection{Norm equivalences}
\label{subsec:appendix-norm-equivalences}
We first establish a local upper bound, then use recursion over subcells
and the inverse estimates to derive norm equivalences for both virtual
families.
Corresponding three-dimensional results for the trimmed family are established in
\cite{BeiraoMascottoMeng2022}, while the scalar $H^m$-conforming case
in arbitrary dimension is treated in \cite{ChenHuangWei2022}.

\begin{lemma}[Local upper bound]
\label{lem:upper}
Under condition {\rm (M1)}, let $0\le k<d$.
Every $\omega_h\in\Vcal_{r+1}^{-,k}(f)$ satisfies
\begin{equation}
\label{eq:trimmed-upper}
  \norm{\omega_h}_f
  \lesssim h_f\norm{\rd\omega_h}_f
  +\norm{\Pi_{r+1,f}^{\kappa,d-k}(\star\omega_h)}_f
  +h_f^{1/2}\norm{\trc\omega_h}_{\partial f}.
\end{equation}
\end{lemma}

\begin{proof}
Write $\omega_h=\rd\phi+\delta\xi$ using the orthogonal Hodge
decomposition of Lemma~\ref{lem:relative-hodge-decomposition}.
Since $\rd\xi=0$ and $\trc\xi=0$, the scaled Poincar\'e estimate
\eqref{eq:cell-poincare} and additive trace estimate
\eqref{eq:cell-div-curl-trace-additive} give
\[
  \norm{\xi}_f+h_f^{1/2}\norm{\ntr\xi}_{\partial f}
  \lesssim h_f\norm{\delta\xi}_f.
\]
The integration-by-parts formula \eqref{eq:green-formula} and
\eqref{eq:relative-coexact-potential} therefore yield
\begin{equation}
\label{eq:coexact}
\begin{aligned}
  \norm{\delta\xi}_f^2
  &=\ip{\rd\omega_h}{\xi}_f
    -\ip{\trc\omega_h}{\ntr\xi}_{\partial f}\\
  &\lesssim
    \bigl(h_f\norm{\rd\omega_h}_f
    +h_f^{1/2}\norm{\trc\omega_h}_{\partial f}\bigr)
    \norm{\delta\xi}_f.
\end{aligned}
\end{equation}
This proves the claim when $k=0$.

For $k\ge1$, the definition of $\Vcal_{r+1}^{-,k}(f)$ in
\eqref{eq:trimmed-codim-one} and \eqref{eq:trimmed-general} implies
that there exists $\mu\in\Kcal_{r+1}^{d-k}(f)$ such that
$\delta\omega_h=\delta\star^{-1}\mu$.
By \eqref{eq:relative-exact-potential} and scaled finite-dimensional
norm equivalence on $\Kcal_{r+1}^{d-k}(f)$, we obtain
\begin{equation}
\label{eq:upper-mu-control}
  \delta\omega_h=\delta\rd\phi\in\Pcal_r\Lambda^{k-1}(f),
  \qquad
  \norm{\mu}_f\lesssim h_f\norm{\delta\omega_h}_f.
\end{equation}
Applying the polynomial negative-norm inverse estimate
\eqref{eq:polynomial-negative-inverse} and the differentiation bound
\eqref{eq:negative-differentiation}, both with $j=1$, then gives
\begin{equation}
\label{eq:upper-codifferential-control}
  h_f\norm{\delta\omega_h}_f
  \lesssim\norm{\delta\omega_h}_{-1,f}
  \lesssim\norm{\rd\phi}_f.
\end{equation}
Using $\trc\phi=0$, integration by parts, the Hodge decomposition,
the projection identity \eqref{eq:hodge-projection-moments}, and
\eqref{eq:upper-mu-control}--\eqref{eq:upper-codifferential-control},
we obtain
\[
\begin{aligned}
  \norm{\rd\phi}_f^2
  &=\ip{\phi}{\delta\rd\phi}_f
   =\int_f(\omega_h-\delta\xi)\wedge\mu\\
  &=\int_f\star^{-1}
      \Pi_{r+1,f}^{\kappa,d-k}(\star\omega_h)\wedge\mu
    -\int_f\delta\xi\wedge\mu\\
  &\lesssim
    \bigl(\norm{\Pi_{r+1,f}^{\kappa,d-k}(\star\omega_h)}_f
    +\norm{\delta\xi}_f\bigr)\norm{\rd\phi}_f.
\end{aligned}
\]
Combining this estimate with \eqref{eq:coexact} and the
triangle inequality proves \eqref{eq:trimmed-upper}.
\end{proof}

\begin{theorem}[Norm equivalences]
\label{thm:norm-equiv}
\label{thm:appendix-trimmed-norm}
Under conditions {\rm (M1)} and {\rm (M2)}, let $T\in\Th$
and $0\le k\le n$. The following norm equivalences hold.
\begin{enumerate}
\renewcommand{\theenumi}{\roman{enumi}}
\renewcommand{\labelenumi}{\textup{(\theenumi)}}
\item For every $\omega_h\in\Vcal_{r+1}^{-,k}(T)$,
\begin{equation}
\label{eq:trimmed-global-upper}
\begin{aligned}
  \norm{\omega_h}_T^2
  &\eqsim\sum_{f\in\Delta_k(T)}
    h_T^{n-k}\norm{\trc_f\omega_h}_f^2\\
  &\quad+\sum_{f\in\Delta_{k+1}(T)}
    \left(h_T^{n-k+1}\norm{\trc_f\rd\omega_h}_f^2
    +h_T^{n-k-1}
    \norm{\Pi_{r+1,f}^{\kappa,1}(\star\trc_f\omega_h)}_f^2\right)\\
  &\quad+\sum_{d=k+2}^{n}\sum_{f\in\Delta_d(T)}
    \left(\begin{aligned}
      &h_T^{n-d+2}
       \norm{\Pi_{r+1,f}^{\kappa,d-k-1}(\star\trc_f\rd\omega_h)}_f^2\\
      &{}+h_T^{n-d}
       \norm{\Pi_{r+1,f}^{\kappa,d-k}(\star\trc_f\omega_h)}_f^2
    \end{aligned}\right).
\end{aligned}
\end{equation}
\item For every $\omega_h\in\Vcal_r^k(T)$,
\begin{equation}
\label{eq:complete-upper}
  \norm{\omega_h}_T^2
  \eqsim\sum_{d=k}^{n}\sum_{f\in\Delta_d(T)}
    h_T^{n-d}
    \norm{\Pi_{r,f}^{-,d-k}(\star\trc_f\omega_h)}_f^2.
\end{equation}
\end{enumerate}
Empty sums are zero, and both right-hand sides are computable from the
respective DOFs. The hidden constants depend only on $n$, $r$, $\rho$,
and $\gamma$.
\end{theorem}

\begin{proof}
For $k=n$, both equivalences are identities since
$\Vcal_r^n(T)=\Vcal_{r+1}^{-,n}(T)=\Pcal_r\Lambda^n(T)$.
Henceforth, let $k<n$.

\emph{Local bounds for the trimmed family.}
Let $\omega\in\Vcal_{r+1}^{-,k}(f)$, with $d>k$.
For $d\ge k+2$, applying Lemma~\ref{lem:upper} to the closed form
$\rd \omega$ gives
\[
  h_f\norm{\rd \omega}_f
  \lesssim
  h_f\norm{\Pi_{r+1,f}^{\kappa,d-k-1}(\star\rd \omega)}_f
  +h_f^{3/2}\norm{\trc\rd \omega}_{\partial f}.
\]
Trace commutativity, Corollary~\ref{cor:local-differential-trace-inverse}
on each facet $e$ of $f$, and $h_e\eqsim h_f$ from {\rm (M2)} yield
$h_f\norm{\trc\rd \omega}_{\partial f}
\lesssim\norm{\trc \omega}_{\partial f}$.
Combining these estimates with Lemma~\ref{lem:upper} applied to $\omega$ gives
\begin{equation}
\label{eq:trimmed-local-norm-upper}
\begin{aligned}
  \norm{\omega}_f\lesssim{}
    h_f\norm{\Pi_{r+1,f}^{\kappa,d-k-1}(\star\rd \omega)}_f
  +\norm{\Pi_{r+1,f}^{\kappa,d-k}(\star \omega)}_f
    +h_f^{1/2}\norm{\trc \omega}_{\partial f}.
\end{aligned}
\end{equation}
For $d=k+1$, \eqref{eq:trimmed-local-norm-upper} holds with the first term replaced by
$h_f\norm{\rd \omega}_f$.

\emph{Local bounds for the complete family.}
Let $\omega\in\Vcal_r^k(f)$, with $d>k$.
The defining constraints \eqref{eq:complete-space} and
the projection identities
\eqref{eq:l2-orthogonal-projection}--\eqref{eq:hodge-projection-moments} give
\begin{equation*}
\begin{aligned}
  \Pi_{r+1,f}^{\kappa,d-k}(\star \omega)
  &=\Pi_{r+1,f}^{\kappa,d-k}(\star\Proj{r}{f}{k}\omega),\\
  \Pi_{r+1,f}^{\kappa,d-k-1}(\star\rd \omega)
  &=\Pi_{r+1,f}^{\kappa,d-k-1}(\star\Proj{r}{f}{k+1}\rd \omega).
\end{aligned}
\end{equation*}
When $d=k+1$, the second identity follows directly from
$\rd \omega\in\Pcal_r\Lambda^d(f)$.
Thus, the local trimmed bound \eqref{eq:trimmed-local-norm-upper},
with the stated modification for $d=k+1$, gives
\begin{equation}
\label{eq:complete-local-projected-bound}
  \norm{\omega}_f
  \lesssim h_f\norm{\Proj{r}{f}{k+1}\rd \omega}_f
  +\norm{\Proj{r}{f}{k}\omega}_f
  +h_f^{1/2}\norm{\trc \omega}_{\partial f}.
\end{equation}
For $\eta\in\Pcal_r\Lambda^{d-k-1}(f)$, rewriting
\eqref{eq:projection-computability-derivative} using
the projection identity \eqref{eq:hodge-projection-moments} gives
\[
\begin{aligned}
  \int_f\Proj{r}{f}{k+1}\rd \omega\wedge\eta
  &=(-1)^{k+1}\int_f
    \star^{-1}\Pi_{r,f}^{-,d-k}(\star \omega)\wedge\rd\eta
    +\int_{\partial f}\trc \omega\wedge\trc\eta\\
  &\le\norm{\Pi_{r,f}^{-,d-k}(\star \omega)}_f\norm{\rd\eta}_f
    +\norm{\trc \omega}_{\partial f}\norm{\trc\eta}_{\partial f}.
\end{aligned}
\]
Using the polynomial inverse and trace estimates
\eqref{eq:polynomial-estimates}, we divide by $\norm{\eta}_f$ and take
the supremum over nonzero $\eta$ to obtain
\begin{equation}
\label{eq:complete-local-derivative-projection-bound}
  h_f\norm{\Proj{r}{f}{k+1}\rd \omega}_f
  \lesssim\norm{\Pi_{r,f}^{-,d-k}(\star \omega)}_f
  +h_f^{1/2}\norm{\trc \omega}_{\partial f}.
\end{equation}
To estimate $\norm{\Proj{r}{f}{k}\omega}_f$, use the isomorphism
\eqref{eq:polynomial-koszul-isomorphism} to write any
$\eta\in\Pcal_r\Lambda^{d-k}(f)$ uniquely as $\eta=\rd\mu+\nu$,
with $\mu\in\Kcal_{r+1}^{d-k-1}(f)$ and $\nu\in\Kcal_r^{d-k}(f)$.
Rewriting \eqref{eq:projection-computability-stokes} using
the projection identity \eqref{eq:hodge-projection-moments} then gives
\begin{equation}
\label{eq:complete-local-moment-bound}
\begin{aligned}
  (-1)^{k+1}\int_f\Proj{r}{f}{k}\omega\wedge\eta
  &=\int_f\Proj{r}{f}{k+1}\rd \omega\wedge\mu
    -\int_{\partial f}\trc \omega\wedge\trc\mu\\
  &\quad+(-1)^{k+1}\int_f
    \star^{-1}\Pi_{r,f}^{-,d-k}(\star \omega)\wedge\nu\\
  &\le\norm{\Proj{r}{f}{k+1}\rd \omega}_f\norm{\mu}_f
    +\norm{\trc \omega}_{\partial f}\norm{\trc\mu}_{\partial f}\\
  &\quad+\norm{\Pi_{r,f}^{-,d-k}(\star \omega)}_f\norm{\nu}_f.
\end{aligned}
\end{equation}
Under {\rm (M1)}, scaling and polynomial norm comparison on inner and
outer balls extend the finite-dimensional arguments
in~\cite[Lemmas~2 and~9]{DiPietroDroniou2023} to the polynomial
isomorphism \eqref{eq:polynomial-koszul-isomorphism} in any dimension,
yielding
\begin{equation}
\label{eq:scaled-polynomial-koszul-control}
  h_f^{-1}\norm{\mu}_f+\norm{\nu}_f
  \lesssim\norm{\rd\mu}_f+\norm{\nu}_f
  \lesssim\norm{\eta}_f.
\end{equation}
The hidden constants depend only on $n$, $r$, and $\rho$.
The polynomial trace estimate \eqref{eq:polynomial-estimates} also gives
$h_f^{-1/2}\norm{\trc\mu}_{\partial f}\lesssim\norm{\eta}_f$.
Using \eqref{eq:scaled-polynomial-koszul-control} and the trace bound,
we also divide \eqref{eq:complete-local-moment-bound} by $\norm{\eta}_f$
and take the supremum over nonzero $\eta$ to obtain
\begin{equation}
\label{eq:complete-local-form-projection-bound}
  \norm{\Proj{r}{f}{k}\omega}_f
  \lesssim h_f\norm{\Proj{r}{f}{k+1}\rd \omega}_f
  +\norm{\Pi_{r,f}^{-,d-k}(\star \omega)}_f
  +h_f^{1/2}\norm{\trc \omega}_{\partial f}.
\end{equation}
Combining \eqref{eq:complete-local-derivative-projection-bound} and
\eqref{eq:complete-local-form-projection-bound} with
\eqref{eq:complete-local-projected-bound} proves
\begin{equation}
\label{eq:complete-local-norm-upper}
  \norm{\omega}_f
  \lesssim\norm{\Pi_{r,f}^{-,d-k}(\star \omega)}_f
  +h_f^{1/2}\norm{\trc \omega}_{\partial f}.
\end{equation}

\emph{Recursion and reverse estimates.}
Conditions {\rm (M1)}--{\rm (M2)} ensure $h_f\eqsim h_T$ and
uniformly bounded multiplicities.
For the complete family, apply \eqref{eq:complete-local-norm-upper}
to $\omega=\trc_f\omega_h$, square the bound, and multiply by $h_T^{n-d}$.
Starting from $f=T$ and repeating this on its facets gives, for $k\le n-2$,
\begin{equation*}
\begin{aligned}
  \norm{\omega_h}_T^2
  &\lesssim\norm{\Pi_{r,T}^{-,n-k}(\star\omega_h)}_T^2
    +h_T\sum_{f\in\Delta_{n-1}(T)}\norm{\trc_f\omega_h}_f^2\\
  &\lesssim\sum_{d=n-1}^{n}\sum_{f\in\Delta_d(T)}
    h_T^{n-d}\norm{\Pi_{r,f}^{-,d-k}(\star\trc_f\omega_h)}_f^2
    +h_T^2\sum_{g\in\Delta_{n-2}(T)}\norm{\trc_g\omega_h}_g^2\\
  &\lesssim\cdots\lesssim
    \sum_{d=k+1}^{n}\sum_{f\in\Delta_d(T)}
    h_T^{n-d}\norm{\Pi_{r,f}^{-,d-k}(\star\trc_f\omega_h)}_f^2
    +h_T^{n-k}\sum_{f\in\Delta_k(T)}\norm{\trc_f\omega_h}_f^2\\
  &=\sum_{d=k}^{n}\sum_{f\in\Delta_d(T)}
    h_T^{n-d}\norm{\Pi_{r,f}^{-,d-k}(\star\trc_f\omega_h)}_f^2.
\end{aligned}
\end{equation*}
The last equality uses $\trc_f\omega_h\in\Pcal_r\Lambda^k(f)$
on $k$-cells. For $k=n-1$, the recursion stops after the first step.
This proves the upper bound in \eqref{eq:complete-upper}.
For the trimmed family, we apply \eqref{eq:trimmed-local-norm-upper}
in the same way, using its modified form for $d=k+1$ and stopping
at $k$-cells, to obtain the upper bound in \eqref{eq:trimmed-global-upper}.

For the reverse bounds, iterating the trace inverse estimate and
using the differential inverse estimate of
Corollary~\ref{cor:local-differential-trace-inverse} gives
\begin{equation}
\label{eq:iterated-virtual-inverse}
\begin{aligned}
  h_T^{(n-d)/2}\norm{\trc_f\omega_h}_f
  &\lesssim\norm{\omega_h}_T,
  &&k\le d\le n,\\
  h_T^{(n-d+2)/2}\norm{\trc_f\rd\omega_h}_f
  &\lesssim\norm{\omega_h}_T,
  &&k+1\le d\le n.
\end{aligned}
\end{equation}
These estimates apply to both families by
$\Vcal_r^k(T)\subset\Vcal_{r+1}^{-,k}(T)$.
Projection stability and the uniformly bounded number of subcells of $T$
complete the proof.
\end{proof}

\begin{remark}[Analytical bridge between polytopal complexes]
\label{rem:continuous-data-norm-bridge}
As discussed in \cite[Section~2.4]{BeiraoDiPietroDroniou2026}, norm equivalence plays
a key role in linking conforming VEM and fully discrete DDR analyses.
Here we make this link explicit for our complete complex
$\Vcal_{r,h}^{\bullet}$ and the DDR-inspired complex
$\underline{\Vcal}_{r,h}^{\bullet}$.
Following \cite[Remark~2]{pietro_uniform_2025}, we equip the local
data space with the scaled DDR norm
\begin{equation}
\label{eq:complete-local-data-norm}
  \norm{\underline\omega_T}_{\mathrm{data},k,T}^2
  :=\sum_{d=k}^{n}\sum_{f\in\Delta_d(T)}
    h_T^{n-d}\norm{\omega_f}_f^2,
  \qquad\underline\omega_T=(\omega_f)_f.
\end{equation}
Since the DOF components are precisely the projected moments in
\eqref{eq:complete-upper}, Theorem~\ref{thm:norm-equiv} gives
\begin{equation}
\label{eq:continuous-data-norm-equivalence}
  \norm{\mathcal J_{r,T}^k\omega_h}_{\mathrm{data},k,T}
  \eqsim\norm{\omega_h}_T,
  \qquad\omega_h\in\Vcal_r^k(T),
\end{equation}
where $\mathcal J_{r,T}^k$ collects the DOF components
\eqref{eq:global-complete-isomorphism} on the subcells of $T$.
Together with the cochain identities
\eqref{eq:fully-discrete-commuting-diagrams}, this equivalence allows
estimates to be transferred between the two formulations.
For example, this equivalence allows us to recover the discrete Poincar\'e
inequality of Theorem~\ref{thm:complete-orthogonal-poincare} from its
fully discrete counterpart (Appendix~\ref{subsec:appendix-analytical-transfer}),
transfer low-regularity approximation estimates from the fully discrete
setting to VEM (Remark~\ref{rem:low-regularity-quasi-interpolation}), and use
the inverse DOF map as a stable conforming lifting in fully discrete error
analysis (Appendix~\ref{subsec:appendix-analytical-transfer}).

\end{remark}

\begin{remark}[Relation to the FES interpretation of DDR]
\label{rem:fes-ddr-interpretation-comparison}
Reference~\cite{ChristiansenRapetti2025} provides a conforming FES
interpretation of the DDR space $\underline{\Vcal}_{r,h}^{k}$,
sharing the moment DOFs of our complete space $\Vcal_{r,h}^{k}$.
Its local spaces are defined through weak formulations with auxiliary
Lagrange multipliers and constrained test
spaces~\cite[Eqs.~(4.15)--(4.21)]{ChristiansenRapetti2025}, whereas our
local complete spaces inherit the trimmed family's explicit polynomial
conditions on $\delta\rd\omega$ and $\delta\omega$ under additional
enhancement constraints.
Since we establish the norm equivalence directly for our VEM spaces,
we do not investigate whether they coincide with the local FES spaces.

The main analytical challenge is uniform control of the norm-equivalence
constants.
In~\cite[Eqs.~(5.3)--(5.5)]{ChristiansenRapetti2025}, fixed-cell
equivalence follows from unisolvence without quantifying the constants'
geometric dependence; mesh-uniformity relies on the conjectural Hypothesis of
Stability of Constants (HSC), concerning stability under small geometric
perturbations~\cite[Remark~2.1 and Section~5]{ChristiansenRapetti2025}.
Using the explicit polynomial differential structure, we establish
Theorem~\ref{thm:norm-equiv} under (M1)--(M2), with constants
depending only on $\rho$ and $\gamma$ for fixed dimension and
polynomial degree.
\end{remark}

\begin{corollary}[Stability of the cochain retraction]
\label{cor:cochain-retraction-stability}
Under conditions {\rm (M1)} and {\rm (M2)}, let $r\ge0$ and
$0\le k\le n$.
The cochain retraction $P_{r,h}^k$ defined by
\eqref{eq:complete-moment-reduction} satisfies
\begin{equation}
\label{eq:cochain-retraction-stability}
  \norm{P_{r,h}^k \omega_h^-}_T\lesssim\norm{\omega_h^-}_T,
  \qquad
  \norm{P_{r,h}^k \omega_h^-}\lesssim\norm{\omega_h^-},
\end{equation}
for every $\omega_h^-\in\Vcal_{r+1,h}^{-,k}$ and $T\in\Th$.
\end{corollary}

\begin{proof}
Fix $T\in\Th$ and $\omega_h^-\in\Vcal_{r+1,h}^{-,k}$, and set
$\omega_h=P_{r,h}^k\omega_h^-$. Writing the defining moment conditions
\eqref{eq:complete-moment-reduction} in terms of the $L^2$ projections
as in \eqref{eq:hodge-projection-moments} gives
\[
  \Pi_{r,f}^{-,d-k}(\star\trc_f\omega_h)
  =\Pi_{r,f}^{-,d-k}(\star\trc_f\omega_h^-)
  \qquad\text{for every }f\in\Delta_d(T),\quad k\le d\le n.
\]
Since $\omega_h|_T\in\Vcal_r^k(T)$, we can apply \eqref{eq:complete-upper}
to $\omega_h$. Using the equality above and then $L^2$-projection stability, we obtain
\[
\begin{aligned}
  \norm{\omega_h}_T^2
  &\lesssim\sum_{d=k}^{n}\sum_{f\in\Delta_d(T)}h_T^{n-d}
    \norm{\Pi_{r,f}^{-,d-k}(\star\trc_f\omega_h^-)}_f^2\\
  &\le\sum_{d=k}^{n}\sum_{f\in\Delta_d(T)}h_T^{n-d}
    \norm{\trc_f\omega_h^-}_f^2
  \lesssim\norm{\omega_h^-}_T^2.
\end{aligned}
\]
The last inequality follows from the first trace estimate in \eqref{eq:iterated-virtual-inverse}
and the uniformly bounded number of subcells of $T$.
Summing the squared local bounds over $T\in\Th$ gives the global estimate.
\end{proof}

\subsection{Canonical interpolation estimates}
\label{subsec:section4-interpolation}
The interpolation estimates follow from the upper bounds in the norm
equivalences, scaled Sobolev trace estimates, and commuting polynomial approximation.

\begin{theorem}[Canonical interpolation estimates]
\label{thm:section4-interpolation}
\label{thm:appendix-trimmed-interpolation}
Under conditions {\rm (M1)} and {\rm (M2)}, let $T\in\Th$,
$0\le k\le n$, and $s>(n-k)/2$.
For every $\omega\in H^s\Lambda^k(\Omega)$, the following estimates hold.
\begin{enumerate}
\renewcommand{\theenumi}{\roman{enumi}}
\renewcommand{\labelenumi}{\textup{(\theenumi)}}
\item The complete interpolant satisfies
\begin{equation}
\label{eq:section4-local-interpolation}
  \norm{\omega-\Ih{k}\omega}_T
  \lesssim h_T^{\min\{s,r+1\}}\norm{\omega}_{s,T}.
\end{equation}
\item If $k<n$ and $\rd\omega\in H^t\Lambda^{k+1}(\Omega)$
for some $t>(n-k-1)/2$, the trimmed interpolant satisfies
\begin{equation}
\label{eq:section4-local-trimmed-interpolation}
  \norm{\omega-\mathcal I_{r+1,h}^{-,k}\omega}_T
  \lesssim h_T^{\min\{s,r+1\}}\norm{\omega}_{s,T}
  +h_T^{\min\{t+1,r+1\}}\norm{\rd\omega}_{t,T}.
\end{equation}
For $k=n$, it satisfies
\[
  \norm{\omega-\mathcal I_{r+1,h}^{-,n}\omega}_T
  \lesssim h_T^{\min\{s,r+1\}}\norm{\omega}_{s,T}.
\]
\end{enumerate}
Under the additional assumptions in {\rm (ii)}, both interpolants satisfy
\begin{align}
\label{eq:section4-derivative-interpolation}
  \norm{\rd\omega-\rd\Ih{k}\omega}_T
  &\lesssim h_T^{\min\{t,r+1\}}\norm{\rd\omega}_{t,T},\\
\label{eq:section4-trimmed-derivative-interpolation}
  \norm{\rd\omega-\rd\mathcal I_{r+1,h}^{-,k}\omega}_T
  &\lesssim h_T^{\min\{t,r+1\}}\norm{\rd\omega}_{t,T}.
\end{align}
\end{theorem}

\begin{proof}
For $k=n$, both interpolants coincide with the $L^2$ projection onto
$\Pcal_r\Lambda^n(T)$, and the estimates follow from
polynomial approximation \eqref{eq:section4-poly-approx}.

Let $k<n$, set
$\omega_{r,T}^{\pi}:=Q_{r,T}^k\omega$
as in Lemma~\ref{lem:polynomial-approximation},
and write
$\zeta_T:=\omega-\omega_{r,T}^{\pi}$.
The polynomial inclusions \textup{(T1)} and \textup{(C1)},
the upper bounds in the norm equivalences
\eqref{eq:trimmed-global-upper}--\eqref{eq:complete-upper},
the interpolation and commuting properties, and $L^2$-projection stability give
\begin{equation}
\label{eq:canonical-interpolation-data-bound}
\begin{aligned}
  \norm{\mathcal I_{r+1,h}^{-,k}\omega-\omega_{r,T}^{\pi}}_T
  =\norm{\mathcal I_{r+1,h}^{-,k}\zeta_T}_T
  &\lesssim
  \sum_{d=k}^{n}\sum_{f\in\Delta_d(T)}
  h_T^{(n-d)/2}\norm{\trc_f\zeta_T}_f\\
  &\quad+\sum_{d=k+1}^{n}\sum_{f\in\Delta_d(T)}
  h_T^{(n-d+2)/2}\norm{\trc_f\rd\zeta_T}_f.
\end{aligned}
\end{equation}
The corresponding estimate for the complete interpolant contains only
the first sum.

To estimate these traces, choose a chain of cells
$T=f_n\supset f_{n-1}\supset\cdots\supset f_d=f$,
and set $q:=n-d$.
Condition {\rm (M2)} gives $\gamma h_T\le h_{f_j}\le h_T$
for each positive-dimensional cell in the chain, so all scales
can be replaced by $h_T$.
Let $\zeta$ be a differential form on $T$ with Sobolev regularity
$\beta>q/2$.
Each application of the trace estimate lowers the available Sobolev order
by $1/2$, so $\beta>q/2$ ensures that
\eqref{eq:facet-sobolev-trace} applies at every step.
For $q\ge2$, iterating \eqref{eq:facet-sobolev-trace} while retaining
the scaled seminorms gives\nopagebreak
\begin{equation}
\label{eq:trace-iterated}
\begin{aligned}
  \norm{\trc_f \zeta}_f
  &\lesssim h_T^{-1/2}\bigl(
    \norm{\trc_{f_{d+1}}\zeta}_{f_{d+1}}
    +h_T^{\beta-(q-1)/2}
      \seminorm{\trc_{f_{d+1}}\zeta}_{\beta-(q-1)/2,f_{d+1}}\bigr)\\
  &\lesssim h_T^{-1}\bigl(
    \norm{\trc_{f_{d+2}}\zeta}_{f_{d+2}}
    +h_T^{\beta-(q-2)/2}
      \seminorm{\trc_{f_{d+2}}\zeta}_{\beta-(q-2)/2,f_{d+2}}\bigr)\\
  &\lesssim\cdots\lesssim h_T^{-q/2}
    \bigl(\norm{\zeta}_T+h_T^\beta\seminorm{\zeta}_{\beta,T}\bigr)
  =h_T^{-q/2}\norm{\zeta}_T
    +h_T^{\beta-q/2}\seminorm{\zeta}_{\beta,T}.
\end{aligned}
\end{equation}
For $q=1$, the same final bound follows from a single application of
\eqref{eq:facet-sobolev-trace}; for $q=0$, it is immediate.
When $d=0$, the step from $f_1$ to $f_0$ uses the endpoint version of
\eqref{eq:facet-sobolev-trace}, with the left-hand side interpreted
as the absolute value at the vertex.
The constants may depend on $\gamma$ as well as $\rho$.
Estimate \eqref{eq:trace-iterated} and the
uniformly bounded number of subcells of $T$ bound the two sums in
\eqref{eq:canonical-interpolation-data-bound} by, respectively,
\[
  \norm{\zeta_T}_T+h_T^s\seminorm{\zeta_T}_{s,T},
  \qquad
  h_T\norm{\rd\zeta_T}_T+h_T^{t+1}\seminorm{\rd\zeta_T}_{t,T}.
\]
Combining these trace bounds with the triangle inequality and
the polynomial approximation estimates
\eqref{eq:section4-poly-approx}--\eqref{eq:section4-poly-approx-derivative}
gives
\[
\begin{aligned}
  \norm{\omega-\mathcal I_{r+1,h}^{-,k}\omega}_T
  &\le \norm{\zeta_T}_T
       +\norm{\mathcal I_{r+1,h}^{-,k}\zeta_T}_T\\
  &\lesssim \norm{\zeta_T}_T+h_T^s\seminorm{\zeta_T}_{s,T}
       +h_T\norm{\rd\zeta_T}_T
       +h_T^{t+1}\seminorm{\rd\zeta_T}_{t,T}\\
  &\lesssim h_T^{\min\{s,r+1\}}\norm{\omega}_{s,T}
       +h_T^{\min\{t+1,r+1\}}\norm{\rd\omega}_{t,T}.
\end{aligned}
\]
For Sobolev orders above the polynomial degrees, the corresponding
polynomial seminorms vanish; when $r=0$, we use $\rd \omega_{r,T}^{\pi}=0$.
This proves \eqref{eq:section4-local-trimmed-interpolation}.
The complete estimate \eqref{eq:section4-local-interpolation}
follows in the same way, using only the first sum in
\eqref{eq:canonical-interpolation-data-bound}.

Finally, under the additional regularity assumption on $\rd\omega$,
the commuting identities \eqref{eq:virtual-interpolator-commutativity} give
\[
\begin{aligned}
  \rd\omega-\rd\Ih{k}\omega
  &=\rd\omega-\Ih{k+1}\rd\omega,\\
  \rd\omega-\rd\mathcal I_{r+1,h}^{-,k}\omega
  &=\rd\omega-\mathcal I_{r+1,h}^{-,k+1}\rd\omega.
\end{aligned}
\]
Applying the $L^2$ interpolation estimates
\eqref{eq:section4-local-interpolation}--\eqref{eq:section4-local-trimmed-interpolation}
in degree $k+1$ (using the top-degree bound in \textup{(ii)} when
$k+1=n$) to the closed form $\rd\omega$ proves
\eqref{eq:section4-derivative-interpolation}--\eqref{eq:section4-trimmed-derivative-interpolation}.
\end{proof}

\begin{remark}[Low-regularity approximation]
\label{rem:low-regularity-quasi-interpolation}
Motivated by the three-dimensional DDR construction in
\cite[(4.1)]{ChaumontFreletDroniouLemaire2025}, we compose the
Falk--Winther cochain projection
\cite[Theorems~4.7 and~5.2]{FalkWinther2014} on a suitable auxiliary
conforming simplicial refinement $\mathcal S_h$ with our complete
canonical interpolator. Set
$\mathcal W_{\mathcal S}^k:=H\Lambda^k(\rd;\Omega)
\cap\Pcal_{r+1}^-\Lambda^k(\mathcal S_h)$.
Extending the commuting canonical interpolator to this space using the same moment DOFs,
we define
\[
  \widetilde{\mathcal I}_{r,h}^k:=\Ih{k}\pi_{\mathcal S}^k:
  H\Lambda^k(\rd;\Omega)\longrightarrow\Vcal_{r,h}^k,
  \qquad
  \pi_{\mathcal S}^k:H\Lambda^k(\rd;\Omega)
  \longrightarrow\mathcal W_{\mathcal S}^k.
\]
For $k<n$, it satisfies
\[
  \rd\widetilde{\mathcal I}_{r,h}^k\omega
  =\Ih{k+1}\rd\pi_{\mathcal S}^k\omega
  =\Ih{k+1}\pi_{\mathcal S}^{k+1}\rd\omega
  =\widetilde{\mathcal I}_{r,h}^{k+1}\rd\omega.
\]
For $n=3$, choose $\mathcal S_h$ as in
Remark~\ref{rem:auxiliary-simplicial-submesh}, which satisfies the
submesh assumptions in~\cite[Section~2.1]{ChaumontFreletDroniouLemaire2025},
and write
$\underline{\mathcal I}_{r,T}^k:=\mathcal J_{r,T}^k\Ih{k}$.
By the interpolation bound in
\cite[Lemma~16]{ChaumontFreletDroniouLemaire2025} for $k<3$
and $L^2$-projection stability for $k=3$, we have
\[
  \bigl\|\underline{\mathcal I}_{r,T}^k\zeta_T
         \bigr\|_{\mathrm{data},k,T}
  \lesssim\norm{\zeta_T}_T,
  \qquad
  \zeta_T\in\left.\mathcal W_{\mathcal S}^k\right|_T.
\]
Set $\zeta_T:=\pi_{\mathcal S}^k\omega-\Proj{r}{T}{k}\omega$.
The polynomial inclusion \textup{(C1)}, the norm equivalence
\eqref{eq:continuous-data-norm-equivalence}, and the preceding bound yield
\begin{equation*}
\begin{aligned}
  \norm{\omega-\widetilde{\mathcal I}_{r,h}^k\omega}_T
  &\le\norm{\omega-\Proj{r}{T}{k}\omega}_T+\norm{\Ih{k}\zeta_T}_T\\
  &\lesssim\norm{\omega-\Proj{r}{T}{k}\omega}_T
       +\norm{\omega-\pi_{\mathcal S}^k\omega}_T.
\end{aligned}
\end{equation*}
Thus the local $L^2$ error is controlled by polynomial approximation
and the auxiliary finite element projection error.
For $k<n$, the corresponding interpolation error estimate for the
exterior derivative follows by applying the above bound to $\rd\omega$
and using the commuting property.
We make no attempt here to extend this approximation analysis
to arbitrary dimension and retain canonical interpolation in
Section~\ref{sec:hodge-laplace-vem}.
\end{remark}

\section{A mixed virtual element method for the Hodge Laplacian}
\label{sec:hodge-laplace-vem}
We consider the Hodge--Laplace problem
\eqref{eq:strong-hodge-laplace} on $\Omega$ with
$V^\bullet(\Omega)=H\Lambda^\bullet(\rd;\Omega)$.
For $j\in L^2\Lambda^k(\Omega)$, its equivalent mixed weak formulation seeks
$(\sigma,u,p)\in H\Lambda^{k-1}(\rd;\Omega)
\times H\Lambda^k(\rd;\Omega)\times\Harm^k(\Omega)$ such that
\begin{align}
\label{eq:continuous-mixed-1}
  \langle\sigma,\tau\rangle
  -\langle u,\rd\tau\rangle
  &=0
  &&\forall\tau\in H\Lambda^{k-1}(\rd;\Omega),
  \\
\label{eq:continuous-mixed-2}
  \langle\rd\sigma,v\rangle
  +\langle\rd u,\rd v\rangle
  +\langle p,v\rangle
  &=\langle j,v\rangle
  &&\forall v\in H\Lambda^k(\rd;\Omega),
  \\
\label{eq:continuous-mixed-harmonic}
  \langle u,q\rangle&=0
  &&\forall q\in\Harm^k(\Omega),
\end{align}
where
\begin{equation*}
  \Harm^k(\Omega)=\{\omega\in H\Lambda^k(\rd;\Omega):\rd \omega=0,\quad
    \langle \omega,\rd\tau\rangle=0
    \quad\forall\tau\in H\Lambda^{k-1}(\rd;\Omega)\}.
\end{equation*}
The mixed solution satisfies
$\sigma=\delta u$ and $p=\Pi_{\Harm^k(\Omega)}j$,
with $u$ solving \eqref{eq:strong-hodge-laplace}.
\par
Throughout this section, the polynomial degree $r$ is fixed and the
mesh regularity conditions in Assumption~\ref{ass:mesh-regularity}
are in force.
We focus on the stability and convergence of the discretization of
\eqref{eq:continuous-mixed-1}--\eqref{eq:continuous-mixed-harmonic}
based on the complete virtual complex
\eqref{eq:complete-virtual-derham-complex}, which reduces the number
of DOFs. The norm equivalence and interpolation estimates established
for the trimmed family in
Theorems~\ref{thm:norm-equiv} and~\ref{thm:section4-interpolation}
also support a parallel analysis.

\subsection{Discrete inner products and Poincar\'e inequalities}
\label{subsec:discrete-inner-products}
Following the standard VEM construction \cite{BeiraoEtAl2013},
for each $T\in\Th$, we define the local discrete $L^2$ inner product
on $\Vcal_r^k(T)$ by
\begin{equation}
\label{eq:discrete-bilinear}
[\omega_h,\zeta_h]_{k,T}
:=
\left\langle
  \Proj{r}{T}{k}\omega_h,
  \Proj{r}{T}{k}\zeta_h
\right\rangle_T
+S_T^k\bigl(
  (I-\Proj{r}{T}{k})\omega_h,
  (I-\Proj{r}{T}{k})\zeta_h
\bigr),
\end{equation}
where the stabilization is defined using the weights and projected
moments in \eqref{eq:complete-upper}:
\begin{equation}
\label{eq:stabilization}
S_T^k(\omega_h,\zeta_h)
:=
\sum_{d=k}^{n}
\sum_{f\in\Delta_d(T)}
h_T^{n-d}
\left\langle
  \Pi_{r,f}^{-,d-k}(\star\trc_f\omega_h),
  \Pi_{r,f}^{-,d-k}(\star\trc_f\zeta_h)
\right\rangle_f .
\end{equation}

The local discrete inner product \eqref{eq:discrete-bilinear} satisfies
the following two properties.
\begin{itemize}
\item[(i)] \emph{Polynomial consistency.}
For all $\omega_h\in\Vcal_r^k(T)$ and
$\mu\in\Pcal_r\Lambda^k(T)$,
\begin{equation}
\label{Consistency}
  [\omega_h,\mu]_{k,T}
  =\langle\omega_h,\mu\rangle_T.
\end{equation}
\item[(ii)] \emph{Uniform stability.}
For all $\omega_h\in\Vcal_r^k(T)$,
\begin{equation}
\label{Coercivity}
  \alpha_*\norm{\omega_h}_T^2
  \le [\omega_h,\omega_h]_{k,T}
  \le \alpha^*\norm{\omega_h}_T^2,
\end{equation}
where $0<\alpha_*\le\alpha^*$ are constants independent of $h$.
\end{itemize}
Polynomial consistency \eqref{Consistency} follows directly from
\eqref{eq:discrete-bilinear}, since $\Proj{r}{T}{k}$ is the
$L^2$-orthogonal projector onto $\Pcal_r\Lambda^k(T)$.
For stability, Theorem~\ref{thm:norm-equiv} gives
$S_T^k(\zeta_h,\zeta_h)\eqsim\norm{\zeta_h}_T^2$ for
$\zeta_h\in\ker\Proj{r}{T}{k}$. Applying this equivalence to
$\zeta_h=(I-\Proj{r}{T}{k})\omega_h$ and using the orthogonal splitting
induced by $\Proj{r}{T}{k}$ proves \eqref{Coercivity}.
Finally, the Cauchy--Schwarz inequality and the upper bound in
\eqref{Coercivity} yield
\begin{equation}
\label{Stability}
  |[\omega_h,\zeta_h]_{k,T}|
  \le\alpha^*\norm{\omega_h}_T\norm{\zeta_h}_T,
  \qquad \omega_h,\zeta_h\in\Vcal_r^k(T).
\end{equation}

Summing the local inner products over the mesh defines the global
discrete $L^2$ inner product on $\Vcal_{r,h}^k$:
\begin{equation}
\label{eq:global-discrete-bilinear}
  [\omega_h,\zeta_h]_{k,h}
  :=\sum_{T\in\Th}
  [\omega_h|_T,\zeta_h|_T]_{k,T},
\end{equation}
with induced norm $\norm{\omega_h}_h:=[\omega_h,\omega_h]_{k,h}^{1/2}$
for $\omega_h\in\Vcal_{r,h}^k$, $0\le k\le n$.

We next establish discrete counterparts of the Hodge decomposition
\eqref{eq:hodge-decomposition} and the Poincar\'e inequality
\eqref{eq:gaffney}. To this end, denote the spaces of discrete exact
and closed forms by
\[
  \mathcal B_{r,h}^k:=\rd\Vcal_{r,h}^{k-1},
  \qquad
  \mathcal Z_{r,h}^k:=\ker(\rd|_{\Vcal_{r,h}^k}).
\]
The associated discrete harmonic space is defined by
\begin{equation}
\label{eq:section4-harmonic-spaces}
  \Harm_{r,h}^k
  :=\{\omega_h\in\Vcal_{r,h}^k:\rd \omega_h=0,\quad
    [\omega_h,\rd\tau_h]_{k,h}=0
    \quad\forall\tau_h\in\Vcal_{r,h}^{k-1}\}.
\end{equation}
Since \eqref{eq:complete-virtual-derham-complex} is a complex,
$\mathcal B_{r,h}^k\subseteq\mathcal Z_{r,h}^k$, and
definition \eqref{eq:section4-harmonic-spaces} yields
the $[\cdot,\cdot]_{k,h}$-orthogonal decomposition
\begin{equation}
\label{eq:section4-discrete-hodge}
  \Vcal_{r,h}^k
  =\mathcal B_{r,h}^k\oplus_h\Harm_{r,h}^k
    \oplus_h(\mathcal Z_{r,h}^k)^{\perp_h}.
\end{equation}
Here the last orthogonal complement is taken within $\Vcal_{r,h}^k$.

\begin{theorem}[Discrete Poincar\'e inequality]
\label{thm:complete-orthogonal-poincare}
Under Assumption~\ref{ass:mesh-regularity}, the complete virtual complex
satisfies the following estimates for $0\le k\le n$.
For every $\omega_h\in\Vcal_{r,h}^k$, there exists
$\widetilde\omega_h\in\Vcal_{r,h}^k$ such that
\begin{equation}
\label{eq:section4-poincare}
  \rd \widetilde\omega_h=\rd \omega_h,
  \qquad \norm{\widetilde\omega_h}\lesssim\norm{\rd \omega_h}.
\end{equation}
Moreover, the discrete orthogonal complement satisfies
\begin{equation}
\label{eq:section4-orthogonal-poincare}
  \norm{\omega_h}\le C_{\mathrm P}\norm{\rd \omega_h}
  \qquad\forall \omega_h\in(\mathcal Z_{r,h}^k)^{\perp_h}.
\end{equation}
The constant $C_{\mathrm P}\ge1$ is independent of $h$.
\end{theorem}

\begin{proof}
Estimate \eqref{eq:section4-poincare} follows from the trimmed-family
result of Theorem~\ref{thm:appendix-trimmed-poincare} in
Appendix~\ref{sec:appendix-discrete-poincare} through the stable cochain
retraction of Corollary~\ref{cor:cochain-retraction-stability}.
For $k<n$, take the representative $\widetilde\omega_h^-$ supplied there for
$\omega_h\in\Vcal_{r,h}^k\subset\Vcal_{r+1,h}^{-,k}$ and set
$\widetilde\omega_h:=P_{r,h}^k\widetilde\omega_h^-$. The cochain identities
\eqref{eq:P-cochain-reduction} and stability
\eqref{eq:cochain-retraction-stability} give
\[
  \rd \widetilde\omega_h=P_{r,h}^{k+1}\rd \widetilde\omega_h^-
  =P_{r,h}^{k+1}\rd \omega_h=\rd \omega_h,
  \qquad
  \norm{\widetilde\omega_h}\lesssim\norm{\widetilde\omega_h^-}
  \lesssim\norm{\rd \omega_h}.
\]
To prove \eqref{eq:section4-orthogonal-poincare}, let $k<n$ and
$\omega_h\in(\mathcal Z_{r,h}^k)^{\perp_h}$, and choose $\widetilde\omega_h$ as in
\eqref{eq:section4-poincare}.
Since $\widetilde\omega_h-\omega_h\in\mathcal Z_{r,h}^k$ and
$\omega_h\in(\mathcal Z_{r,h}^k)^{\perp_h}$, we have
$[\omega_h,\widetilde\omega_h-\omega_h]_{k,h}=0$. Hence
\[
\begin{aligned}
  \norm{\omega_h}_h^2
  &=[\omega_h,\widetilde\omega_h]_{k,h}
   \le\norm{\omega_h}_h\norm{\widetilde\omega_h}_h\\
  &\lesssim\norm{\omega_h}_h\norm{\widetilde\omega_h}
   \lesssim\norm{\omega_h}_h\norm{\rd \omega_h}
   \lesssim\norm{\omega_h}_h\norm{\rd \omega_h}_h.
\end{aligned}
\]
The first and last $\lesssim$ bounds follow from
\eqref{Coercivity}, and the middle one from
\eqref{eq:section4-poincare}.
Absorbing the constants above into $C_{\mathrm P}$, we obtain
\begin{equation}
\label{eq:discrete-norm-poincare}
  \norm{\omega_h}_h\le C_{\mathrm P}\norm{\rd \omega_h}_h,
  \qquad \omega_h\in(\mathcal Z_{r,h}^k)^{\perp_h}.
\end{equation}
Applying \eqref{Coercivity} once more gives
\eqref{eq:section4-orthogonal-poincare},
with the resulting constant still denoted by $C_{\mathrm P}$.
For $k=n$, take $\widetilde\omega_h=0$; the second assertion is also immediate since
$(\mathcal Z_{r,h}^n)^{\perp_h}=\{0\}$.
\end{proof}
Alternatively, both estimates
\eqref{eq:section4-poincare} and
\eqref{eq:section4-orthogonal-poincare}
follow from the fully discrete DDR Poincar\'e inequalities
through the cochain isomorphism and norm equivalence;
see Appendix~\ref{subsec:appendix-analytical-transfer}.

\subsection{Discrete mixed problem and inf--sup stability}
\label{subsec:mixed-problems}
We discretize the mixed problem
\eqref{eq:continuous-mixed-1}--\eqref{eq:continuous-mixed-harmonic}
on the product space
\[
  \mathbb V_{r,h}^k
  :=\Vcal_{r,h}^{k-1}\times\Vcal_{r,h}^{k}\times\Harm_{r,h}^k.
\]
Using the discrete inner products \eqref{eq:global-discrete-bilinear},
we define the bilinear form
$\mathcal A_h:\mathbb V_{r,h}^k\times\mathbb V_{r,h}^k\to\mathbb R$ by
\begin{equation}
\label{eq:section4-Ah}
\begin{aligned}
  \mathcal A_h(\boldsymbol u_h,\boldsymbol v_h)
  &:=[\sigma_h,\tau_h]_{k-1,h}-[u_h,\rd\tau_h]_{k,h}
    +[\rd\sigma_h,v_h]_{k,h}
  \\
  &\quad+[\rd u_h,\rd v_h]_{k+1,h}
    +[p_h,v_h]_{k,h}+[u_h,q_h]_{k,h}.
\end{aligned}
\end{equation}
for $\boldsymbol u_h=(\sigma_h,u_h,p_h)$ and
$\boldsymbol v_h=(\tau_h,v_h,q_h)$ in $\mathbb V_{r,h}^k$.
Let $j_h$ be the elementwise $L^2$ projection of $j$, with
$j_h|_T:=\Proj{r}{T}{k}j$. We define the computable load functional by
\begin{equation}
\label{eq:discrete-load}
  \ell_h(v_h)
  :=\sum_{T\in\Th}
  \left\langle j_h,\Proj{r}{T}{k}v_h\right\rangle_T
  =\langle j_h,v_h\rangle.
\end{equation}
The discrete mixed problem then reads: find
$\boldsymbol u_h=(\sigma_h,u_h,p_h)\in\mathbb V_{r,h}^k$ such that
\begin{equation}
\label{eq:discrete-mixed}
  \mathcal A_h(\boldsymbol u_h,\boldsymbol v_h)=\ell_h(v_h)
  \qquad\forall\boldsymbol v_h\in\mathbb V_{r,h}^k.
\end{equation}
The DOF isomorphism in Remark~\ref{rem:fully-discrete-identification}
yields an algebraically equivalent fully discrete formulation,
with the discrete inner products and load identified in
Appendix~\ref{subsec:appendix-analytical-transfer}.
That appendix also details how the norm equivalences transfer
the stability and convergence estimates established below.

We carry out the analysis in the conforming setting, equipping
the continuous and discrete product spaces with the norm
\begin{equation}
\label{eq:discrete-product-norm}
  \norm{(\tau,v,q)}_V^2
  :=\norm{\tau}^2+\norm{\rd\tau}^2
    +\norm{v}^2+\norm{\rd v}^2+\norm{q}^2.
\end{equation}

The norm equivalence \eqref{Coercivity}, the discrete Hodge decomposition
\eqref{eq:section4-discrete-hodge}, and the Poincar\'e inequality
\eqref{eq:section4-orthogonal-poincare} yield the following uniform
inf--sup bound.

\begin{theorem}[Discrete inf--sup stability]
\label{thm:section4-inf-sup}
There exists a constant $\beta>0$, independent of $h$, such that
for $0\le k\le n$ and every $\boldsymbol u_h\in\mathbb V_{r,h}^k$,
\begin{equation}
\label{eq:section4-inf-sup}
  \beta\norm{\boldsymbol u_h}_V
  \le
  \sup_{\boldsymbol v_h\in\mathbb V_{r,h}^k\setminus\{0\}}
  \frac{\mathcal A_h(\boldsymbol u_h,\boldsymbol v_h)}
       {\norm{\boldsymbol v_h}_V}.
\end{equation}
\end{theorem}

\begin{proof}
For $\boldsymbol u_h=(\sigma_h,u_h,p_h)$, decompose
$u_h=u_{\mathrm B,h}+u_{\mathrm H,h}+u_{\perp,h}$ according to
\eqref{eq:section4-discrete-hodge}. For $k\ge1$, let
$\phi_h\in(\mathcal Z_{r,h}^{k-1})^{\perp_h}$ be the unique form
satisfying $\rd\phi_h=u_{\mathrm B,h}$; for $k=0$, set $\phi_h=0$.
Applying \eqref{eq:discrete-norm-poincare}
with a common constant $C_{\mathrm P}$ and using
$\rd u_{\perp,h}=\rd u_h$, we obtain
\begin{equation}
\label{eq:infsup-lifting-bounds}
  \rd\phi_h=u_{\mathrm B,h},\qquad
  \norm{\phi_h}_h\le C_{\mathrm P}\norm{u_{\mathrm B,h}}_h,\qquad
  \norm{u_{\perp,h}}_h\le C_{\mathrm P}\norm{\rd u_h}_h.
\end{equation}
Choose $\boldsymbol v_h=(\tau_h,v_h,q_h)\in\mathbb V_{r,h}^k$ by
\[
  \tau_h:=\sigma_h-C_{\mathrm P}^{-2}\phi_h,\qquad
  v_h:=u_h+\rd\sigma_h+p_h,\qquad
  q_h:=u_{\mathrm H,h}-p_h.
\]
By \eqref{eq:infsup-lifting-bounds}, the orthogonal decomposition
\eqref{eq:section4-discrete-hodge}, and the norm equivalence
\eqref{Coercivity},
\begin{equation}
\label{eq:infsup-test-bound}
  \norm{\boldsymbol v_h}_V\lesssim\norm{\boldsymbol u_h}_V.
\end{equation}
Using symmetry, $\rd^2=0$, $\rd p_h=0$, and the orthogonality in
\eqref{eq:section4-discrete-hodge}, we obtain
\begin{equation}
\label{eq:infsup-bilinear-expansion}
\begin{aligned}
  \mathcal A_h(\boldsymbol u_h,\boldsymbol v_h)
  &=\norm{\sigma_h}_h^2+\norm{\rd\sigma_h}_h^2
    +\norm{\rd u_h}_h^2\\
  &\quad+C_{\mathrm P}^{-2}\norm{u_{\mathrm B,h}}_h^2
    -C_{\mathrm P}^{-2}[\sigma_h,\phi_h]_{k-1,h}
    +\norm{p_h}_h^2+\norm{u_{\mathrm H,h}}_h^2.
\end{aligned}
\end{equation}
The only term of indefinite sign is the cross term. By the
Cauchy--Schwarz inequality, \eqref{eq:infsup-lifting-bounds}, and
Young's inequality,
\[
\begin{aligned}
  C_{\mathrm P}^{-2}|[\sigma_h,\phi_h]_{k-1,h}|
  &\le C_{\mathrm P}^{-2}\norm{\sigma_h}_h\norm{\phi_h}_h
  \le C_{\mathrm P}^{-1}\norm{\sigma_h}_h\norm{u_{\mathrm B,h}}_h
  \\
  &\le\tfrac12\norm{\sigma_h}_h^2
    +\tfrac{1}{2C_{\mathrm P}^2}\norm{u_{\mathrm B,h}}_h^2.
\end{aligned}
\]
Substituting this estimate into
\eqref{eq:infsup-bilinear-expansion} yields
\begin{align*}
  \mathcal A_h(\boldsymbol u_h,\boldsymbol v_h)
  &\ge\tfrac12\norm{\sigma_h}_h^2+\norm{\rd\sigma_h}_h^2
    +\norm{\rd u_h}_h^2
    +\tfrac{1}{2C_{\mathrm P}^2}\norm{u_{\mathrm B,h}}_h^2
    +\norm{p_h}_h^2+\norm{u_{\mathrm H,h}}_h^2
  \\
  &\gtrsim\norm{\boldsymbol u_h}_V^2.
\end{align*}
The final bound follows from \eqref{Coercivity},
\eqref{eq:section4-discrete-hodge}, and \eqref{eq:infsup-lifting-bounds}.
Combining this bound with \eqref{eq:infsup-test-bound} proves
\eqref{eq:section4-inf-sup}.
\end{proof}

By \eqref{Stability} and \eqref{eq:discrete-load},
\[
  |\mathcal A_h(\boldsymbol u_h,\boldsymbol v_h)|
  \lesssim\norm{\boldsymbol u_h}_V\norm{\boldsymbol v_h}_V,
  \qquad |\ell_h(v_h)|\le\norm{j_h}\norm{v_h}
  \le\norm{j}\norm{v_h}.
\]
By the inf--sup bound \eqref{eq:section4-inf-sup}, the discrete problem
\eqref{eq:discrete-mixed} has a unique solution satisfying
\[
  \norm{\boldsymbol u_h}_V\le\beta^{-1}\norm{j_h}
  \le\beta^{-1}\norm{j}.
\]

\subsection{Convergence analysis}
\label{subsec:convergence-analysis}
As noted in Remark~\ref{rem:low-regularity-quasi-interpolation}, we
use canonical interpolation and assume sufficient regularity of the
exact solution and the relevant harmonic forms.

\subsubsection{Approximation of harmonic spaces}
\label{subsubsec:harmonic-approximation}
Note that the discrete harmonic space $\Harm_{r,h}^k$
need not be contained in the continuous harmonic space $\Harm^k$,
even in conforming finite element methods
\cite[Sections~3.3--3.4]{ArnoldFalkWinther2010}.
Following the FEEC analysis of harmonic gaps, we quantify
this discrepancy to account for the discrete harmonic constraint.
For closed subspaces $E,F\subset L^2\Lambda^k(\Omega)$, define
\begin{equation}
\label{eq:harmonic-gap-definition}
  \delta(E,F):=\sup_{\substack{\omega\in E\\\norm{\omega}\le1}}
       \inf_{\zeta\in F}\norm{\omega-\zeta},
  \qquad
  \operatorname{gap}(E,F):=\max\{\delta(E,F),\delta(F,E)\}.
\end{equation}
By Theorem~\ref{thm:veec-cohomology} and the Hodge decompositions
\eqref{eq:hodge-decomposition} and \eqref{eq:section4-discrete-hodge},
$\dim\Harm^k=\dim\Harm_{r,h}^k<\infty$.
Since the two harmonic spaces have the same finite dimension,
\cite[Lemma~2]{DemlowHirani2014} implies
\begin{equation}
\label{eq:harmonic-gap-symmetry}
  \operatorname{gap}(\Harm^k,\Harm_{r,h}^k)
  =\delta(\Harm^k,\Harm_{r,h}^k)
  =\delta(\Harm_{r,h}^k,\Harm^k).
\end{equation}
\par
For $0\le\ell\le n$ and $\omega\in L^2\Lambda^\ell(\Omega)$, define its polynomial
approximation error by
\begin{equation}
\label{eq:polynomial-approximation-error}
  E_h^\ell(\omega):=\left(\sum_{T\in\Th}
    \norm{\omega-\Proj{r}{T}{\ell}\omega}_T^2\right)^{1/2}.
\end{equation}
For $\omega\in L^2\Lambda^\ell(\Omega)$ and
$\omega_h,\zeta_h\in\Vcal_{r,h}^\ell$, define the consistency error by
\begin{equation}
\label{eq:section5-consistency-splitting}
\begin{aligned}
  \mathfrak c_h^\ell(\omega,\omega_h;\zeta_h)
  &:=[\omega_h,\zeta_h]_{\ell,h}-\langle\omega,\zeta_h\rangle\\
  &=\sum_{T\in\Th}\bigl([\omega_h-\Proj{r}{T}{\ell}\omega,\zeta_h]_{\ell,T}
       +\langle\Proj{r}{T}{\ell}\omega-\omega,\zeta_h\rangle_T\bigr),
\end{aligned}
\end{equation}
where the second equality follows from polynomial consistency
\eqref{Consistency}. Applying \eqref{Stability}, the triangle inequality,
and the Cauchy--Schwarz inequality over the mesh gives
\begin{equation}
\label{eq:section5-consistency-estimate}
  |\mathfrak c_h^\ell(\omega,\omega_h;\zeta_h)|
  \le\bigl(\alpha^*\norm{\omega-\omega_h}
       +(\alpha^*+1)E_h^\ell(\omega)\bigr)\norm{\zeta_h}.
\end{equation}

\begin{theorem}[Approximation of discrete harmonic spaces]
\label{thm:canonical-harmonic-approximation}
There is a constant $C>0$, independent of $h$, such that
\begin{equation}
\label{eq:harmonic-gap-best-approximation}
  \operatorname{gap}(\Harm^k,\Harm_{r,h}^k)
  \le C\sup_{\substack{q\in\Harm^k\\\norm{q}\le1}}
  \left(\inf_{\zeta_h\in\mathcal Z_{r,h}^k}\norm{q-\zeta_h}+E_h^k(q)\right).
\end{equation}
If $\Harm^k\subset H^{s_{\Harm}}\Lambda^k(\Omega)$ for some
$s_{\Harm}>(n-k)/2$, then
\begin{equation}
\label{eq:canonical-harmonic-gap-rate}
  \operatorname{gap}(\Harm^k,\Harm_{r,h}^k)
  \le C C_{\Harm}h^{\min\{s_{\Harm},r+1\}}.
\end{equation}
Here $C_{\Harm}$ depends only on $\Omega$, $k$, and $s_{\Harm}$.
For $k=0$, $\operatorname{gap}(\Harm^0,\Harm_{r,h}^0)=0$,
without any additional regularity assumption.
\end{theorem}

\begin{proof}
Let $q\in\Harm^k$ and $\zeta_h\in\mathcal Z_{r,h}^k$.
The discrete Hodge decomposition
\eqref{eq:section4-discrete-hodge} gives
$\zeta_h=q_h+b_h$, where $q_h\in\Harm_{r,h}^k$ and
$b_h\in\mathcal B_{r,h}^k$.
Since $\mathcal B_{r,h}^k\subset\mathcal B^k$,
continuous and discrete orthogonality give
$\langle q,b_h\rangle=[q_h,b_h]_{k,h}=0$, and hence
\[
  [b_h,b_h]_{k,h}
  =[\zeta_h,b_h]_{k,h}
  =\mathfrak c_h^k(q,\zeta_h;b_h).
\]
Combining \eqref{Coercivity} with
\eqref{eq:section5-consistency-estimate} gives
\[
  \norm{b_h}
  \le \frac{\alpha^*}{\alpha_*}\norm{q-\zeta_h}
     +\frac{\alpha^*+1}{\alpha_*}E_h^k(q).
\]
Since $q-q_h=q-\zeta_h+b_h$, the triangle inequality yields
\begin{equation*}
  \inf_{q_h'\in\Harm_{r,h}^k}\norm{q-q_h'}
  \le\norm{q-q_h}
  \lesssim\norm{q-\zeta_h}+E_h^k(q).
\end{equation*}
Taking the infimum over $\zeta_h\in\mathcal Z_{r,h}^k$
and the supremum over $\norm{q}\le1$ in $\Harm^k$,
then using \eqref{eq:harmonic-gap-symmetry}, proves
\eqref{eq:harmonic-gap-best-approximation}.

Since $\Harm^k$ is finite dimensional, the assumed Sobolev regularity implies
$\norm{q}_{H^{s_{\Harm}}(\Omega)}
\le C_{\Harm}\norm{q}$ for all $q\in\Harm^k$,
with $C_{\Harm}$ independent of $h$.
Moreover, $\rd q=0$ and \eqref{eq:commuting-regularity-domain}
allow us to apply the commuting identity
\eqref{eq:virtual-interpolator-commutativity}, giving
$\Ih{k}q\in\mathcal Z_{r,h}^k$
(automatically when $k=n$).
The interpolation estimate
\eqref{eq:section4-local-interpolation} and polynomial approximation
\eqref{eq:section4-poly-approx} give
\begin{equation}
\label{eq:harmonic-interpolation-approximation}
  \norm{q-\Ih{k}q}+E_h^k(q)
  \lesssim h^{\min\{s_{\Harm},r+1\}}\norm{q}_{H^{s_{\Harm}}(\Omega)}
  \le C_{\Harm}h^{\min\{s_{\Harm},r+1\}}\norm{q}.
\end{equation}
Choosing $\zeta_h=\Ih{k}q$ in
\eqref{eq:harmonic-gap-best-approximation} and taking the
supremum over $q\in\Harm^k$ with $\norm{q}\le1$ proves
\eqref{eq:canonical-harmonic-gap-rate}.
In degree zero, both harmonic spaces consist of functions constant
on each connected component of $\Omega$, so the polynomial projection
errors vanish.
\end{proof}

In three dimensions, the commuting quasi-interpolator of
Remark~\ref{rem:low-regularity-quasi-interpolation} can replace
the canonical interpolator in the harmonic approximation step
\eqref{eq:harmonic-interpolation-approximation} and in the mixed
error analysis below. Both arguments therefore extend to
$H\Lambda^\ell(\rd;\Omega)$ in each relevant form degree $\ell$,
without the additional Sobolev regularity required by
canonical interpolation.

\subsubsection{Error estimates for the mixed problem}
\label{subsubsec:mixed-error-estimates}
Let $\boldsymbol u=(\sigma,u,p)$ solve
\eqref{eq:continuous-mixed-1}--\eqref{eq:continuous-mixed-harmonic}.
Using the orthogonal projection notation
\eqref{eq:l2-orthogonal-projection}, set
$u_{\mathrm B}:=\Pi_{\mathcal B^k}u$.
For a form $\omega$ of degree $\ell$ in the canonical interpolation domain,
write
\[
  e_h^\ell(\omega):=\bigl(\norm{\omega-\Ih{\ell}\omega}^2
       +\norm{\rd(\omega-\Ih{\ell}\omega)}^2\bigr)^{1/2}.
\]
For the error analysis, define
\begin{equation*}
  \boldsymbol u_I:=(\sigma_I,u_I,p_I)
  :=(\Ih{k-1}\sigma,\Ih{k}u,\Pi_{\Harm_{r,h}^k}p),
\end{equation*}
and introduce the error quantities
\begin{equation}
\label{eq:canonical-error-indicators}
\begin{aligned}
  \mathcal E_I&:=e_h^{k-1}(\sigma)+e_h^k(u),\\
  \mathcal E_{\mathrm{pol}}&:=E_h^{k-1}(\sigma)+E_h^k(u)
                      +E_h^k(\rd\sigma)+E_h^{k+1}(\rd u),\\
  \mathcal E_{\Harm}&:=\operatorname{gap}(\Harm^k,\Harm_{r,h}^k)
                    \bigl(\norm{p}+\norm{u_{\mathrm B}}\bigr)+E_h^k(p).
\end{aligned}
\end{equation}
The quantity $\mathcal E_{\Harm}$ accounts for the approximation
and consistency errors associated with $p$ and for the discrete
harmonic constraint.

\begin{theorem}[Convergence of the mixed method]
\label{thm:canonical-convergence}
Let $\boldsymbol u_h$ solve \eqref{eq:discrete-mixed}.
Suppose that $\sigma,u,\rd\sigma,\rd u$ have $H^s$
regularity in their respective form degrees, where
$s>n/2$ for $k=0$ and $s>(n-k+1)/2$ for $k\ge1$.
Then
\begin{equation}
\label{eq:canonical-error-bound}
  \norm{\boldsymbol u-\boldsymbol u_h}_V
  \le C\bigl(\mathcal E_I+\mathcal E_{\mathrm{pol}}
             +\mathcal E_{\Harm}+\norm{j-j_h}\bigr),
\end{equation}
where $C$ is independent of $h$.
If, in addition,
$j\in H^{\min\{s,r+1\}}\Lambda^k(\Omega)$ and
$\Harm^k\subset H^{s_{\Harm}}\Lambda^k(\Omega)$
for some $s_{\Harm}>(n-k)/2$, then
\begin{equation}
\label{eq:canonical-convergence-rate}
  \norm{\boldsymbol u-\boldsymbol u_h}_V
  \le C\Bigl(
    h^{\min\{s,r+1\}} C_{\sigma,u,j}
    +C_{\Harm}h^{\min\{s_{\Harm},r+1\}}
      \bigl(\norm{p}+\norm{u_{\mathrm B}}\bigr)
  \Bigr),
\end{equation}
where
\[
  C_{\sigma,u,j}:=\norm{\sigma}_{H^s(\Omega)}
    +\norm{u}_{H^s(\Omega)}+\norm{\rd\sigma}_{H^s(\Omega)}
    +\norm{\rd u}_{H^s(\Omega)}+\norm{j}_{H^{\min\{s,r+1\}}(\Omega)}.
\]
Here $C_{\Harm}$ is the constant in
Theorem~\ref{thm:canonical-harmonic-approximation},
and $s_{\Harm}$ need not equal $s$.
The constants are independent of $h$.
When $\Harm^k=\{0\}$ or $k=0$, no harmonic regularity
assumption is needed, and the harmonic term in
\eqref{eq:canonical-convergence-rate} can be omitted.
\end{theorem}

\begin{proof}
The assumed regularity ensures that $u$ and $\sigma$
belong to their canonical interpolation domains
\eqref{eq:commuting-regularity-domain}.
For $\boldsymbol v_h\in\mathbb V_{r,h}^k$, define the residual
\begin{equation*}
  \mathcal E_h(\boldsymbol v_h)
  :=\ell_h(v_h)-\mathcal A_h(\boldsymbol u_I,\boldsymbol v_h).
\end{equation*}
By \eqref{eq:discrete-mixed},
$\mathcal A_h(\boldsymbol u_h-\boldsymbol u_I,\boldsymbol v_h)
=\mathcal E_h(\boldsymbol v_h)$.
The inf--sup bound \eqref{eq:section4-inf-sup} therefore gives
\begin{equation}
\label{eq:canonical-residual-infsup}
  \norm{\boldsymbol u_h-\boldsymbol u_I}_V
  \le\beta^{-1}
  \sup_{\boldsymbol v_h\in\mathbb V_{r,h}^k\setminus\{0\}}
  \frac{|\mathcal E_h(\boldsymbol v_h)|}{\norm{\boldsymbol v_h}_V}.
\end{equation}
We now estimate the residual.

Let $\mathcal A$ denote the bilinear form obtained from
\eqref{eq:section4-Ah} by replacing each discrete inner product
with the exact $L^2$ inner product, and set
$\ell(v):=\langle j,v\rangle$.
For $\boldsymbol v_h=(\tau_h,v_h,q_h)\in\mathbb V_{r,h}^k$,
testing \eqref{eq:continuous-mixed-1}--\eqref{eq:continuous-mixed-2}
with $\tau_h$ and $v_h$ and adding gives
\begin{equation*}
  \mathcal A(\boldsymbol u,\boldsymbol v_h)
  =\ell(v_h)+\langle u,q_h\rangle.
\end{equation*}
The additional term $\langle u,q_h\rangle$ is controlled below
by the harmonic gap estimated in
Theorem~\ref{thm:canonical-harmonic-approximation}.
Expressing each difference of corresponding inner products in
$\mathcal A(\boldsymbol u,\boldsymbol v_h)
-\mathcal A_h(\boldsymbol u_I,\boldsymbol v_h)$
in terms of the consistency errors
\eqref{eq:section5-consistency-splitting} gives the residual expansion
\begin{equation}
\label{eq:canonical-residual-expansion}
\begin{aligned}
  \mathcal E_h(\boldsymbol v_h)
  &=\langle j_h-j,v_h\rangle
     +\mathcal A(\boldsymbol u,\boldsymbol v_h)
     -\mathcal A_h(\boldsymbol u_I,\boldsymbol v_h)
     -\langle u,q_h\rangle\\
  &=\langle j_h-j,v_h\rangle
     -\mathfrak c_h^{k-1}(\sigma,\sigma_I;\tau_h)
     +\mathfrak c_h^k(u,u_I;\rd\tau_h)\\
  &\quad-\mathfrak c_h^k(\rd\sigma,\rd\sigma_I;v_h)
     -\mathfrak c_h^{k+1}(\rd u,\rd u_I;\rd v_h)\\
  &\quad-\mathfrak c_h^k(p,p_I;v_h)
     -\mathfrak c_h^k(u,u_I;q_h)-\langle u,q_h\rangle.
\end{aligned}
\end{equation}
To bound the harmonic contributions, we first use
\eqref{eq:harmonic-gap-symmetry} to obtain
\begin{equation}
\label{eq:harmonic-variable-comparison}
  \norm{p-p_I}
  \le\operatorname{gap}(\Harm^k,\Harm_{r,h}^k)\norm{p}.
\end{equation}
For $q_h\in\Harm_{r,h}^k\subset\mathcal Z^k$,
the continuous Hodge decomposition
\eqref{eq:hodge-decomposition} and $u\perp\Harm^k$ give
\begin{equation}
\label{eq:discrete-harmonic-constraint-error}
\begin{aligned}
  |\langle u,q_h\rangle|
  &=|\langle u_{\mathrm B},q_h-\Pi_{\Harm^k}q_h\rangle|\\
  &\le\norm{u_{\mathrm B}}\norm{q_h-\Pi_{\Harm^k}q_h}\\
  &\le\operatorname{gap}(\Harm^k,\Harm_{r,h}^k)
       \norm{u_{\mathrm B}}\norm{q_h},
\end{aligned}
\end{equation}
where the last inequality follows from
$\delta(\Harm_{r,h}^k,\Harm^k)
=\operatorname{gap}(\Harm^k,\Harm_{r,h}^k)$
in \eqref{eq:harmonic-gap-symmetry}.
Applying \eqref{eq:section5-consistency-estimate} to each consistency
error in \eqref{eq:canonical-residual-expansion} and using
\eqref{eq:harmonic-variable-comparison} and
\eqref{eq:discrete-harmonic-constraint-error}, we obtain
\begin{equation}
\label{eq:canonical-residual-bound}
  |\mathcal E_h(\boldsymbol v_h)|
  \lesssim\bigl(\mathcal E_I+\mathcal E_{\mathrm{pol}}+\mathcal E_{\Harm}
                    +\norm{j-j_h}\bigr)\norm{\boldsymbol v_h}_V.
\end{equation}
Combining \eqref{eq:canonical-residual-infsup} with
\eqref{eq:canonical-residual-bound} yields
\begin{equation*}
  \norm{\boldsymbol u_h-\boldsymbol u_I}_V
  \lesssim\mathcal E_I+\mathcal E_{\mathrm{pol}}
      +\mathcal E_{\Harm}+\norm{j-j_h}.
\end{equation*}
The triangle inequality and
$\norm{\boldsymbol u-\boldsymbol u_I}_V
\le\mathcal E_I+\operatorname{gap}(\Harm^k,\Harm_{r,h}^k)\norm{p}$
then prove
\eqref{eq:canonical-error-bound}.

For the convergence rate, summing the squared local interpolation
estimates \eqref{eq:section4-local-interpolation} and
\eqref{eq:section4-derivative-interpolation} over the mesh and
using polynomial approximation \eqref{eq:section4-poly-approx} gives
\[
  \mathcal E_I+\mathcal E_{\mathrm{pol}}+\norm{j-j_h}
  \lesssim h^{\min\{s,r+1\}} C_{\sigma,u,j}.
\]
The harmonic-gap estimate \eqref{eq:canonical-harmonic-gap-rate},
together with polynomial approximation of $p\in\Harm^k$, gives independently
\[
  \mathcal E_{\Harm}
  \lesssim C_{\Harm}h^{\min\{s_{\Harm},r+1\}}
    \bigl(\norm{p}+\norm{u_{\mathrm B}}\bigr).
\]
Substituting these bounds into \eqref{eq:canonical-error-bound} proves
\eqref{eq:canonical-convergence-rate}.
\end{proof}

\appendix
\section{Scaled trace and Poincar\'e estimates}
\label{sec:appendix-local}
This appendix proves the scaled trace estimates and the Poincar\'e
inequality stated in Subsection~\ref{subsec:mesh-auxiliary-estimates},
whose assumptions and constant conventions apply throughout.
We first recall two identities used in the proofs.
For a vector field $Y$ on $\operatorname{aff}(f)$, let $Y^\flat$
denote the one-form defined by $Y^\flat(v)=Y\cdot v$.
For $\omega\in\Lambda^k(f)$ and $\zeta\in\Lambda^\ell(f)$,
contraction satisfies the graded Leibniz rule
\begin{equation}
\label{eq:contraction-wedge-rule}
  (\omega\wedge\zeta)\lrcorner Y
  =(\omega\lrcorner Y)\wedge\zeta
   +(-1)^k\omega\wedge(\zeta\lrcorner Y).
\end{equation}
Let $\boldsymbol n$ denote the outward unit normal within $\operatorname{aff}(f)$.
For $\zeta\in\Lambda^\ell(f)$ smooth up to the boundary, the normal trace
satisfies \cite[Section~2.2]{Awanou2023}
\begin{equation}
\label{eq:normal-trace-contraction}
  \ntr\zeta=\trc(\zeta\lrcorner\boldsymbol n),
  \qquad 1\le \ell\le d.
\end{equation}

\subsection{Proof of the scaled Rellich-type trace estimates}
\label{app:rellich-proof}
\begin{proof}[Proof of Lemma~\ref{lem:scaled-traces}\textup{(\ref{item:trace-rellich})}]
For $k=0$ or $k=d$, the graph-space assumptions imply
$\omega\in H^1\Lambda^k(f)$, so both traces belong to $L^2(\partial f)$.
For $1\le k\le d-1$, the assumed $L^2$ tangential trace implies
$\ntr\omega\in L^2(\partial f)$ by
\cite[Theorem~6.3]{MitreaMitrea2003} with $p=2$.

Let $c_f$ be the center of $B_f$, and set
\[
  X(x):=x-c_f\quad(x\in f),
  \qquad a:=X\cdot\boldsymbol n=\ntr(X^\flat)\quad\text{on }\partial f.
\]
For $0\le k\le d$, we claim the following radial Rellich-type identity:
\begin{equation}
\label{eq:radial-rellich-identity}
\begin{aligned}
  &\frac12\int_{\partial f}
    a\bigl(|\trc\omega|^2-|\ntr\omega|^2\bigr)
  -\ip{\trc\omega}{\trc(X^\flat)\wedge\ntr\omega}_{\partial f}\\
  &\qquad=
  \left(\frac d2-k\right)\norm{\omega}_f^2
  +\ip{\rd\omega}{X^\flat\wedge\omega}_f
  +\ip{\omega}{X^\flat\wedge\delta\omega}_f.
\end{aligned}
\end{equation}
For a related Rellich-type identity, see
\cite[Eq.~(6.2)]{Mitrea2010}.
By \cite[Proposition~2.5]{MitreaTaylorVasy2005}, $C^2$ forms are dense
in the norm
\[
  \norm{\omega}_f+\norm{\rd\omega}_f
  +\norm{\delta\omega}_f
  +\norm{\trc\omega}_{\partial f}
  +\norm{\ntr\omega}_{\partial f}.
\]
It therefore suffices to prove \eqref{eq:radial-rellich-identity}
for $C^2$ forms.

For the radial field $X$, $\nabla X=\mathrm{Id}$ and
$\delta(X^\flat)=-\operatorname{div}X=-d$.
Applying the generalized Goldberg formula
\cite[Theorem~3.2, Eq.~(3.11)]{DeNicolaYudin2016}
with $X^\flat$ as the fixed one-form gives
\begin{equation}
\label{eq:rellich-codifferential-product}
\begin{aligned}
  \delta(X^\flat\wedge\omega)
  &=\delta(X^\flat)\,\omega-X^\flat\wedge\delta\omega
    -\nabla_X\omega+i_{\mathrm{Id}}\omega\\
  &=-X^\flat\wedge\delta\omega
    -(X\cdot\nabla)\omega-(d-k)\omega.
\end{aligned}
\end{equation}
Here $i_{\mathrm{Id}}\omega=k\omega$, since
$(i_{\mathrm{Id}}\omega)(v_1,\ldots,v_k)
=\sum_{i=1}^{k}\omega(v_1,\ldots,\mathrm{Id}v_i,\ldots,v_k)
=k\omega(v_1,\ldots,v_k)$.
In orthonormal coordinates, $\nabla_X\omega=(X\cdot\nabla)\omega$
denotes coefficientwise differentiation along $X$.

On $\partial f$, the normal-trace identity
\eqref{eq:normal-trace-contraction} and the graded Leibniz rule
\eqref{eq:contraction-wedge-rule} give
\begin{equation}
\label{eq:rellich-normal-product}
\begin{aligned}
  \ntr(X^\flat\wedge\omega)
  =\trc\bigl((X^\flat\wedge\omega)\lrcorner\boldsymbol n\bigr)
  =a\trc\omega-\trc(X^\flat)\wedge\ntr\omega.
\end{aligned}
\end{equation}
We apply the integration-by-parts formula \eqref{eq:green-formula}
with test form $X^\flat\wedge\omega$.
Substituting \eqref{eq:rellich-codifferential-product}--\eqref{eq:rellich-normal-product}
gives
\begin{equation}
\label{eq:rellich-green-expansion}
\begin{aligned}
  \ip{\rd\omega}{X^\flat\wedge\omega}_f
  &=\ip{\omega}{\delta(X^\flat\wedge\omega)}_f
    +\ip{\trc\omega}{\ntr(X^\flat\wedge\omega)}_{\partial f}\\
  &=-\ip{\omega}{X^\flat\wedge\delta\omega}_f
    -\ip{(X\cdot\nabla)\omega}{\omega}_f
    -(d-k)\norm{\omega}_f^2\\
  &\quad+\int_{\partial f}a|\trc\omega|^2
    -\ip{\trc\omega}{\trc(X^\flat)\wedge\ntr\omega}_{\partial f}.
\end{aligned}
\end{equation}
The orthogonal decomposition of $\omega$ into tangential and
normal components gives
\begin{equation}
\label{eq:boundary-tangential-normal-splitting}
  |\omega|^2=|\trc\omega|^2+|\ntr\omega|^2
  \qquad\text{a.e. on }\partial f.
\end{equation}
Differentiating coefficientwise, applying the divergence theorem
with $\operatorname{div}X=d$ and $a=X\cdot\boldsymbol n$, and then using
\eqref{eq:boundary-tangential-normal-splitting}, we obtain
\begin{equation}
\label{eq:rellich-radial-energy}
\begin{aligned}
  2\ip{(X\cdot\nabla)\omega}{\omega}_f
  &=\int_f X\cdot\nabla|\omega|^2
  =\int_{\partial f}a|\omega|^2-d\norm{\omega}_f^2\\
  &=\int_{\partial f}
    a\bigl(|\trc\omega|^2+|\ntr\omega|^2\bigr)
    -d\norm{\omega}_f^2.
\end{aligned}
\end{equation}
Substituting \eqref{eq:rellich-radial-energy} into
\eqref{eq:rellich-green-expansion} and rearranging proves
\eqref{eq:radial-rellich-identity}.

Condition {\rm (M1)} gives
\begin{equation}
\label{eq:rellich-multiplier-bounds}
  |X|\le h_f\quad\text{in }f,
  \qquad
  a\ge\rho h_f\quad\text{a.e. on }\partial f.
\end{equation}
Using \eqref{eq:rellich-multiplier-bounds} in
\eqref{eq:radial-rellich-identity} and applying Cauchy--Schwarz gives
\begin{equation*}
\begin{aligned}
  \rho h_f\norm{\ntr\omega}_{\partial f}^{2}
  &\le h_f\norm{\trc\omega}_{\partial f}^{2}
    +2h_f\norm{\trc\omega}_{\partial f}
           \norm{\ntr\omega}_{\partial f}\\
  &\quad+|d-2k|\norm{\omega}_f^{2}
    +2h_f\norm{\omega}_f
       \bigl(\norm{\rd\omega}_f+\norm{\delta\omega}_f\bigr).
\end{aligned}
\end{equation*}
Young's inequality yields
\begin{equation*}
  2h_f\norm{\trc\omega}_{\partial f}
       \norm{\ntr\omega}_{\partial f}
  \le \frac{\rho h_f}{2}\norm{\ntr\omega}_{\partial f}^{2}
    +\frac{2h_f}{\rho}\norm{\trc\omega}_{\partial f}^{2}.
\end{equation*}
Absorbing the normal-trace term and dividing by $h_f$ gives
\eqref{eq:cell-div-curl-trace}.
Since $|d-2k|\le d$, the hidden constant depends only on $\rho$
for fixed dimension.
Moreover, Young's inequality gives
\begin{equation*}
\begin{aligned}
  \norm{\omega}_f
  \bigl(\norm{\rd\omega}_f+\norm{\delta\omega}_f\bigr)
  \le \frac{1}{2h_f}\norm{\omega}_f^2
  +\frac{h_f}{2}
    \bigl(\norm{\rd\omega}_f+\norm{\delta\omega}_f\bigr)^2.
\end{aligned}
\end{equation*}
Substituting this bound into \eqref{eq:cell-div-curl-trace} and taking
square roots yields \eqref{eq:cell-div-curl-trace-additive}.
Applying both estimates to $\star_f\omega$ gives the bounds with
the two traces interchanged.
\end{proof}

\subsection{Proof of the scaled negative-norm trace estimates}
\label{app:facet-negative-trace}

The scaled trace and dual norms are defined by
\begin{equation}
\label{eq:scaled-boundary-dual-norm}
\begin{aligned}
  \norm{\psi}_{1/2,\partial f}
  &:=\inf_{\substack{\zeta\in H^1\Lambda^{k+1}(f)\\\ntr \zeta=\psi}}
    \bigl(h_f^{-1}\norm{\zeta}_f+\seminorm{\zeta}_{1,f}\bigr),\\
  \norm{\chi}_{-1/2,\partial f}
  &:=\sup_{0\ne\psi\in H_\perp^{1/2}\Lambda^k(\partial f)}
    \frac{|\langle\chi,\psi\rangle_{\partial f}|}
      {\norm{\psi}_{1/2,\partial f}}.
\end{aligned}
\end{equation}

\begin{proof}[Proof of Lemma~\ref{lem:scaled-traces}\textup{(\ref{item:trace-weak})}]
For any $\psi\in H_\perp^{1/2}\Lambda^k(\partial f)$ and any lifting
$\zeta\in H^1\Lambda^{k+1}(f)$ satisfying $\ntr \zeta=\psi$,
the integration-by-parts formula \eqref{eq:green-formula} and the
Cauchy--Schwarz inequality give
\[
\begin{aligned}
  |\langle\trc\omega,\psi\rangle_{\partial f}|
  &=|\langle\rd\omega,\zeta\rangle_f
      -\langle\omega,\delta \zeta\rangle_f|\\
  &\lesssim
    \bigl(\norm{\omega}_f+h_f\norm{\rd\omega}_f\bigr)
    \bigl(h_f^{-1}\norm{\zeta}_f+\seminorm{\zeta}_{1,f}\bigr).
\end{aligned}
\]
Taking the infimum over all such $\zeta$ and then the supremum over
nonzero $\psi$ proves the right-hand inequality of
\eqref{eq:weak-trace}.

Under {\rm (M1)}, the trace-lifting argument in
\cite[Sections~2.3 and~2.7]{BrennerSung2018} applies in any dimension
through a uniformly bi-Lipschitz map to a ball, providing a linear right
inverse $\mathcal R_f:H^{1/2}(\partial f)\to H^1(f)$ of the scalar
boundary trace $\gamma_f$.
Combining its scaled bound with the boundary-interpolation argument
in \cite[Section~2.5]{BrennerSung2018} gives, for
$v\in H^1(\partial f)$,
\begin{equation}
\label{eq:uniform-trace-right-inverse}
\begin{gathered}
  \gamma_f\mathcal R_f v=v,\\
\begin{aligned}
  h_f^{-1}\norm{\mathcal R_f v}_f
  +\seminorm{\mathcal R_f v}_{1,f}
  &\lesssim h_f^{-1/2}\norm{v}_{\partial f}
    +\seminorm{v}_{1/2,\partial f}\\
  &\lesssim h_f^{-1/2}\norm{v}_{\partial f}
    +h_f^{1/2}\seminorm{v}_{1,\partial f}.
\end{aligned}
\end{gathered}
\end{equation}
The hidden constants depend only on $\rho$ for fixed dimension.
Fix $e\in\Delta_{d-1}(f)$ and $\varphi\in H_0^1\Lambda^k(e)$.
Let $\boldsymbol n_e$ be the outward unit normal to $e$ within
$\operatorname{aff}(f)$.
Extend $\varphi$ to ambient vectors by orthogonal projection onto
the tangent space of $e$, and define the ambient boundary form
\[
  \Psi_e:=
  \begin{cases}
    \boldsymbol n_e^\flat\wedge\varphi,&\text{on }e,\\
    0,&\text{on }\partial f\setminus e.
  \end{cases}
\]
Since $\varphi$ has zero trace on $\partial e$, the coefficients
of $\Psi_e$ belong to $H^1(\partial f)$.
Applying $\mathcal R_f$ coefficientwise, set
\[
  L_e\varphi:=\mathcal R_f\Psi_e
  \in H^1\Lambda^{k+1}(f),
  \qquad
  \gamma_fL_e\varphi=\Psi_e.
\]
Since $\varphi\lrcorner\boldsymbol n_e=0$, identities
\eqref{eq:contraction-wedge-rule}--\eqref{eq:normal-trace-contraction} give
\begin{equation*}
  \left.\ntr(L_e\varphi)\right|_e
  =(\boldsymbol n_e^\flat\wedge\varphi)\lrcorner\boldsymbol n_e
  =\varphi,
  \qquad
  \left.\ntr(L_e\varphi)\right|_{\partial f\setminus e}=0.
\end{equation*}
Applying \eqref{eq:uniform-trace-right-inverse} coefficientwise
to $\Psi_e$ gives
\begin{equation}
\label{eq:facet-trace-lifting-bound}
  h_f^{-1}\norm{L_e\varphi}_f+\seminorm{L_e\varphi}_{1,f}
  \lesssim h_f^{-1/2}\norm{\varphi}_e
    +h_f^{1/2}\seminorm{\varphi}_{1,e}
  \lesssim h_f^{1/2}\seminorm{\varphi}_{1,e},
\end{equation}
where the last inequality follows from the Poincar\'e inequality on
$e$ and $h_e\le h_f$.

For $0\ne\varphi\in H_0^1\Lambda^k(e)$, the scaled norm definitions
\eqref{eq:scaled-boundary-dual-norm} and the lifting estimate
\eqref{eq:facet-trace-lifting-bound} give
\[
\begin{aligned}
  |\langle\trc_e\omega,\varphi\rangle_e|
  &=|\langle\trc\omega,\ntr L_e\varphi\rangle_{\partial f}|\\
  &\le\norm{\trc\omega}_{-1/2,\partial f}
    \bigl(h_f^{-1}\norm{L_e\varphi}_f
      +\seminorm{L_e\varphi}_{1,f}\bigr)\\
  &\lesssim h_f^{1/2}\norm{\trc\omega}_{-1/2,\partial f}
    \seminorm{\varphi}_{1,e}.
\end{aligned}
\]
Here the first equality defines the facet pairing by weak-trace duality.
Dividing by $\seminorm{\varphi}_{1,e}$ and taking the supremum gives,
in particular, the left-hand inequality of \eqref{eq:weak-trace}.
\end{proof}

\subsection{Proof of the scaled Sobolev trace estimate on facets}
\label{subsec:appendix-scaled-estimates}

\begin{proof}[Proof of Lemma~\ref{lem:scaled-traces}\textup{(\ref{item:trace-facets})}]
Let $v$ be a scalar coefficient of $\omega$ in fixed orthonormal
coordinates on $\operatorname{aff}(f)$.
For $d\ge2$, choose a covering of $\partial f$ by local Lipschitz
graph neighborhoods with uniformly bounded overlap, as in
\cite[Definition~2.6]{TanakaEtAl2015}.
Under {\rm (M1)}, \cite[Proposition~2.4.4]{HenrotPierre2018}
gives a uniform cone condition at scale $h_f$.
Remark~2.4.8 of the same reference then gives neighborhood radii
at least $c(\rho)h_f$ and graph Lipschitz constants bounded by
$C(\rho)$.
These uniform bounds allow us to apply the Stein extension theorem
\cite[Theorem~5.24]{AdamsFournier2003}.
Together with Sobolev interpolation
\cite[Chapter~7]{AdamsFournier2003} and scaling, it provides a
linear extension operator $\mathscr E_f$ satisfying
$(\mathscr E_fv)|_f=v$ and, for each fixed $s\ge0$,
\begin{equation}
\label{eq:scaled-sobolev-extension}
  \norm{\mathscr E_fv}_{L^2(\mathbb R^d)}
  +h_f^s\seminorm{\mathscr E_fv}_{s,\mathbb R^d}
  \lesssim\norm{v}_f+h_f^s\seminorm{v}_{s,f}.
\end{equation}
The hidden constant depends only on $\rho$ for fixed dimension
and Sobolev order.
For $s>1/2$, the classical trace theorem
\cite[Theorem~8.1]{Wloka1987} gives
\[
  \norm{\gamma_0 w}_{H^{s-1/2}(\mathbb R^{d-1})}
  \lesssim\norm{w}_{H^s(\mathbb R^d)},
  \qquad w\in H^s(\mathbb R^d),
\]
where $\gamma_0$ denotes the trace on $\{x_d=0\}$.
After a translation and rotation identifying $\operatorname{aff}(e)$ with this
hyperplane, scaling by $h_f$, applying the estimate to $\mathscr E_fv$,
and restricting the trace to $e$ yield
\[
\begin{aligned}
  \norm{\trc_e v}_e
    +h_f^{s-1/2}\seminorm{\trc_e v}_{s-1/2,e}
  &\lesssim
    h_f^{-1/2}\norm{\mathscr E_fv}_{L^2(\mathbb R^d)}
    +h_f^{s-1/2}\seminorm{\mathscr E_fv}_{s,\mathbb R^d}\\
  &\lesssim
    h_f^{-1/2}\norm{v}_f+h_f^{s-1/2}\seminorm{v}_{s,f},
\end{aligned}
\]
where the second inequality follows from
\eqref{eq:scaled-sobolev-extension}.
For $d=1$, the same bound follows from the scaled endpoint trace
estimate.
Applying these scalar estimates componentwise proves
\eqref{eq:facet-sobolev-trace}.
\end{proof}

\subsection{Proof of the scaled Poincar\'e inequality}
\label{subsec:appendix-cell-poincare}
\begin{proof}[Proof of Lemma~\ref{lem:cell-poincare}]
First assume $\trc\omega=0$.
Let $E_f^0$ denote extension by zero to $\operatorname{aff}(f)$,
and let $c_f$ be the center of $B_f$.
In orthonormal coordinates on $\operatorname{aff}(f)$, fix
$\theta_0\in C_c^\infty(B(0,1/2))$ with unit integral and set
\[
  \theta_f(x):=(\rho h_f)^{-d}
  \theta_0\!\left(\frac{x-c_f}{\rho h_f}\right).
\]
Then $\theta_f$ is supported in $B_f$ and has unit integral.
Let $T_{\ell,f}$ be the Bogovski\u{\i}-type integral operators
constructed with the averaging function $\theta_f$ in
\cite[Definition~3.1]{CostabelMcIntosh2010}.
Since $f$ is star-shaped with respect to $B_f$, which contains
$\operatorname{supp}\theta_f$, Theorem~3.2 and Corollary~3.4 of the
same reference, together with scaling, yield
\[
  \mathcal T_f^\ell:L^2\Lambda^\ell(f)\longrightarrow
  H_0^1\Lambda^{\ell-1}(f),
  \qquad
  \mathcal T_f^\ell \zeta:=(T_{\ell,f}E_f^0\zeta)|_f,
  \quad 1\le\ell\le d,
\]
with
\begin{equation}
\label{eq:uniform-bogovskii-bound}
  h_f^{-1}\norm{\mathcal T_f^\ell \zeta}_f
  +\seminorm{\mathcal T_f^\ell \zeta}_{1,f}
  \lesssim\norm{\zeta}_f.
\end{equation}
For fixed dimension and $\theta_0$, the hidden constant depends
only on $\rho$.

Since $\trc\omega=0$, extension by zero satisfies
$\rd E_f^0\omega=E_f^0\rd\omega$.
The homotopy identity \cite[Eq.~(3.13)]{CostabelMcIntosh2010},
together with the assumption $\int_f\omega=0$ when $k=d$,
therefore yields
\[
  \omega=\rd\mathcal T_f^k\omega
          +\mathcal T_f^{k+1}\rd\omega,
\]
where the first term is omitted for $k=0$ and the second for $k=d$.
Taking the $L^2(f)$ inner product with $\omega$ and applying the
integration-by-parts formula \eqref{eq:green-formula}, we obtain
\[
\begin{aligned}
  \norm{\omega}_f^2
  &=\ip{\delta\omega}{\mathcal T_f^k\omega}_f
    +\ip{\omega}{\mathcal T_f^{k+1}\rd\omega}_f\\
  &\lesssim h_f\norm{\omega}_f
       \bigl(\norm{\delta\omega}_f+\norm{\rd\omega}_f\bigr).
\end{aligned}
\]
The boundary term vanishes because $\mathcal T_f^k\omega$ has zero
trace for $k\ge1$, and the inequality follows from Cauchy--Schwarz
and \eqref{eq:uniform-bogovskii-bound}.
With the same endpoint conventions, this proves
\eqref{eq:cell-poincare} when $\trc\omega=0$.
The case $\ntr\omega=0$ follows by applying the preceding estimate
to $\star_f\omega$.
\end{proof}

\section{Equivalence with fully discrete polytopal methods}
\label{sec:appendix-fully-discrete}
This appendix establishes cochain isomorphisms from our trimmed
and complete virtual complexes to the fully discrete VEM-inspired
and DDR-inspired complexes of~\cite{Bonaldi2025}, respectively.
For the complete family, it further identifies the corresponding
discrete inner products, loads, and mixed formulations, and
details the transfer of stability and convergence estimates
through the norm equivalences.

\subsection{Reconstruction identities and cochain isomorphisms}
\label{subsec:appendix-fully-discrete-realizations}
Let $f\in\Delta_d(\Mh)$ and $0\le k\le d$.
Define the local DOF maps
$\mathcal J_{r+1,f}^{-,k}$ on $\Vcal_{r+1}^{-,k}(f)$
and $\mathcal J_{r,f}^k$ on $\Vcal_r^k(f)$
by applying the component formulas
\eqref{eq:global-fully-discrete-isomorphisms} and
\eqref{eq:global-complete-isomorphism}, respectively,
to all $e\in\Delta_\ell(f)$ with $k\le\ell\le d$.
Their ranges define the local data spaces:
\[
  \underline{\Vcal}_{r+1,f}^{-,k}
  :=\mathcal J_{r+1,f}^{-,k}\bigl(\Vcal_{r+1}^{-,k}(f)\bigr),
  \qquad
  \underline{\Vcal}_{r,f}^{k}
  :=\mathcal J_{r,f}^{k}\bigl(\Vcal_r^k(f)\bigr).
\]
\paragraph{Trimmed family.}
We denote the data on $e$ by $\omega_e$ when $\ell=k$ and
by $(\omega_e,D_{\omega,e})$ when $\ell>k$.
For $f\in\Delta_{k+1}(\Mh)$ and
$\underline\omega_f\in\underline{\Vcal}_{r+1,f}^{-,k}$, define
$\mathsf d_{r+1,f}^{-,k}\underline\omega_f\in\Pcal_r\Lambda^{k+1}(f)$
by \cite[Eq.~(4.4)]{Bonaldi2025}:
\begin{equation}
\label{eq:trimmed-minimal-derivative-reconstruction}
  \int_f\mathsf d_{r+1,f}^{-,k}\underline\omega_f\wedge(\mu_0+\nu)
  =\int_{\partial f}\star_{\partial f}^{-1}\omega_{\partial f}
       \wedge\trc\mu_0
   +\int_f\star_f^{-1}D_{\omega,f}\wedge\nu,
\end{equation}
for $(\mu_0,\nu)\in\Pcal_0\Lambda^0(f)\times\Kcal_r^0(f)$.
Here $\star_{\partial f}^{-1}\omega_{\partial f}$ denotes the
facetwise assembly of $\star_e^{-1}\omega_e$, using the fixed facet
orientations. Boundary integrals here and below use the orientation
induced by $f$, as in \eqref{eq:integral.relative.orientation}.
The scalar decomposition \eqref{eq:polynomial-koszul-decomposition-k0}
makes these moments unisolvent.
The global operator is then \cite[Eq.~(4.6)]{Bonaldi2025}
\begin{equation}
\label{eq:global-trimmed-fully-discrete-derivative}
  \bigl(\underline{\rd}_{r+1,h}^{-,k}\underline\omega_h\bigr)_f
  :=\begin{cases}
    \star_f\mathsf d_{r+1,f}^{-,k}\underline\omega_f,
       &\dim f=k+1,\\
    (D_{\omega,f},0),&\dim f\ge k+2.
  \end{cases}
\end{equation}

\paragraph{Complete family.}
We denote the data component on $e$ by $\omega_e$.
Following \cite[Eqs.~(3.4)--(3.5)]{Bonaldi2025}, define the
potential $\mathsf P_{r,f}^k\underline\omega_f\in\Pcal_r\Lambda^k(f)$
and derivative
$\mathsf d_{r,f}^k\underline\omega_f\in\Pcal_r\Lambda^{k+1}(f)$
recursively for
$\underline\omega_f\in\underline{\Vcal}_{r,f}^k$.
If $d=k$, set $\mathsf P_{r,f}^k\underline\omega_f:=\star_f^{-1}\omega_f$.
If $d>k$, assemble
$\mathsf P_{r,\partial f}^k\underline\omega_{\partial f}$ from the
facet potentials $\mathsf P_{r,e}^k\underline\omega_e$,
$e\in\Delta_{d-1}(f)$, and first define the exterior derivative by
\begin{equation}
\label{eq:ddr-local-derivative-definition}
  \int_f\mathsf d_{r,f}^k\underline\omega_f\wedge\eta
  =(-1)^{k+1}\int_f\star_f^{-1}\omega_f\wedge\rd\eta
   +\int_{\partial f}
       \mathsf P_{r,\partial f}^k\underline\omega_{\partial f}
       \wedge\trc\eta,
  \quad\eta\in\Pcal_r\Lambda^{d-k-1}(f).
\end{equation}
Then define the potential by
\begin{equation}
\label{eq:ddr-local-potential-definition}
\begin{aligned}
  (-1)^{k+1}\int_f
    \mathsf P_{r,f}^k\underline\omega_f\wedge(\rd\mu+\nu)
  &=\int_f\mathsf d_{r,f}^k\underline\omega_f\wedge\mu
    -\int_{\partial f}
       \mathsf P_{r,\partial f}^k\underline\omega_{\partial f}
       \wedge\trc\mu\\
  &\quad+(-1)^{k+1}\int_f\star_f^{-1}\omega_f\wedge\nu,
\end{aligned}
\end{equation}
for $(\mu,\nu)\in\Kcal_{r+1}^{d-k-1}(f)\times\Kcal_r^{d-k}(f)$.
The isomorphism \eqref{eq:polynomial-koszul-isomorphism} maps
$(\mu,\nu)$ to $\rd\mu+\nu\in\Pcal_r\Lambda^{d-k}(f)$, so the
nondegenerate wedge pairing determines a unique potential.
Following \cite[Eq.~(3.8)]{Bonaldi2025}, assemble the local derivatives as
\begin{equation}
\label{eq:global-complete-fully-discrete-derivative}
  \bigl(\underline{\rd}_{r,h}^{k}\underline\omega_h\bigr)_f
  :=\Pi_{r,f}^{-,d-k-1}
       \bigl(\star_f\mathsf d_{r,f}^{k}\underline\omega_f\bigr),
  \qquad f\in\Delta_d(\Mh),\quad d\ge k+1.
\end{equation}
\begin{lemma}[Local reconstruction identities]
\label{lem:local-reconstruction-identities}
Let $f\in\Delta_d(\Mh)$ and $0\le k\le d$.
If $d=k+1$, then every $\omega\in\Vcal_{r+1}^{-,k}(f)$ satisfies
\begin{equation}
\label{eq:local-trimmed-reconstruction-identification}
  \mathsf d_{r+1,f}^{-,k}\mathcal J_{r+1,f}^{-,k}\omega
  =\rd\omega.
\end{equation}
For every $\omega\in\Vcal_r^k(f)$,
\begin{equation}
\label{eq:local-vem-ddr-reconstruction-identification}
  \mathsf d_{r,f}^k\mathcal J_{r,f}^k\omega
  =\Proj{r}{f}{k+1}\rd\omega\quad(k<d),
  \qquad
  \mathsf P_{r,f}^k\mathcal J_{r,f}^k\omega
  =\Proj{r}{f}{k}\omega.
\end{equation}
\end{lemma}
\begin{proof}
For the trimmed family, set
$\underline\omega_f:=\mathcal J_{r+1,f}^{-,k}\omega$.
By \textup{(T3)} and \eqref{eq:trimmed-top},
$\rd\omega\in\Pcal_r\Lambda^{k+1}(f)$.
The DOF components in~\eqref{eq:global-fully-discrete-isomorphisms} give
$\star_e^{-1}\omega_e=\trc_e\omega$ on each facet and
$D_{\omega,f}=\Pi_{r,f}^{\kappa,0}(\star_f\rd\omega)$.
Using these identities, the projection identity \eqref{eq:hodge-projection-moments}, and the
integration-by-parts formula \eqref{eq:stokes-integration-by-parts},
we rewrite \eqref{eq:trimmed-minimal-derivative-reconstruction} as
\[
\begin{aligned}
  \int_f\mathsf d_{r+1,f}^{-,k}\underline\omega_f\wedge(\mu_0+\nu)
  &=\int_{\partial f}\trc\omega\wedge\trc\mu_0
    +\int_f\rd\omega\wedge\nu
   =\int_f\rd\omega\wedge(\mu_0+\nu).
\end{aligned}
\]
Since the tests span $\Pcal_r\Lambda^0(f)$, this proves
\eqref{eq:local-trimmed-reconstruction-identification}.

For the complete family, set
$\underline\omega_f:=\mathcal J_{r,f}^k\omega$ and argue by induction
on $d$, with $k$ fixed.
If $d=k$, \eqref{eq:trimmed-top} gives
$\omega\in\Pcal_r\Lambda^k(f)$, and hence
$\mathsf P_{r,f}^k\underline\omega_f=\star_f^{-1}\omega_f
=\omega=\Proj{r}{f}{k}\omega$.
If $d>k$, first let $\eta\in\Pcal_r\Lambda^{d-k-1}(f)$.
The induction hypothesis gives
$\mathsf P_{r,\partial f}^k\underline\omega_{\partial f}
=\Pi_{r,\partial f}^k\trc\omega$.
Since $\omega_f=\Pi_{r,f}^{-,d-k}(\star_f\omega)$,
the projection identity \eqref{eq:hodge-projection-moments},
the derivative reconstruction defined in
\eqref{eq:ddr-local-derivative-definition}, and the moment identity
for $\rd\omega$ in \eqref{eq:projection-computability-derivative} give
\[
\begin{aligned}
  \int_f\mathsf d_{r,f}^k\underline\omega_f\wedge\eta
  &=(-1)^{k+1}\int_f
       \star_f^{-1}\Pi_{r,f}^{-,d-k}(\star_f\omega)\wedge\rd\eta
    +\int_{\partial f}
       \Pi_{r,\partial f}^k\trc\omega\wedge\trc\eta\\
  &=\int_f\Proj{r}{f}{k+1}\rd\omega\wedge\eta.
\end{aligned}
\]
Hence $\mathsf d_{r,f}^k\underline\omega_f
=\Proj{r}{f}{k+1}\rd\omega$.
Next, for $(\mu,\nu)\in\Kcal_{r+1}^{d-k-1}(f)\times\Kcal_r^{d-k}(f)$,
we substitute this derivative reconstruction identity and the induction
hypothesis for the facet potentials into
\eqref{eq:ddr-local-potential-definition}, then apply
the projection identity \eqref{eq:hodge-projection-moments} to obtain
\[
\begin{aligned}
  (-1)^{k+1}\int_f
    \mathsf P_{r,f}^k\underline\omega_f\wedge(\rd\mu+\nu)
  &=\int_f\Proj{r}{f}{k+1}\rd\omega\wedge\mu
    -\int_{\partial f}
       \Pi_{r,\partial f}^k\trc\omega\wedge\trc\mu\\
  &\quad+(-1)^{k+1}\int_f
       \star_f^{-1}\Pi_{r,f}^{-,d-k}(\star_f\omega)\wedge\nu\\
  &=(-1)^{k+1}\int_f
       \Proj{r}{f}{k}\omega\wedge(\rd\mu+\nu).
\end{aligned}
\]
The last equality uses the moment identity
\eqref{eq:projection-computability-stokes}, and the polynomial
isomorphism \eqref{eq:polynomial-koszul-isomorphism} then gives
$\mathsf P_{r,f}^k\underline\omega_f=\Proj{r}{f}{k}\omega$.
\end{proof}

\begin{remark}[Improved scalar projection]
\label{rem:degree-r-plus-one-scalar-potential}
For $k=0$, the improved DDR potential in
\cite[Remark~3.7]{Bonaldi2025} has polynomial degree $r+1$ and uses
the same discrete unknowns. Denote this reconstruction by
\[
  \widehat{\mathsf P}_{r+1,f}^0:
  \underline{\Vcal}_{r,f}^0\longrightarrow\Pcal_{r+1}\Lambda^0(f).
\]
We also introduce a modified scalar family for which this reconstruction
coincides with the degree-$r+1$ $L^2$ projection.
For $r\ge0$ and $f\in\Delta_d(\Mh)$, define
\begin{equation*}
\begin{aligned}
\widetilde{\Vcal}_r^0(f):=
\biggl\{v\in\Vcal_{r+2}^{-,0}(f):\;&
  \int_e
  \bigl(\trc_e\rd v-\Proj{r}{e}{1}\trc_e\rd v\bigr)
  \wedge\mu=0,\\
&\quad\mu\in\kappa\bigl(
  \mathcal H_r\Lambda^\ell(e)\oplus
  \mathcal H_{r+1}\Lambda^\ell(e)\bigr),\\
&\quad e\in\Delta_\ell(f),\quad
  1\le\ell\le d\biggr\}.
\end{aligned}
\end{equation*}
This family contains $\Pcal_{r+1}\Lambda^0(f)$, coincides with it
for $d\le1$, and is closed under tangential traces.
For $v\in\widetilde{\Vcal}_r^0(f)$, let $\underline v_f$ collect
the same moment data \eqref{eq:complete-dofs} as for the complete
scalar family. Using the Koszul splitting
\eqref{eq:koszul-degree-splitting} and the modified moment constraints,
one can verify, by the same integration-by-parts argument as in
Lemma~\ref{lem:local-reconstruction-identities}, that
\begin{equation*}
  \mathsf d_{r,f}^0\underline v_f
  =\Proj{r}{f}{1}\rd v\quad(d\ge1),
  \qquad
  \widehat{\mathsf P}_{r+1,f}^0\underline v_f
  =\Proj{r+1}{f}{0}v.
\end{equation*}
For the corresponding three-dimensional construction and projection
formulas, see
\cite[Section~6.2.2]{BeiraoDassiDiPietroDroniou2022}.
To simplify the analysis, we retain $\Vcal_r^0(f)$ in the complex and
use degree-$r$ $L^2$ projections in every form degree.
\end{remark}

\begin{theorem}[Cochain isomorphisms]
\label{thm:fully-discrete-cochain-isomorphisms}
For every $r\ge0$, the trimmed and complete DOF maps are cochain
isomorphisms onto their respective fully discrete complexes.
With the discrete exterior derivatives defined in
\eqref{eq:global-trimmed-fully-discrete-derivative} and
\eqref{eq:global-complete-fully-discrete-derivative}, they satisfy,
for $0\le k<n$,
\begin{equation}
\label{eq:global-cochain-isomorphisms}
  \underline{\rd}_{r+1,h}^{-,k}\mathcal J_{r+1,h}^{-,k}
  =\mathcal J_{r+1,h}^{-,k+1}\rd,
  \qquad
  \underline{\rd}_{r,h}^{k}\mathcal J_{r,h}^{k}
  =\mathcal J_{r,h}^{k+1}\rd.
\end{equation}
\end{theorem}
\begin{proof}
For the trimmed family, let $\omega_h\in\Vcal_{r+1,h}^{-,k}$ and set
$\underline\omega_f=\mathcal J_{r+1,f}^{-,k}\trc_f\omega_h$.
On $f\in\Delta_{k+1}(\Mh)$, identity
\eqref{eq:local-trimmed-reconstruction-identification} gives
$\mathsf d_{r+1,f}^{-,k}\underline\omega_f=\trc_f\rd\omega_h$.
For $\dim f\ge k+2$, both sides of the cochain identity have component
$(D_{\omega,f},0)$ by
\eqref{eq:global-fully-discrete-isomorphisms},
\eqref{eq:global-trimmed-fully-discrete-derivative}, and $\rd^2=0$.

For the complete family, let $\omega_h\in\Vcal_{r,h}^{k}$ and
$f\in\Delta_d(\Mh)$, with $d\ge k+1$.
\eqref{eq:global-complete-fully-discrete-derivative},
\eqref{eq:local-vem-ddr-reconstruction-identification},
and the projection identities
\eqref{eq:l2-orthogonal-projection}--\eqref{eq:hodge-projection-moments} give
\[
\begin{aligned}
  \bigl(\underline{\rd}_{r,h}^{k}\mathcal J_{r,h}^{k}\omega_h\bigr)_f
  &=\Pi_{r,f}^{-,d-k-1}
    \bigl(\star_f\Proj{r}{f}{k+1}\trc_f\rd\omega_h\bigr)\\
  &=\Pi_{r,f}^{-,d-k-1}\bigl(\star_f\trc_f\rd\omega_h\bigr)
   =\bigl(\mathcal J_{r,h}^{k+1}\rd\omega_h\bigr)_f.
\end{aligned}
\]
Thus both cochain identities hold, and global unisolvence gives bijectivity.
\end{proof}

\subsection{Equivalent formulations and analytical transfer}
\label{subsec:appendix-analytical-transfer}
\label{rem:analytical-fully-discrete-link}
For the complete family, the DDR stabilization of
\cite[Remark~3.14]{Bonaldi2025}, with the weighted data inner product
associated with \eqref{eq:complete-local-data-norm}, gives the local inner product
\begin{equation}
\label{eq:fully-discrete-L2-inner-product}
\begin{aligned}
  [\underline\omega_T,\underline\zeta_T]_{k,T}^{\mathrm{fd}}
  &:=\ip{\mathsf P_{r,T}^k\underline\omega_T}
          {\mathsf P_{r,T}^k\underline\zeta_T}_T\\
  &\quad+\sum_{d=k}^{n}\sum_{f\in\Delta_d(T)}h_T^{n-d}
    \left\langle
      \bigl((I-\mathcal J_{r,T}^k\mathsf P_{r,T}^k)
          \underline\omega_T\bigr)_f,
      \bigl((I-\mathcal J_{r,T}^k\mathsf P_{r,T}^k)
          \underline\zeta_T\bigr)_f
    \right\rangle_f.
\end{aligned}
\end{equation}
Let $\mathcal L_{r,T}^k:=(\mathcal J_{r,T}^k)^{-1}$ denote the local
conforming VEM lifting.
Substituting the reconstruction identity
\eqref{eq:local-vem-ddr-reconstruction-identification} into
\eqref{eq:fully-discrete-L2-inner-product} and using
\eqref{eq:discrete-bilinear}--\eqref{eq:stabilization}, we obtain,
for all $\underline\omega_T,\underline\zeta_T\in\underline{\Vcal}_{r,T}^k$,
\begin{equation}
\label{eq:fully-discrete-inner-product-identification}
  [\underline\omega_T,\underline\zeta_T]_{k,T}^{\mathrm{fd}}
  =[\mathcal L_{r,T}^k\underline\omega_T,\mathcal L_{r,T}^k\underline\zeta_T]_{k,T}.
\end{equation}
Under this identification, the fully discrete inner product inherits
the following two properties from its VEM counterpart.
\begin{itemize}
\item[(i)] \emph{Polynomial consistency.}
For every $\mu\in\Pcal_r\Lambda^k(T)\subseteq\Vcal_r^k(T)$, its fully
discrete interpolant satisfies
$\underline{\mathcal I}_{r,T}^k\mu=\mathcal J_{r,T}^k\mu$, so that
$\mathcal L_{r,T}^k\underline{\mathcal I}_{r,T}^k\mu=\mu$.
Hence, by \eqref{eq:fully-discrete-inner-product-identification},
\eqref{Consistency}, and \eqref{eq:local-vem-ddr-reconstruction-identification},
\[
  [\underline\omega_T,\underline{\mathcal I}_{r,T}^k\mu]_{k,T}^{\mathrm{fd}}
  =[\mathcal L_{r,T}^k\underline\omega_T,\mu]_{k,T}
  =\ip{\mathcal L_{r,T}^k\underline\omega_T}{\mu}_T
  =\ip{\mathsf P_{r,T}^k\underline\omega_T}{\mu}_T.
\]
This recovers the DDR consistency identity \cite[Eq.~(3.16)]{Bonaldi2025}.
\item[(ii)] \emph{Uniform stability.}
By \eqref{eq:fully-discrete-inner-product-identification} and \eqref{Coercivity},
\[
  \alpha_*\norm{\mathcal L_{r,T}^k\underline\omega_T}_{L^2(T)}^2
  \le [\underline\omega_T,\underline\omega_T]_{k,T}^{\mathrm{fd}}
  \le \alpha^*\norm{\mathcal L_{r,T}^k\underline\omega_T}_{L^2(T)}^2.
\]
Thus the norm induced by the local DDR inner product is uniformly
equivalent to the $L^2(T)$ norm of the conforming VEM lifting.
\end{itemize}
The global inner product and the squared global data norm
$\norm{\cdot}_{\mathrm{data},k}^{2}$ are obtained by summing the local inner
products in \eqref{eq:fully-discrete-L2-inner-product} and the squared local
data norms in \eqref{eq:complete-local-data-norm}, respectively, over $T\in\Th$.

Let $\boldsymbol{\mathcal J}_h$ denote the componentwise DOF map in degrees
$(k-1,k,k)$, set
$\underline{\mathbb V}_{r,h}^k
:=\boldsymbol{\mathcal J}_h\mathbb V_{r,h}^k$,
and denote its inverse by $\boldsymbol{\mathcal L}_h$.
The inner-product and cochain identities identify the discrete
harmonic factors. Define $\underline{\mathcal A}_h$ by
\eqref{eq:section4-Ah} using the fully discrete inner products and exterior
derivatives, and set
$\underline\ell_h(\underline v_h)
:=\sum_{T\in\Th}\ip{j_h}{\mathsf P_{r,T}^k\underline v_T}_T$.
Then, for all $\underline{\boldsymbol u}_h,\underline{\boldsymbol v}_h
\in\underline{\mathbb V}_{r,h}^k$,
\begin{equation}
\label{eq:fully-discrete-mixed-equivalence}
\begin{aligned}
  \underline{\mathcal A}_h
      (\underline{\boldsymbol u}_h,\underline{\boldsymbol v}_h)
  &=\mathcal A_h(\boldsymbol{\mathcal L}_h\underline{\boldsymbol u}_h,
                 \boldsymbol{\mathcal L}_h\underline{\boldsymbol v}_h),\\
  \underline\ell_h(\underline v_h)&=\ell_h(\mathcal L_{r,h}^k\underline v_h).
\end{aligned}
\end{equation}
Consequently, $\underline{\boldsymbol u}_h$ solves the fully discrete mixed
problem if and only if its lifting
$\boldsymbol{\mathcal L}_h\underline{\boldsymbol u}_h$ solves \eqref{eq:discrete-mixed}.
Under conditions {\rm (M1)}--{\rm (M2)}, the local stability estimate in
\textup{(ii)} and the cochain identities give
\begin{equation}
\label{eq:fully-discrete-graph-norm-equivalence}
  \norm{\underline{\boldsymbol u}_h}_{\mathrm{fd},V}
  \eqsim\norm{\boldsymbol{\mathcal L}_h\underline{\boldsymbol u}_h}_V,
  \qquad\underline{\boldsymbol u}_h\in\underline{\mathbb V}_{r,h}^k,
\end{equation}
where $\norm{\cdot}_{\mathrm{fd},V}$ is the product graph norm defined as in
\eqref{eq:discrete-product-norm} using the fully discrete inner products.

We now use these identities and norm equivalences to transfer
stability and error estimates between the two formulations.

\paragraph{Poincar\'e inequalities and inf--sup stability.}
To recover the discrete Poincar\'e inequality for the complete virtual complex
(Theorem~\ref{thm:complete-orthogonal-poincare}), we first derive
\eqref{eq:section4-poincare} from the following fully discrete
estimate~\cite[Remark~2 and Theorem~4]{pietro_uniform_2025}.
The auxiliary submesh constructed in
Remark~\ref{rem:auxiliary-simplicial-submesh} ensures that the mesh
assumptions required for this estimate are satisfied.

\begin{lemma}[Fully discrete Poincar\'e inequality]
\label{lem:fully-discrete-poincare}
Under Assumption~\ref{ass:mesh-regularity}, let $r\ge0$
and $0\le k<n$.
For every $\underline \omega_h\in\underline{\Vcal}_{r,h}^k$,
there exists $\underline{\widetilde\omega}_h\in\underline{\Vcal}_{r,h}^k$
such that
\[
  \underline{\rd}_{r,h}^k\underline{\widetilde\omega}_h
  =\underline{\rd}_{r,h}^k\underline \omega_h,
  \qquad
  \norm{\underline{\widetilde\omega}_h}_{\mathrm{data},k}
  \lesssim
  \bigl\|\underline{\rd}_{r,h}^k\underline \omega_h
  \bigr\|_{\mathrm{data},k+1},
\]
with a constant independent of $h$.
\end{lemma}

Given $\underline\omega_h\in\underline{\Vcal}_{r,h}^k$, let
$\underline{\widetilde\omega}_h$ be the form supplied by
Lemma~\ref{lem:fully-discrete-poincare}, and set
\[
  \omega_h=\mathcal L_{r,h}^k\underline\omega_h,\qquad
  \widetilde\omega_h=\mathcal L_{r,h}^k\underline{\widetilde\omega}_h.
\]
By the cochain identity \eqref{eq:global-cochain-isomorphisms}, the two
liftings have the same exterior derivative,
\[
  \rd\widetilde\omega_h
   =\mathcal L_{r,h}^{k+1}
       \underline{\rd}_{r,h}^k\underline{\widetilde\omega}_h
   =\mathcal L_{r,h}^{k+1}
       \underline{\rd}_{r,h}^k\underline \omega_h
   =\rd \omega_h.
\]
Moreover, \eqref{eq:continuous-data-norm-equivalence} and
Lemma~\ref{lem:fully-discrete-poincare} give
\[
  \norm{\widetilde\omega_h}
   \eqsim\norm{\underline{\widetilde\omega}_h}_{\mathrm{data},k}
   \lesssim\bigl\|\underline{\rd}_{r,h}^k\underline \omega_h
            \bigr\|_{\mathrm{data},k+1}
   \eqsim\norm{\rd \omega_h}.
\]
Since $\mathcal L_{r,h}^k$ maps onto $\Vcal_{r,h}^k$, this proves
\eqref{eq:section4-poincare}. By \eqref{eq:global-cochain-isomorphisms}
and \eqref{eq:fully-discrete-inner-product-identification}, the lifting also
maps the orthogonal complement of the fully discrete kernel onto its VEM
counterpart. Hence \cite[Corollary~5]{pietro_uniform_2025} gives
\eqref{eq:section4-orthogonal-poincare}.
For $k=n$, the orthogonal complement is trivial.
Finally, \eqref{eq:fully-discrete-mixed-equivalence} and
\eqref{eq:fully-discrete-graph-norm-equivalence} transfer the inf--sup
bound of Theorem~\ref{thm:section4-inf-sup} to the fully discrete
formulation with a mesh-independent constant.

\paragraph{Adjoint consistency error.}
We next use the conforming lifting and the VEM consistency estimate
to bound the fully discrete adjoint consistency error.
For simplicity, let $1\le k<n$ and assume that $\Omega$
has trivial topology.
Analogously to the VEM residual analysis
\eqref{eq:canonical-residual-expansion}--\eqref{eq:canonical-residual-bound},
the fully discrete residual obtained with the interpolated load
$[\underline{\mathcal I}_{r,h}^k j,\underline v_h]_{k,h}^{\mathrm{fd}}$
in place of $\underline\ell_h(\underline v_h)$ takes the form
\cite[Eq.~(3.27)]{Bonaldi2025}
\[
  \mathcal R_h(\underline\tau_h,\underline v_h)
  =\widetilde{\mathcal E}_h^k(\rd u;\underline v_h)
   -\widetilde{\mathcal E}_h^{k-1}(u;\underline\tau_h),
\]
where the adjoint consistency error is defined by
\begin{equation}
\label{eq:conforming-adjoint-consistency-error}
\begin{aligned}
  \widetilde{\mathcal E}_h^\ell(\omega;\underline \zeta_h)
  :=[\underline{\mathcal I}_{r,h}^\ell\delta\omega,
       \underline \zeta_h]_{\ell,h}^{\mathrm{fd}}-[\underline{\mathcal I}_{r,h}^{\ell+1}\omega,
       \underline{\rd}_{r,h}^\ell\underline \zeta_h]_{\ell+1,h}^{\mathrm{fd}},
  \qquad\underline \zeta_h\in\underline{\Vcal}_{r,h}^\ell.
\end{aligned}
\end{equation}
Here $0\le \ell<n$,
$\underline{\mathcal I}_{r,h}^{\ell}
:=\mathcal J_{r,h}^{\ell}\Ih{\ell}$, and
$\omega\in\overset{\circ}{H}\Lambda^{\ell+1}(\delta;\Omega)$;
all forms are assumed to lie in the required canonical
interpolation domains.
Set $\omega_I:=\Ih{\ell+1}\omega$ and
$(\delta\omega)_I:=\Ih{\ell}\delta\omega$.

To estimate \eqref{eq:conforming-adjoint-consistency-error}, lift the
discrete test form by setting $\zeta_h=\mathcal L_{r,h}^{\ell}\underline\zeta_h$,
as in the DDR adjoint consistency analysis of
\cite[Section~6.2 and Appendix~B]{DiPietroDroniou2023} and
\cite[Theorem~8 and Corollary~9]{DiPietroDroniouPitassi2025}.
Under Assumption~\ref{ass:mesh-regularity}, \eqref{eq:global-cochain-isomorphisms} and
\eqref{eq:continuous-data-norm-equivalence} give
\begin{equation}
\label{eq:inverse-dof-conforming-lifting}
\begin{aligned}
  \rd \zeta_h=\mathcal L_{r,h}^{\ell+1}
       \underline{\rd}_{r,h}^\ell\underline \zeta_h,\quad
  \norm{\zeta_h}\eqsim\norm{\underline \zeta_h}_{\mathrm{data},\ell},
  \quad
  \norm{\rd \zeta_h}\eqsim
  \bigl\|\underline{\rd}_{r,h}^\ell\underline \zeta_h
  \bigr\|_{\mathrm{data},\ell+1}.
\end{aligned}
\end{equation}
Since $\zeta_h$ is conforming and $\omega$ has zero normal trace,
the integration-by-parts formula \eqref{eq:green-formula} gives
$\langle\delta\omega,\zeta_h\rangle=\langle\omega,\rd \zeta_h\rangle$.
Using this cancellation together with
\eqref{eq:fully-discrete-inner-product-identification} and
\eqref{eq:inverse-dof-conforming-lifting}, we obtain
\begin{equation*}
\begin{aligned}
  \widetilde{\mathcal E}_h^\ell(\omega;\underline \zeta_h)
  &=[(\delta\omega)_I,\zeta_h]_{\ell,h}
       -[\omega_I,\rd \zeta_h]_{\ell+1,h}\\
  &=\bigl([(\delta\omega)_I,\zeta_h]_{\ell,h}
       -\langle\delta\omega,\zeta_h\rangle\bigr)
    -\bigl([\omega_I,\rd \zeta_h]_{\ell+1,h}
       -\langle\omega,\rd \zeta_h\rangle\bigr)\\
  &=\mathfrak c_h^\ell(\delta\omega,(\delta\omega)_I;\zeta_h)
       -\mathfrak c_h^{\ell+1}(\omega,\omega_I;\rd \zeta_h).
\end{aligned}
\end{equation*}
Applying the VEM consistency estimate \eqref{eq:section5-consistency-estimate}
to each term and using \eqref{eq:inverse-dof-conforming-lifting} yields
\begin{equation*}
\begin{aligned}
  |\widetilde{\mathcal E}_h^\ell(\omega;\underline \zeta_h)|
  &\lesssim\bigl(\norm{\delta\omega-(\delta\omega)_I}
                   +E_h^\ell(\delta\omega)\bigr)
       \norm{\underline \zeta_h}_{\mathrm{data},\ell}\\
  &\quad+\bigl(\norm{\omega-\omega_I}+E_h^{\ell+1}(\omega)\bigr)
       \bigl\|\underline{\rd}_{r,h}^\ell\underline \zeta_h
             \bigr\|_{\mathrm{data},\ell+1}.
\end{aligned}
\end{equation*}
Together with the canonical VEM interpolation and polynomial approximation
estimates, this bound yields an $O(h^{r+1})$ adjoint consistency estimate
relative to the fully discrete graph norm of the test data, under sufficient
regularity.

\section{Discrete Poincar\'e inequality for the trimmed complex}
\label{sec:appendix-discrete-poincare}
\label{subsec:appendix-discrete-poincare}
We establish a mesh-uniform Poincar\'e inequality for the trimmed
virtual complex. The proof combines a cellular cochain estimate
with the local norm equivalence.

\begin{remark}[Auxiliary simplicial submesh]
\label{rem:auxiliary-simplicial-submesh}
Under {\rm (M1)}--{\rm (M2)}, $\Mh$ admits a conforming,
uniformly shape-regular simplicial submesh $\mathcal S_h$.
Starting from the edges, it is constructed recursively by joining
the centre of $B_f$ to the simplices in the triangulated boundary
of each cell $f$, using the same triangulation on shared subcells.
Each $d$-cell $f$, $d\ge1$, is partitioned into a uniformly bounded
number of $d$-simplices $S$ with $h_S\eqsim h_f$, with constants
depending only on $n$, $\rho$, and $\gamma$.
This submesh satisfies the assumptions in
\cite[Section~2.1.1]{pietro_uniform_2025}.
\end{remark}

\begin{lemma}[Cellular Poincar\'e inequality,
  {\cite[Lemma~6]{pietro_uniform_2025}}]
\label{lem:cellular-poincare}
Under conditions {\rm (M1)} and {\rm (M2)}, let $0\le k<n$ and
$\xi=(\xi_g)_{g\in\Delta_{k+1}(\Mh)}
\in\operatorname{im}\delta_{\Mh}^k$,
where $C^\ell(\Mh)$ denotes the space of real cellular $\ell$-cochains
for a fixed orientation of the mesh cells, and
$\delta_{\Mh}^k:C^k(\Mh)\to C^{k+1}(\Mh)$ is the coboundary
\begin{equation}
\label{eq:cellular-coboundary}
  (\delta_{\Mh}^k c)(g)
  :=\sum_{f\in\Delta_k(g)}\varepsilon_{gf}c(f),
  \qquad g\in\Delta_{k+1}(\Mh),
\end{equation}
with $\varepsilon_{gf}$ the oriented incidence number.
Then there exists a $k$-cochain
$\lambda=(\lambda_f)_{f\in\Delta_k(\Mh)}$
such that $\delta_{\Mh}^k\lambda=\xi$ and
\begin{equation}
\label{eq:cellular-poincare-bound}
  \sum_{T\in\Th}h_T^{n-2k}
  \sum_{f\in\Delta_k(T)}\lambda_f^2
  \lesssim
  \sum_{T\in\Th}h_T^{n-2k-2}
  \sum_{g\in\Delta_{k+1}(T)}\xi_g^2.
\end{equation}
\end{lemma}

\begin{theorem}[Discrete Poincar\'e inequality for the trimmed complex]
\label{thm:appendix-trimmed-poincare}
\label{thm:section4-discrete-poincare}
Under conditions {\rm (M1)} and {\rm (M2)}, let $r\ge0$ and $0\le k<n$.
For every $\omega_h\in\Vcal_{r+1,h}^{-,k}$, there exists
$\widetilde\omega_h\in\Vcal_{r+1,h}^{-,k}$ such that
\begin{equation}
\label{eq:section4-trimmed-primitive-target}
  \rd \widetilde\omega_h=\rd\omega_h,
  \qquad
  \norm{\widetilde\omega_h}\lesssim\norm{\rd\omega_h}.
\end{equation}
\end{theorem}

\begin{proof}
Define $\alpha\in C^k(\Mh)$ and $\xi\in C^{k+1}(\Mh)$ by
\[
  \alpha_f:=\int_f\trc_f\omega_h,
  \quad f\in\Delta_k(\Mh),
  \qquad
  \xi_g:=\int_g\trc_g\rd\omega_h,
  \quad g\in\Delta_{k+1}(\Mh).
\]
Stokes' formula gives
\[
  \xi_g=\sum_{f\in\Delta_k(g)}\varepsilon_{gf}\alpha_f
  =(\delta_{\Mh}^k\alpha)_g.
\]
Lemma~\ref{lem:cellular-poincare} therefore provides
$\lambda\in C^k(\Mh)$ satisfying
$\delta_{\Mh}^k\lambda=\xi$ and
\eqref{eq:cellular-poincare-bound}.

Adapting the moment construction in the proof of
Lemma~\ref{lem:zero-moment-exactness},
we prescribe $\widetilde\omega_h$ as follows.
For $k>0$, set
\[
  \trc_f \widetilde\omega_h=\frac{\lambda_f}{|f|}\star_f1,
  \qquad f\in\Delta_k(\Mh);
\]
for $k=0$, set $\widetilde\omega_h(f)=\lambda_f$ at each vertex.
On higher-dimensional cells, prescribe
\[
  \mathsf D_{d,h}^{k,-}(\widetilde\omega_h)
  =\mathsf D_{d,h}^{k,-}(\omega_h),
  \qquad
  \mathsf K_{d,h}^{k,-}(\widetilde\omega_h)=0,
  \qquad k+1\le d\le n.
\]
Global unisolvence determines a unique
$\widetilde\omega_h\in\Vcal_{r+1,h}^{-,k}$.

For every $g\in\Delta_{k+1}(\Mh)$, Stokes' formula and
$\delta_{\Mh}^k\lambda=\xi$ give
\[
  \int_g\trc_g\rd \widetilde\omega_h
  =\sum_{f\in\Delta_k(g)}\varepsilon_{gf}\lambda_f
  =\xi_g
  =\int_g\trc_g\rd\omega_h.
\]
Thus the two derivatives have the same
$\mathsf C_h^{k+1,-}$ moments.
By the moment correspondence in Remark~\ref{rem:trimmed-moment-shift},
the prescribed $\mathsf D_{d,h}^{k,-}$ moments imply equality of their
$\mathsf K_{d,h}^{k+1,-}$ moments, while their derivative moments vanish
because $\rd^2=0$.
Since both derivatives belong to $\Vcal_{r+1,h}^{-,k+1}$
by \eqref{eq:trimmed-virtual-derham-complex},
global unisolvence yields
\[
  \rd \widetilde\omega_h=\rd\omega_h.
\]

\emph{Uniform $L^2$ bound.}
Fix $T\in\Th$.
For $k>0$ and $f\in\Delta_k(T)$,
{\rm (M1)} and {\rm (M2)} give $|f|\gtrsim h_T^k$.
The prescribed trace therefore satisfies
\begin{equation}
\label{eq:trimmed-primitive-lowest-trace}
  h_T^{n-k}\norm{\trc_f \widetilde\omega_h}_f^2
  =h_T^{n-k}|f|^{-1}\lambda_f^2
  \lesssim h_T^{n-2k}\lambda_f^2.
\end{equation}
For $k=0$, the corresponding estimate is the identity
$h_T^n|\widetilde\omega_h(f)|^2=h_T^n\lambda_f^2$ at each vertex.
Since the higher-dimensional Koszul value moments of
$\widetilde\omega_h$ vanish, \eqref{eq:trimmed-global-upper} gives
\[
\begin{aligned}
  \norm{\widetilde\omega_h}_T^2
  &\lesssim
  \sum_{f\in\Delta_k(T)}
  h_T^{n-k}\norm{\trc_f \widetilde\omega_h}_f^2\\
  &\quad+h_T^2\left[\begin{aligned}
    &\sum_{g\in\Delta_{k+1}(T)}
     h_T^{n-k-1}\norm{\trc_g\rd \widetilde\omega_h}_g^2\\
    &{}+\sum_{d=k+2}^{n}\sum_{f\in\Delta_d(T)}
     h_T^{n-d}
     \norm{\Pi_{r+1,f}^{\kappa,d-k-1}
       (\star\trc_f\rd \widetilde\omega_h)}_f^2
  \end{aligned}\right].
\end{aligned}
\]
Since $\rd^2\widetilde\omega_h=0$, the bracketed expression is equivalent
to $\norm{\rd \widetilde\omega_h}_T^2$ by
\eqref{eq:trimmed-global-upper} at degree $k+1$.
Together with \eqref{eq:trimmed-primitive-lowest-trace}
and $\rd \widetilde\omega_h=\rd\omega_h$, this yields
\begin{equation*}
\begin{aligned}
  \norm{\widetilde\omega_h}_T^2
  &\lesssim
  \sum_{f\in\Delta_k(T)}h_T^{n-k}
  \norm{\trc_f \widetilde\omega_h}_f^2
  +h_T^2\norm{\rd \widetilde\omega_h}_T^2\\
  &\lesssim h_T^{n-2k}
  \sum_{f\in\Delta_k(T)}\lambda_f^2
  +h_T^2\norm{\rd\omega_h}_T^2.
\end{aligned}
\end{equation*}

Summing over $T\in\Th$ and applying
\eqref{eq:cellular-poincare-bound}, we obtain
\[
\begin{aligned}
  \norm{\widetilde\omega_h}^2
  &\lesssim
  \sum_{T\in\Th}h_T^{n-2k-2}
  \sum_{g\in\Delta_{k+1}(T)}\xi_g^2
  +\sum_{T\in\Th}h_T^2\norm{\rd\omega_h}_T^2\\
  &\lesssim
  \sum_{T\in\Th}h_T^{n-k-1}
  \sum_{g\in\Delta_{k+1}(T)}
  \norm{\trc_g\rd\omega_h}_g^2
  +\sum_{T\in\Th}h_T^2\norm{\rd\omega_h}_T^2\\
  &\lesssim\norm{\rd\omega_h}^2.
\end{aligned}
\]
The second inequality follows from Cauchy--Schwarz and
$|g|\lesssim h_T^{k+1}$.
The last uses the reverse bound in
\eqref{eq:trimmed-global-upper} at degree $k+1$ and
$h_T\le\operatorname{diam}(\Omega)$.
\end{proof}

\raggedbottom
\bibliographystyle{amsplain}
\bibliography{references}

\end{document}